\documentclass[opre,nonblindrev]{informs3} 
\OneAndAHalfSpacedXI 

\usepackage{endnotes}
\let\footnote=\endnote
\let\enotesize=\normalsize
\def\notesname{Endnotes}
\def\makeenmark{$^{\theenmark}$}
\def\enoteformat{\rightskip0pt\leftskip0pt\parindent=1.75em
  \leavevmode\llap{\theenmark.\enskip}}

\usepackage{natbib}
 \bibpunct[, ]{(}{)}{,}{a}{}{,}
 \def\bibfont{\small}

\usepackage{booktabs}

\usepackage{multirow}
\usepackage{algorithm}
\usepackage{algpseudocode}
\usepackage{mathtools}
\usepackage{subfigure}
\usepackage{mathtools}
\usepackage{paralist}
\usepackage{enumitem}
\usepackage{bold-extra}
\usepackage[pdfstartview=FitH, hidelinks]{hyperref}
\usepackage[export]{adjustbox}
\usepackage[usenames,dvipsnames]{xcolor}
\usepackage{comment}
\usepackage[normalem]{ulem}

\usepackage[sc]{mathpazo}

\graphicspath{{figures/}}

\renewcommand{\vec}[1]{%
  \boldsymbol{#1}
}

\newcommand{\xx}{\vec{x}}
\newcommand{\yy}{\vec{y}}
\newcommand{\uu}{\vec{u}}
\newcommand{\vv}{\vec{v}}
\newcommand{\zz}{\vec{z}} 
\newcommand{\ww}{\vec{w}}
\renewcommand{\ee}{\vec{e}}

\newcommand{\cc}{\vec{c}}
\newcommand{\ff}{\vec{f}}
\newcommand{\bb}{\vec{b}}

\newcommand{\vzero}{\boldsymbol{{0}}} 
\newcommand{\itt}{\mathit{{(t)}}}

\newcommand{\cut}{\mathcal{H}}
\newcommand{\hp}{\partial}
\newcommand{\epi}{\mathcal{E}}
\newcommand{\calL}{\mathcal{L}}
\newcommand{\calC}{\mathcal{C}}
\newcommand{\calO}{\mathcal{O}}

\newcommand{\prf}[1]{\proof{Proof.}#1\Halmos
\endproof}

\DeclareMathOperator*{\dom}{dom}
\DeclareMathOperator*{\agmin}{argmin}

\TheoremsNumberedThrough     
\ECRepeatTheorems

\EquationsNumberedThrough  
\MANUSCRIPTNO{OPRE-2021-08-503-R3} 

\begin{document}

\RUNAUTHOR{Hosseini and Turner}
\RUNTITLE{Deepest Benders Cuts}
\TITLE{Deepest Cuts for Benders Decomposition}

\ARTICLEAUTHORS{
\AUTHOR{Mojtaba Hosseini$^\text{a}$, John G. Turner$^\text{b}$}
\AFF{$^\text{a}$Tippie College of Business, University of Iowa, IA,}
\AFF{$^\text{b}$The Paul Merage School of Business, University of California, Irvine, CA,\\ \EMAIL{mojtaba-hosseini@uiowa.edu}, \EMAIL{john.turner@uci.edu}}
}

\ABSTRACT{
Since its inception, Benders Decomposition (BD) has been successfully applied to a wide range of large-scale mixed-integer (linear) problems. 
The key element of BD is the derivation of Benders cuts, which are often not unique.
In this paper, we introduce a novel unifying Benders cut selection technique based on a geometric interpretation of cut ``depth'', produce deepest Benders cuts based on $\ell_p$-norms,
and study their properties. Specifically, we show that deepest cuts resolve infeasibility through minimal deviation (in a distance sense) from the incumbent point, are relatively sparse, and may produce optimality cuts even when classical Benders would require a feasibility cut.
Leveraging the duality between separation and projection, we develop a Guided Projections Algorithm for producing deepest cuts while exploiting the combinatorial structure and decomposability of problem instances. 
We then propose a generalization of our Benders separation problem, which not only brings several well-known cut selection strategies under one umbrella, but also, when endowed with a homogeneous function, enjoys several properties of geometric separation problems.
We show that, when the homogeneous function is linear, the separation problem takes the form of the Minimal Infeasible Subsystems (MIS) problem. 
As such, we provide systematic ways of selecting the normalization coefficients of the MIS method, and introduce a Directed Depth-Maximizing Algorithm for deriving these cuts. 
Inspired by the geometric interpretation of distance-based cuts and the repetitive nature of two-stage stochastic programs, we introduce a tailored algorithm to further facilitate deriving these cuts.
Our computational experiments on various benchmark problems illustrate effectiveness of deepest cuts in reducing both computation time and number of Benders iterations, and producing high quality bounds at early iterations.
}


\maketitle

\section{Introduction}
Since \cite{benders1962partitioning} originally proposed a procedure for solving Mixed-Integer Linear Programming (MILP) problems that temporarily fixes some variables to produce one or more much easier-to-solve subproblems at the expense of additional inference and algorithm iterations, Benders Decomposition (BD) has increasingly attracted the attention of researchers in the last five decades. Of note, BD has proven very effective in tackling several classes of challenging MILP problems through both the classical as well as the generalized and logic-based variants of the BD algorithm.

The inherent capacity of BD for exploiting the structural properties of problems with complicating variables has made it one of the most prominent exact algorithms for solving large-scale optimization problems. 
Over the years, BD has grown in its ability to solve a wide range of challenging problems including variants of facility location problems
\citep{magnanti1981accelerating, fischetti2016benders,fischetti2017redesigning},
supply chain and network design problems \citep{keyvanshokooh2016hybrid,alshamsi2018large, fontaine2018benders, pearce2018disaggregated}, 
hub location problems \citep{contreras2011benders, contreras2012exact, maheo2017benders, taherkhani2020benders}, 
scheduling and routing problems \citep{mercier2008theoretical,adulyasak2015benders,bodur2016mixed,bayram2017shelter, perrykkad2022simultaneous},
healthcare operations \citep{cho2014simultaneous,naderi2021increased}, 
machine learning \citep{rahimi2021efficient}, e-commerce \citep{fontaine2023branch}, and variants of stochastic programming problems \citep{santoso2005stochastic, adulyasak2015benders, bodur2016strengthened,khassiba2020two, crainic2021partial,taherkhani2021robust} among several other applications.

BD, at its core, is a relax and ``learn from mistakes'' procedure \citep{hooker2003logic}. In classical BD, this learning mechanism is naturally manifested through Linear Programming (LP) duality and mistakes are ``corrected'' via Benders feasibility and optimality cuts. These cuts are obtained by solving the dual of the subproblem induced by fixing the complicating variables. The learning mechanism, however, need not be restricted to cuts based on LP duality. \cite{geoffrion1972generalized} laid the foundation for extending BD to general nonlinear optimization problems, \cite{hooker2003logic} introduced logic-based BD for tackling problems with logical constraints, and \cite{codato2006combinatorial} tailored this idea to MILP problems involving big-M constraints. 
By treating the separation problem as a feasibility problem and exploiting the duality between separation and projection, we derive what we call \textit{deepest Benders cuts}. With this perspective, the learning component in our Benders procedure can be viewed as cutting off the maximum amount of detected infeasibility, which, interestingly, coincides with resolving infeasibility in this feasibility problem through minimal deviation from the incumbent point. In essence, this is equivalent to learning a ``projection'' point, which we use for facilitating the derivation of these cuts.

Despite its promising structure, a na\"{i}ve implementation of BD may suffer from slow convergence and other computational deficiencies. A wealth of studies have addressed different drawbacks of BD from different angles (see \citealt{rahmaniani2017benders} and references therein for recent advancements and \citealt{bonami2020implementing,maher2021implementing} for efficient implementation guidelines). As with other cutting-plane algorithms, the convergence rate is directly tied to the effectiveness of the generated cuts. Given that there is typically more than one way to generate a Benders cut, an important theoretical and practical question is how to select the most ``effective" cut(s) in each iteration.  This question has spawned a stream of research, which we contribute to.

In their seminal paper, \cite{magnanti1981accelerating} introduced a general-purpose cut selection strategy for selecting a nondominated (or Pareto-optimal) optimality cut among the alternative optimal solutions of the subproblem. More recently, \cite{fischetti2010note} cast the Benders subproblem as a feasibility problem, and proposed an alternative cut selection criterion that approximately identifies a minimal source of infeasibility from the derived feasibility problem. 
\cite{saharidis2010improving} introduced the Maximum Feasibility Subsystem (MFS) cut generation strategy for accelerating BD in problems where most generated cuts are feasibility (as opposed to optimality) cuts. \cite{sherali2013generating} treated cut generation as a multi-objective optimization problem and proposed generating maximal nondominated cuts, which they produced by perturbing the right-hand-side of the primal subproblem.  
Finally, \cite{bodur2016mixed} and \cite{bodur2016strengthened} proposed methods for sharpening Benders cuts using mixed-integer rounding schemes.

We depart from these studies by 
\begin{inparaenum}[(i)]
    \item taking the ``depth'' of the candidate cuts explicitly into account, 
    \item providing a unifying framework for producing deep optimality and feasibility cuts, and 
    \item introducing Benders distance functions that bring several cut selection strategies under one umbrella.
\end{inparaenum}
We begin Section~\ref{sec:classical-bd} with an outline of the classical BD algorithm, which we then contrast to an alternative decomposition scheme in Section~\ref{sec:unifying-decomposition-scheme}. This paves the way for us to formally define what we mean by ``deep'' Benders cuts. In Section~\ref{sec:deepest-benders-cuts}, we introduce a procedure to produce a so-called ``deepest Benders cut'' by taking the Euclidean depth of the candidate cuts as a measure of cut quality. Then we extend the notion of depth using general $\ell_p$-norms in Section~\ref{sec:lp-norm-deepest-cuts} and study the properties of deepest cuts in Section~\ref{sec:deepest-cuts-primal}.

In Section \ref{sec:normalized-distance-functions} we introduce normalized distance functions which replace the $\ell_p$-norm in the denominator of the distance function with a general positive homogeneous normalization function. As a special case, we introduce projective distance functions in Section \ref{sec:projective} and show that they admit a simple characterization of the distance to a candidate hyperplane through the gradient of the normalization function.
In Section \ref{sec:linear_pseudonorms} we introduce distance functions based on linear normalization functions and present several ways of deriving effective normalization coefficients for these linear normalization functions, which connect our method to other cut selection strategies by \cite{fischetti2010note}, \cite{magnanti1981accelerating} and \cite{conforti2019facet}.

Next, we present algorithms for deriving Benders cuts in a combinatorial fashion in Section \ref{sec:method_structured}. In particular, we introduce our Directed Depth-Maximizing Algorithm (DDMA) for linear normalization functions in Section~\ref{sec:directed-depth-maximizing} and our Guided Projections Algorithm (GPA) for $\ell_p$ norms in Section~\ref{sec:guided_projections_algorithm}. 
We dedicate Appendix~\ref{app:separable_subproblem} to separable subproblems, with a special attention to two-stage stochastic programs, and introduce a tailored algorithm for these problems by using the geometric interpretation of distance functions.
In Section~\ref{sec:computational-experiments}, we run computational experiments on various benchmark problems, including deterministic and stochastic variants of facility/network design problems and network interdiction problem, to test the performance of Benders' cuts with different choices of distance functions. Finally, we summarize our conclusions in Section~\ref{sec:conclusions}.

\subsection{Classical Benders Decomposition}\label{sec:classical-bd}

We begin with a brief outline of the classical BD algorithm. Consider the MILP problem
\begin{equation}
    \begin{split}\label{generic-milp}
        [\mbox{OP}]\;\min\quad & {\cc}^{\top}{\xx}+{\ff}^{\top}{\yy}\\
    \text{s.t.} \quad & A{\xx}+B{\yy}\ge \bb\\
    & {\xx} \ge {\vzero}, {\yy}\in Y,
    \end{split}
\end{equation}
where ${\ff}\in \mathbb{R}^n$, ${\cc}\in \mathbb{R}^{n^\prime}$, ${\bb}\in \mathbb{R}^m$, matrices $A$ and $B$ are conformable, and
$Y\subset \mathbb{Z}^n$ is the domain of the ${\yy}$-variables. In what follows, we reserve $i$ and $j$ for indexing the rows and columns of $B$, respectively.
For the sake of generality, we do not make any specific assumptions about the structure of the problem, except that it is a bounded MILP. We shall discuss structured MILPs in Section~\ref{sec:method_structured}.

Let us define $Q({\yy})= {\ff}^{\top}{\yy}+\tilde{Q}({\yy})$ where $\tilde{Q}({\yy})$ accounts for the contribution of the ${\xx}$-variables to the objective function for given $\yy$ which is obtained by solving the primal subproblem (PSP)
\begin{align}
    [\mbox{PSP}]\;\tilde Q({\yy})= \min\;& \left\{{\cc}^{\top}{\xx}: A{\xx}\ge \bb- B{\yy}, {\xx} \ge {\vzero}\right\}.\label{generic-milp-projection-psp}
\end{align}
Since original problem (OP) is bounded, PSP is also bounded for any ${\yy}\in Y$. 
Let $\dom(Q)$ be the set of ${\yy}$ values that induce a feasible PSP. 
We may restate OP in the space of ${\yy}$-variables as 
\begin{equation}
    \begin{split}\label{generic-milp-projection}
        \min\left\{ Q({\yy}): {\yy}\in Y\cap \dom(Q)\right\}.
    \end{split}
\end{equation}
In the classical BD algorithm, problem \eqref{generic-milp-projection} is first reformulated in epigraph form as
\begin{align}
    \min\; \{\eta : ({\yy},\eta)\in\epi, {\yy}\in Y\},\label{generic-milp-projection-reform}
\end{align}
where $\epi$ is the epigraph of $Q$ defined as
\begin{align*}
    \epi= \left\{({\yy},\eta)\in \mathbb{R}^{n+1}: \eta\ge Q({\yy}), {\yy}\in \dom(Q)\right\}.
\end{align*}
Then, a relaxation of \eqref{generic-milp-projection-reform} is successively tightened by progressively outer-approximating $\epi$ with supporting hyperplanes obtained by evaluating, at given ${\yy}$ values, the dual of \eqref{generic-milp-projection-psp} formulated as 
\begin{align}
    [\mbox{DSP}]\;\tilde Q({\yy})=\max\;\{{\uu}^{\top}({\bb}-B{\yy}): {\uu}^{\top} A\le {\cc}^{\top} ,{\uu} \ge {\vzero}\},\label{generic-milp-projection-dsp}
\end{align}
which is known as the dual subproblem (DSP). From this dual formulation, we can observe that $\tilde Q({\yy})$ is a piece-wise linear convex function of ${\yy}$. Thus, $Q({\yy})={\ff}^{\top}{\yy}+\tilde{Q}({\yy})$ is a piece-wise linear convex function and $\epi$ is a closed convex set.
Let $\mathcal{U}$ denote the polyhedron defining the set of feasible solutions of DSP, with $\mathcal{U}^*$ and $\mathcal{V}^*$ as its extreme points and rays, respectively. For ${\yy}\in \dom(Q)$, the DSP induced by ${\yy}$ is bounded and its optimal value is attained at one of the extreme points of $\mathcal{U}$. Additionally, since $\tilde Q({\yy})$ is the optimal value of DSP, it follows from weak duality that
\begin{align*}
    \eta \ge Q({\yy})={\ff}^{\top}{\yy}+\tilde{Q}({\yy})\ge {\ff}^{\top}{\yy}+\hat {\uu}^{\top}({\bb}-B{\yy})\quad\quad \forall \hat{\uu}\in \mathcal{U}.
\end{align*}
On the other hand, by Farkas lemma, the values of ${\yy}$ that induce an infeasible PSP (i.e., an unbounded DSP) are the ones for which $\hat {\vv}^{\top}({\bb}-B{\yy})>0$ for some (extreme) ray $\hat{\vv}$ of $\mathcal{U}$; hence, $\dom(Q)=\{{\yy}: 0\ge \hat {\vv}^{\top}({\bb}-B{\yy}) \;\; \forall \hat{\vv}\in \mathcal{V}^*\}$. Putting these pieces together, we can rewrite \eqref{generic-milp-projection-reform} as 
\begin{align}
    [\mbox{CMP}]\;\min\quad& \eta\label{classical-mp-obj}\\
    \text{s.t.}\quad& \eta\ge {\ff}^{\top}{\yy}+ \hat {\uu}^{\top}({\bb}-B{\yy}) & \forall \hat{\uu}\in \mathcal{U}^*\label{classical-mp-optimality-cuts}\\
    & 0\ge \hat {\vv}^{\top}({\bb}-B{\yy}) & \forall \hat{\vv}\in \mathcal{V}^*\label{classical-mp-feasibility-cuts}\\
    & \eta\in\mathbb{R}, {\yy}\in Y,\label{classical-mp-domain}
\end{align}
which we refer to as the classical Benders master problem (CMP). Constraint sets \eqref{classical-mp-optimality-cuts} and \eqref{classical-mp-feasibility-cuts} are known as the \textit{Benders optimality} and \textit{feasibility} cuts, respectively. 
The classical BD algorithm solves CMP by initially relaxing these constraints, and at each iteration posts one or more cuts of the form \eqref{classical-mp-optimality-cuts} or \eqref{classical-mp-feasibility-cuts} to this relaxation of CMP until the optimality gap is sufficiently closed. 
Thus, BD can be viewed as an outer-approximation (OA) algorithm applied to the integer convex program \eqref{generic-milp-projection}, with the difference that, OA algorithms typically address mixed-integer smooth convex programs \citep{duran1986outer, fletcher1994solving,bonami2008algorithmic}, whereas $Q(\yy)$ is non-smooth and \eqref{generic-milp-projection} is purely integer \cite[see, e.g., ][for a more detailed discussion]{belotti2013mixed}.

\subsection{A Unifying Decomposition Scheme}\label{sec:unifying-decomposition-scheme}
In classical BD, ${\yy}$ is the only piece of information passed from the master problem to the subproblems, and $\eta$ is merely used to obtain a lower bound on OP. To incorporate $\eta$ when producing the Benders cuts, we first reformulate the original problem \eqref{generic-milp} in epigraph form as
\begin{equation}
    \begin{split}\label{generic-milp-2}
        \min\quad & \eta\\
    \text{s.t.} \quad & \eta \ge {\cc}^{\top}{\xx}+{\ff}^{\top}{\yy}\\
    & A{\xx}+B{\yy}\ge \bb\\
    & {\xx} \ge {\vzero}, {\yy}\in Y,
    \end{split}
\end{equation}
and then apply BD to \eqref{generic-milp-2} by taking $({\yy},\eta)$ as the master problem variables. This alternative decomposition scheme, which was first introduced in \cite{fischetti2010note}, allows us to treat the Benders optimality and feasibility cuts in a unified framework.
Taking this viewpoint, the primal subproblem induced by trial solution $(\hat{\yy},\hat \eta)$ is the following \textit{feasibility subproblem} (FSP)
\begin{equation}
    \begin{split}\label{generic-milp-fsp}
        [\mbox{FSP}]\;\min\quad & 0\\
    \text{s.t.} \quad & -{\cc}^{\top}{\xx} \ge {\ff}^{\top}\hat{\yy}-\hat \eta\\
    & A{\xx}\ge {\bb}-B\hat{\yy}\\
    & {\xx} \ge {\vzero}.
    \end{split}
\end{equation}
While the constraint $-{\cc}^{\top}{\xx} \ge {\ff}^{\top}\hat{\yy}-\hat \eta$ is sometimes prone to numerical issues \citep[cf.][]{bonami2020implementing}, we can alleviate this with a suitable scaling of the cost vectors $\ff$ and $\cc$ (see Appendix \ref{app:coefficient-scaling}).

Assigning the dual variable $\pi_0$ to the first constraint and the dual vector ${\BFpi}$ to the second set of constraints, a Farkas certificate for infeasibility of FSP can be produced using
\begin{align}
    [\mbox{CGSP}]\;\max_{(\BFpi,\pi_0)\in \Pi}\quad & {\BFpi}^{\top}({\bb}-B\hat {\yy})+\pi_0({\ff}^{\top}\hat{\yy}-\hat\eta),\label{generic-milp-cgsp}
\end{align}
which we refer to as the \textit{certificate generating subproblem (CGSP)}, where
$$\Pi=\{({\BFpi},\pi_0): {\BFpi}^{\top} A\le \pi_0 {\cc}^{\top},\; {\BFpi} \ge {\vzero}, \pi_0 \ge 0\}$$ 
is the cone of feasible solutions (rays). 
Given $(\hat{\BFpi},\hat{\pi}_0)\in\Pi$, we denote by
\begin{align*}
    \cut(\hat{\BFpi},\hat{\pi}_0)=\{({\yy},\eta): \hat{\BFpi}^{\top}({\bb}-B{\yy})+\hat{\pi}_0({\ff}^{\top}{\yy}-\eta)\le 0\},\\
    \hp(\hat{\BFpi},\hat{\pi}_0)=\{({\yy},\eta): \hat{\BFpi}^{\top}({\bb}-B{\yy})+\hat{\pi}_0({\ff}^{\top}{\yy}-\eta)=0\},
\end{align*}
the half-space and hyperplane defined by $(\hat{\BFpi},\hat{\pi}_0)$, respectively.
If FSP is feasible, the optimal value of both FSP and CGSP is zero. Otherwise, CGSP is unbounded and a ray $(\hat{\BFpi},\hat{\pi}_0)$ exists such that $\hat{\BFpi}^{\top}({\bb}-B\hat {\yy})+\hat{\pi}_0({\ff}^{\top}\hat{\yy}-\hat\eta)>0$; hence, the infeasible solution $(\hat{\yy},\hat \eta)$ must violate
$({\yy},\eta)\in \cut(\hat{\BFpi},\hat{\pi}_0)$.
Consequently, OP can be restated as the following modified master problem (MP):
\begin{align}
    [\mbox{MP}]\;\min\quad & \eta\label{mp-obj}\\
    \text{s.t.} \quad & ({\yy},\eta)\in \cut(\hat{\BFpi},\hat{\pi}_0) & \forall (\hat{\BFpi},\hat{\pi}_0)\in \Pi\label{mp-cut}\\
    & \eta\in\mathbb{R}, {\yy}\in Y.\label{mp-domain}
\end{align}
With this representation of BD, at iteration $t$, we produce a candidate point $({\yy}^{\itt},\eta^{\itt})$ by solving a relaxation of MP, and test its feasibility using CGSP. If the test proves $({\yy}^{\itt},\eta^{\itt})$ infeasible, we generate a certificate $(\hat{\BFpi},\hat{\pi}_0)$ and add a cut of the form \eqref{mp-cut} to the relaxed MP to avoid reproducing the infeasible $({\yy}^{\itt},\eta^{\itt})$; otherwise, we conclude that $({\yy}^{\itt},\eta^{\itt})$ is an optimal solution for MP.
See Algorithm~\ref{pseudo-code-bd} for details.
Note that cuts of the form \eqref{mp-cut} represent both Benders optimality and feasibility cuts; when $\hat{\pi}_0>0$, the cut corresponds to a classical Benders optimality cut, while $\hat{\pi}_0=0$ corresponds to a classical Benders feasibility cut.

\begin{algorithm}[h]
	\caption{Overview of Benders Decomposition algorithm}
	\label{pseudo-code-bd}
	\begin{algorithmic}[1]
		\State $t \leftarrow 1$, $\hat\Pi_t \leftarrow \emptyset$.
		\State Solve MP with $\hat \Pi_t$ in place of $\Pi$ and obtain master solution $({\yy}^{\itt},\eta^{\itt})$.
		\State Find a certificate $(\hat{\BFpi},\hat{\pi}_0)$ for infeasibility of $({\yy}^{\itt},\eta^{\itt})$ using CGSP \eqref{generic-milp-cgsp}.
		\If{certificate $(\hat{\BFpi},\hat{\pi}_0)$ exists}
		    \State Set $\hat \Pi_{t+1}\leftarrow \hat \Pi_t\cup\{(\hat{\BFpi},\hat{\pi}_0)\}$, $t\leftarrow t+1$ and go to step 2.
		\Else
		    \State Stop. $({\yy}^{\itt},\eta^{\itt})$ is an optimal solution for MP.
		\EndIf
	\end{algorithmic}
\end{algorithm}

At step 3 of Algorithm~\ref{pseudo-code-bd}, CGSP provides a logical answer to whether the current master problem solution $({\yy}^{\itt},\eta^{\itt})$ is feasible (and hence optimal) for MP. But, not every logical answer is equally useful. 
In other words, to prove suboptimality of $({\yy}^{\itt},\eta^{\itt})$, CGSP produces a certificate ${(\hat{\BFpi},\hat{\pi}_0)\in\Pi}$ such that $\hat{\BFpi}^{\top}({\bb}-B {\yy}^{\itt})+\hat{\pi}_0({\ff}^{\top} {\yy}^{\itt}-\eta^{\itt})>0$, without providing further information about how ``far'' $({\yy}^{\itt},\eta^{\itt})$ is from being optimal.
Moreover, not only do we want to discard the trial solution $({\yy}^{\itt},\eta^{\itt})$, but we also want to rule out as many other sub-optimal solutions $({\yy},\eta)$ as possible.
Hence, we may phrase the key question of the BD algorithm as: \textit{How should we select a certificate $(\hat{\BFpi},\hat{\pi}_0)\in\Pi$ that conveys additional information about the sub-optimality of $({\yy}^{\itt},\eta^{\itt})$, so that we may exploit this information to speed up the convergence of the BD algorithm?} Our order of business in this article is to address this question by introducing selection strategies that exploit the properties of promising cuts in a computationally tractable manner.

\section{Deepest Benders Cuts}\label{sec:deepest-benders-cuts}
At each iteration of the BD Algorithm~\ref{pseudo-code-bd}, we wish to separate (if possible) the incumbent point $({\yy}^{\itt},\eta^{\itt})$ from the epigraph $\epi$. Note that we may equivalently define $\epi$ as
$$\epi=\{({\yy},\eta): ({\yy},\eta)\in \cut(\hat{\BFpi},\hat{\pi}_0)\quad \forall (\hat{\BFpi},\hat{\pi}_0)\in \Pi\}.$$ 
In cutting-plane theory, the separation problem produces a hyperplane (or a cut) that lies between a given point and a closed convex set. In our application, we want to separate the incumbent point  $({\yy}^{\itt},\eta^{\itt})$ from the closed convex set $\epi$ using a hyperplane $\hp({\BFpi},\pi_0)$ for some $(\BFpi,\pi_0)\in\Pi$. Note that infinitely many such hyperplanes may exist, thus one needs a selection criterion for producing the cut that ``best'' separates $({\yy}^{\itt},\eta^{\itt})$ from $\epi$. While there is no universal definition of ``best'' cut, a ``good'' cut should satisfy some natural requirements. First, it should be a supporting hyperplane for $\epi$ in the sense that it should touch $\epi$ at some point. We further postulate that the cut must be \textit{deep}, in the sense that it is as far from the given point $({\yy}^{\itt},\eta^{\itt})$ as possible.
Finally, it is desirable for the cut to expose a facet of $\epi$, which is a stronger property than supporting $\epi$. Note that since facets of $\epi$ are not necessarily facets of the convex hull of $\epi$ when $\yy$ is restricted to integer values (i.e., $\text{conv}(\epi\cap \{(\eta,\yy): \yy\in \mathbb{Z}^n\}$), facet-definingness in this context should not be confused with the facet-definingness property in the integer programming sense. Given that any Benders cut, regardless of how it is selected, can at best expose a facet of $\epi$, facet-definingness for Benders cuts in the integer programming sense is orthogonal to cut selection, and further polishing, such as the mixed-integer rounding scheme of \cite{bodur2016mixed}, is necessary to recover this property.

We begin in Section \ref{sec:euclidean-deepest-cuts} with Euclidean distance as our measure of cut depth, then generalize to distances induced by $\ell_p$-norms in Section \ref{sec:lp-norm-deepest-cuts}. In Section \ref{sec:deepest-cuts-primal} we present an alternative primal perspective of deepest cut generation, and derive some important properties.
Specifically, we show deepest cuts not only support $\epi$, but also \begin{inparaenum}[(i)]
    \item minimally resolve infeasibility in the system FSP, 
    \item amount to optimality cuts, and 
    \item are relatively flat, thus help close the gap quickly.
\end{inparaenum}

\subsection{Euclidean Deepest Cuts}
\label{sec:euclidean-deepest-cuts}
As our measure of cut depth, we start with the Euclidean distance from the point $(\hat{\yy},\hat{\eta})$ to the hyperplane $\hp({\BFpi},\pi_0)$. Euclidean norm is the standard norm used in convex analysis, and measuring depth using this norm is also common practice in cutting-plane theory. For example, to produce deep facet-defining cuts for solving mixed-integer programs, \cite{balas1993lift} and \cite{cadoux2010computing} use the Euclidean distance between the optimal vertex of the current relaxation and candidate separating hyperplanes; in a similar spirit, we also call the cuts we generate \textit{deepest Benders cuts}.  

Given that the Euclidean distance from the point $\hat{\zz}$ to the hyperplane ${\BFalpha^{\top} {\zz}+\beta=0}$ is $\frac{|\BFalpha^{\top}\hat{\zz}+\beta|}{\|\BFalpha\|_2}$, the Euclidean distance between $(\hat{\yy},\hat{\eta})$ and the hyperplane $\hp({\BFpi},\pi_0)$, denoted $d(\hat{\yy},\hat{\eta}|{\BFpi},\pi_0)$, is
\begin{align}
    d(\hat{\yy},\hat{\eta}|{\BFpi},\pi_0)=\frac{|{\BFpi}^{\top} ({\bb}-B\hat{\yy})+\pi_0({\ff}^{\top}\hat{\yy}-\hat{\eta})|}{\|({\BFpi}^{\top} B-\pi_0 {\ff}^{\top},\pi_0)\|_2}=\frac{{\BFpi}^{\top} ({\bb}-B\hat{\yy})+\pi_0({\ff}^{\top}\hat{\yy}-\hat{\eta})}{\|({\BFpi}^{\top} B-\pi_0 {\ff}^{\top},\pi_0)\|_2},\label{distance-euclidean}
\end{align}
where the last equality holds because a sub-optimal $(\hat{\yy},\hat{\eta})$ must violate the constraint ${\BFpi}^{\top}({\bb}-B{\yy})+\pi_0({\ff}^{\top}{\yy}-\eta)\le 0$, hence ${\BFpi}^{\top} ({\bb}-B\hat{\yy})+\pi_0({\ff}^{\top}\hat{\yy}-\hat{\eta})\ge 0$. To produce a \textit{deepest} cut, we choose $({\BFpi},\pi_0)\in\Pi$ which maximizes \eqref{distance-euclidean} by solving the \textit{separation subproblem (SSP)}:
$$[\mbox{SSP}]\;d^*(\hat{\yy},\hat{\eta})=\max_{({\BFpi},\pi_0)\in\Pi} \quad d(\hat{\yy},\hat{\eta}|{\BFpi},\pi_0).$$
As we will show in Section~\ref{sec:deepest-cuts-primal}, maximizing the distance of a separating hyperplane from the point $(\hat{\yy},\hat{\eta})$ coincides with  finding the distance of $(\hat{\yy},\hat{\eta})$ from the epigraph $\epi$; thus we call $d^*(\hat{\yy},\hat{\eta})$ the (Euclidean) distance of $(\hat{\yy},\hat{\eta})$ from the epigraph $\epi$. At iteration $t$ of the BD algorithm, if $d^*({\yy}^{\itt},\eta^{\itt})>0$, then we can separate $({\yy}^{\itt},\eta^{\itt})$ from $\epi$, otherwise $({\yy}^{\itt},\eta^{\itt})\in \epi$. Consequently, if $d^*({\yy}^{\itt},\eta^{\itt})=0$, then $({\yy}^{\itt},\eta^{\itt})$ is an optimal solution for MP.

\begin{figure}[t]
	\centering
		\includegraphics[clip,width=0.3\textwidth]{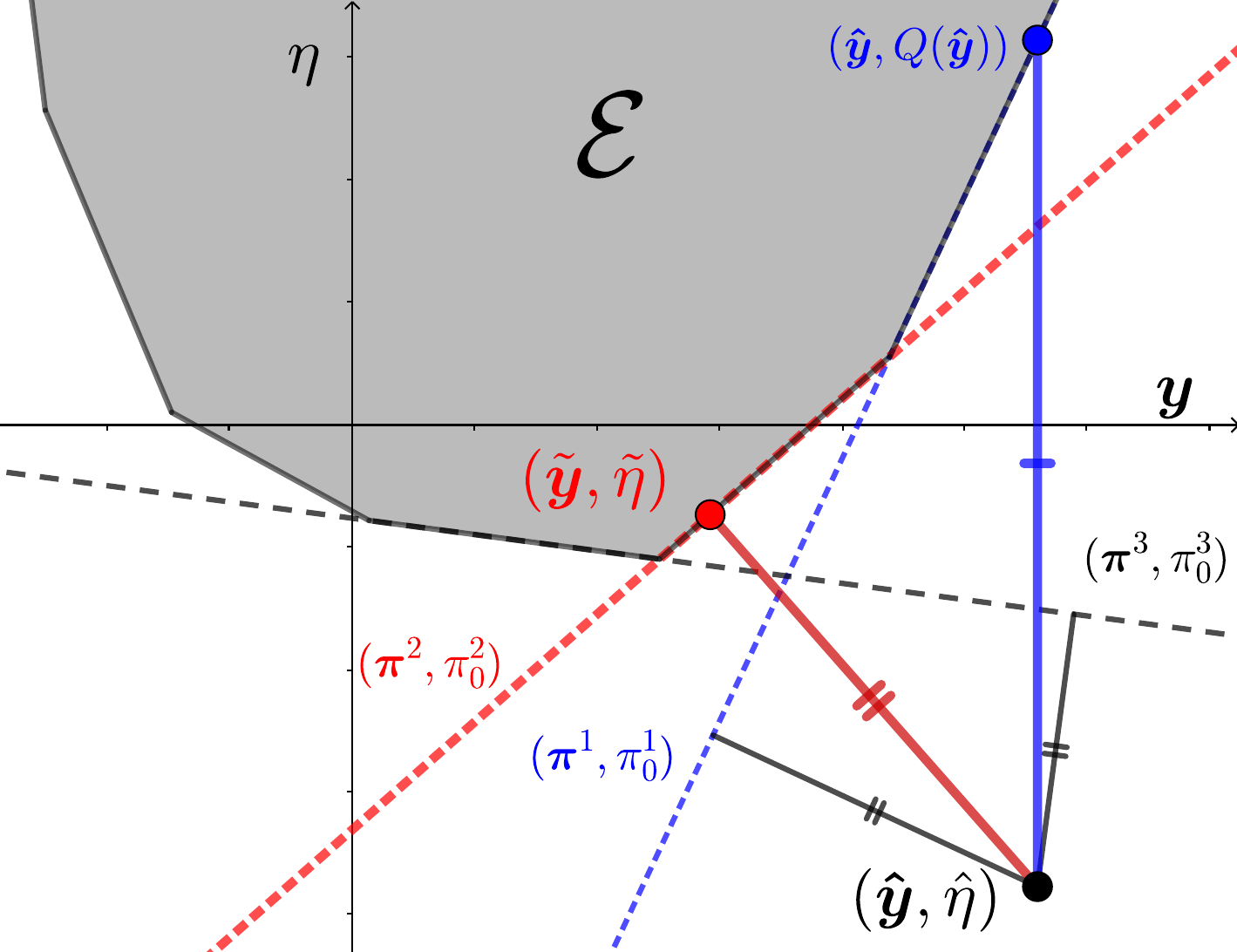}
	\caption{Deepest (red) versus classical (blue) Benders cut selection. Double-crossed segments correspond to Euclidean distances from $(\hat{\yy},\hat{\eta})$ to candidate cuts.}
	\label{fig:cut-selection-deep-classical}
\end{figure}

Figure~\ref{fig:cut-selection-deep-classical} illustrates finding the hyperplane that has the maximum Euclidean distance from the master problem's solution $(\hat{\yy},\hat\eta)$. For demonstration purposes, we assume $y$ is a continuous one-dimensional variable in this toy example. 
The blue single-crossed line indicates the direction that one would take to find the classical Benders cut (see Section \ref{sec:projective}).
As illustrated in Figure~\ref{fig:cut-selection-deep-classical}, the hyperplane produced by classical DSP supports $\epi$ at $(\hat{\yy},Q(\hat{\yy}))$. While deepest cuts also support $\epi$ (see Proposition \ref{prop:deepest-cuts-support}), they do not necessarily do so at $(\hat{\yy},Q(\hat{\yy}))$. More generally, we denote their point of support as $(\tilde{\yy},\tilde{\eta})$, which we call the \textit{projection of} $(\hat{\yy},\hat\eta)$ \textit{onto} $\epi$. 
Dual solutions $({\BFpi}^1,{\pi}^1_0)$, $({\BFpi}^2,{\pi}^2_0)$ and $({\BFpi}^3,{\pi}^3_0)$ (and their convex combinations) are the candidate solutions evaluated based on the Euclidean distance of their associated hyperplanes to the point $(\hat{\yy},\hat\eta)$, and $({\BFpi}^2,{\pi}^2_0)$ is selected as the deepest cut.  It is worth pointing out that in the literature, the question of which Benders cut to select usually arises when the classical DSP admits alternative optimal solutions \citep{magnanti1981accelerating}. However, even when the classical DSP admits a unique optimal solution (as in the given example), the deepest cut may not coincide with the classical Benders cut.

\begin{remark}
    It is possible to define $\eta$ in such a way that it only upper-bounds the $\cc^{\top}\xx$ component of the objective function, but having $\eta$ upper bound $\cc^{\top}\xx+\ff^{\top}\yy$ allows us to take the contribution of $\ff^{\top}\yy$ to the produced cut into account while assessing the depth of the candidate cuts. 
\end{remark}

\subsection{\texorpdfstring{$\ell_p$}{Lp}-norm Deepest Cuts}\label{sec:lp-norm-deepest-cuts}
We now generalize our notion of distance to that induced by any $\ell_p$-norm. Our derivation begins with the observation that the denominator in \eqref{distance-euclidean} is the $\ell_2$-norm of the vector of coefficients $({\BFpi}^{\top} B-\pi_0 {\ff}^{\top},\pi_0)$. If we replace this norm with a general $\ell_p$-norm for $p\ge 1$, and define $d_{\ell_p}$ as
\begin{align}
    d_{\ell_p}(\hat{\yy},\hat{\eta}|{\BFpi},\pi_0)=\frac{{\BFpi}^{\top} ({\bb}-B\hat{\yy})+\pi_0({\ff}^{\top}\hat{\yy}-\hat{\eta})}{\|({\BFpi}^{\top} B-\pi_0 {\ff}^{\top},\pi_0)\|_p},\label{distance-Lp}
\end{align}
then an \textit{$\ell_p$-deepest cut} can be produced by solving the following separation problem
\begin{align}
    [\mbox{SSP}]\;d_{\ell_p}^*(\hat{\yy},\hat{\eta})=\max_{({\BFpi},\pi_0)\in\Pi} \quad&d_{\ell_p}(\hat{\yy},\hat{\eta}|{\BFpi},\pi_0).\label{separation-subproblem-Lp}
\end{align}
As Proposition~\ref{prop:Lq-distance} shows below, $d_{\ell_p}$ measures the distance between $(\hat{\yy},\hat{\eta})$ and the hyperplane $\hp({\BFpi},\pi_0)$ with respect to the dual norm $\ell_q$, where $\frac{1}{p}+\frac{1}{q}=1$. 
The proof of this result, as with all other proofs of propositions and theorems in this paper, is provided in Appendix~\ref{app:proofs}. 

\begin{proposition}\label{prop:Lq-distance}
Given $q\ge 1$ and $p\ge 1$ such that $\ell_p$ is the dual norm of $\ell_q$ (i.e., $\frac{1}{p}+\frac{1}{q}=1$), the $\ell_q$-distance from the point $\hat{\zz}\in\mathbb{R}^{n+1}$ to hyperplane ${\BFalpha}^{\top}{\zz}+\beta=0$ is
\begin{align*}
    \min\limits_{\zz: {\BFalpha}^{\top}{\zz}+\beta=0}\|\zz-\hat{\zz}\|_q=\frac{|{\BFalpha}^{\top}\hat{\zz}+\beta|}{\|{\BFalpha}\|_p}.
\end{align*} 
\end{proposition}
\begin{remark}
    As given in the proof of Proposition~\ref{prop:Lq-distance}, we may extend the definition of deepest cuts by replacing the denominator in \eqref{distance-Lp} with general norms (e.g., a composition of $\ell_p$-norms with different $p$ for different subsets of the components of $({\BFpi}^{\top} B-\pi_0 {\ff}^{\top},\pi_0)$). However, for clarity and simplicity of exposition, we restrict consideration in the remainder of this paper to standard $\ell_p$-norms.
\end{remark}

Some choices of $p$ for $\ell_p$-deepest cuts merit special attention. In particular, for $p=1$ and ${p=\infty}$, $d_{\ell_p}$ defined by \eqref{distance-Lp} measures the $\ell_{\infty}$ and $\ell_1$ distance of $(\hat{\yy},\hat{\eta})$ from the hyperplane $\hp({\BFpi},\pi_0)$, respectively. As given in Appendix~\ref{app:reformulations-lp}, these norms are in general computationally favorable over the $\ell_2$-norm since they result in linear separation subproblems.

As well, note that $\pi_0$ is the coefficient of $\eta$ and ${\BFpi}^{\top}B-\pi_0 {\ff}^{\top}$ is the coefficient of ${\yy}$ in the cut ${\BFpi}^{\top}{\bb}\le ({\BFpi}^{\top}B-\pi_0 {\ff}^{\top}){\yy}+\pi_0\eta$. Therefore, deepest cuts effectively cut off the point $(\hat{\yy},\hat{\eta})$ while minimizing the coefficients of the variables in the produced constraint. In particular, when the $\ell_1$-norm is employed, producing deepest cuts mimics the idea of producing maximally nondominated Benders cuts introduced by \cite{sherali2013generating}, where the cut is maximally nondominated in the sense typically used in the cutting-plane theory for integer programs.

\subsection{A Primal Projection Perspective of the Separation Problem} \label{sec:deepest-cuts-primal}
We now provide another view of deepest cuts, which will be important for analyzing their properties and paves the way for devising algorithms to produce them efficiently. By strong duality, we establish a duality between separation and projection as stated in Theorem~\ref{thrm:deepest-cuts-primal} below.

\begin{figure*}[t]
    \centering
    \subfigure[$\ell_1$-deepest cut ($p=1$, $q=\infty$)]{
 		\includegraphics[clip,width=0.32\textwidth]{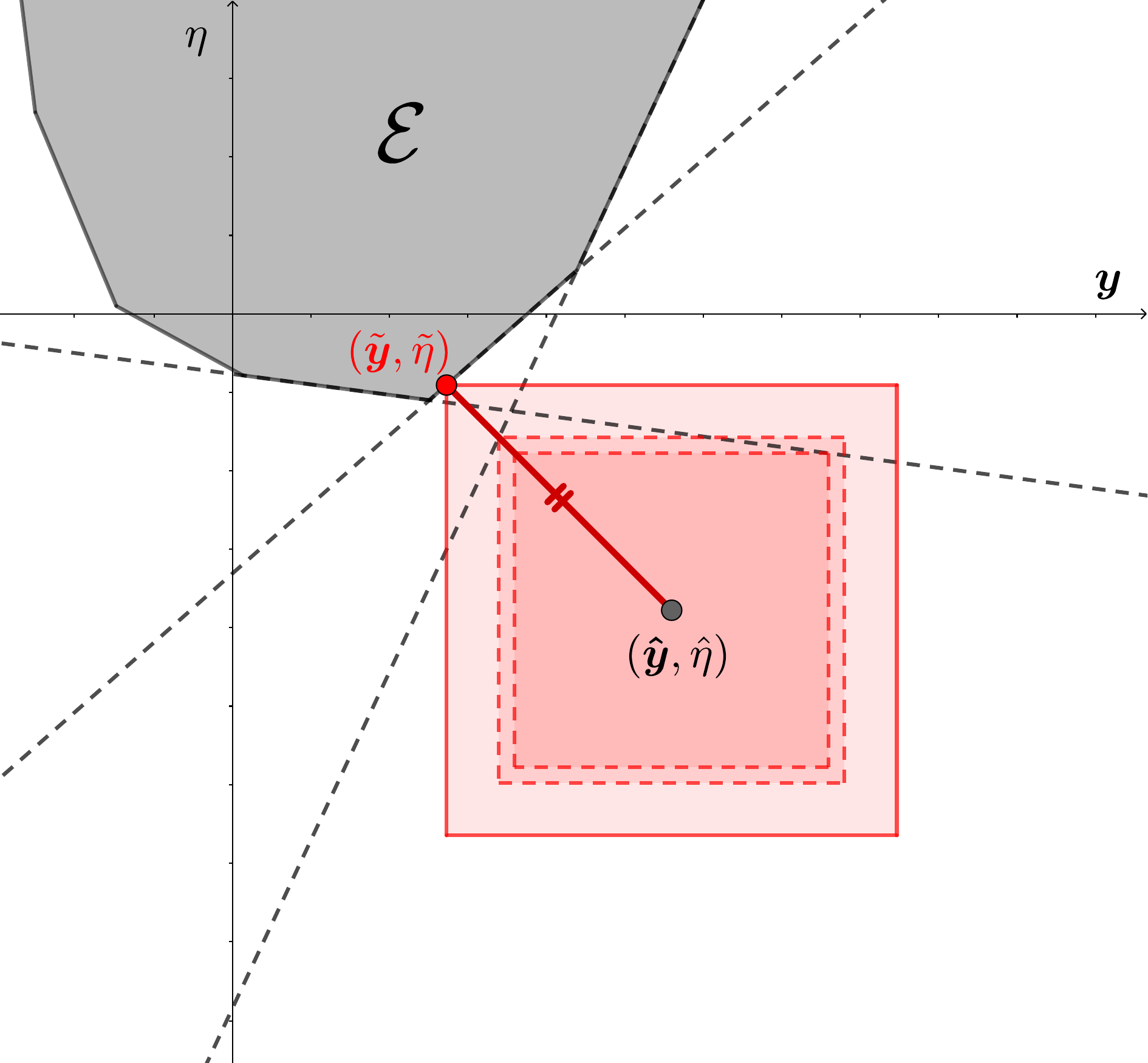}\label{fig:deep_cuts_primal_norm-l1}}
    \subfigure[$\ell_2$-deepest cut ($p=2$, $q=2$)]{
 		\includegraphics[clip,width=0.32\textwidth]{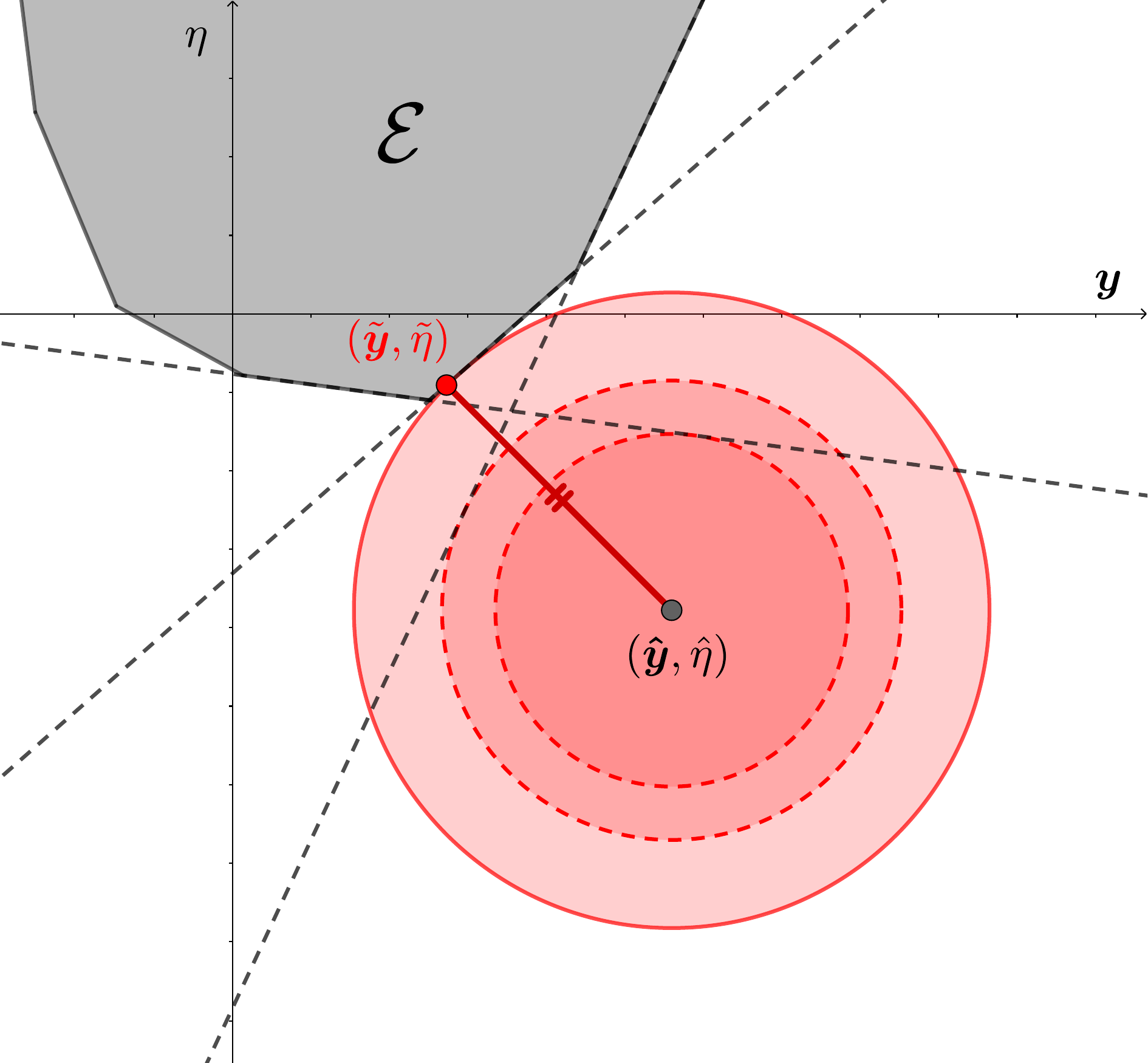}\label{fig:deep_cuts_primal_norm-l2}}
 	\subfigure[$\ell_{\infty}$-deepest cut ($p=\infty$, $q=1$)]{
 	    \includegraphics[clip,width=0.32\textwidth]{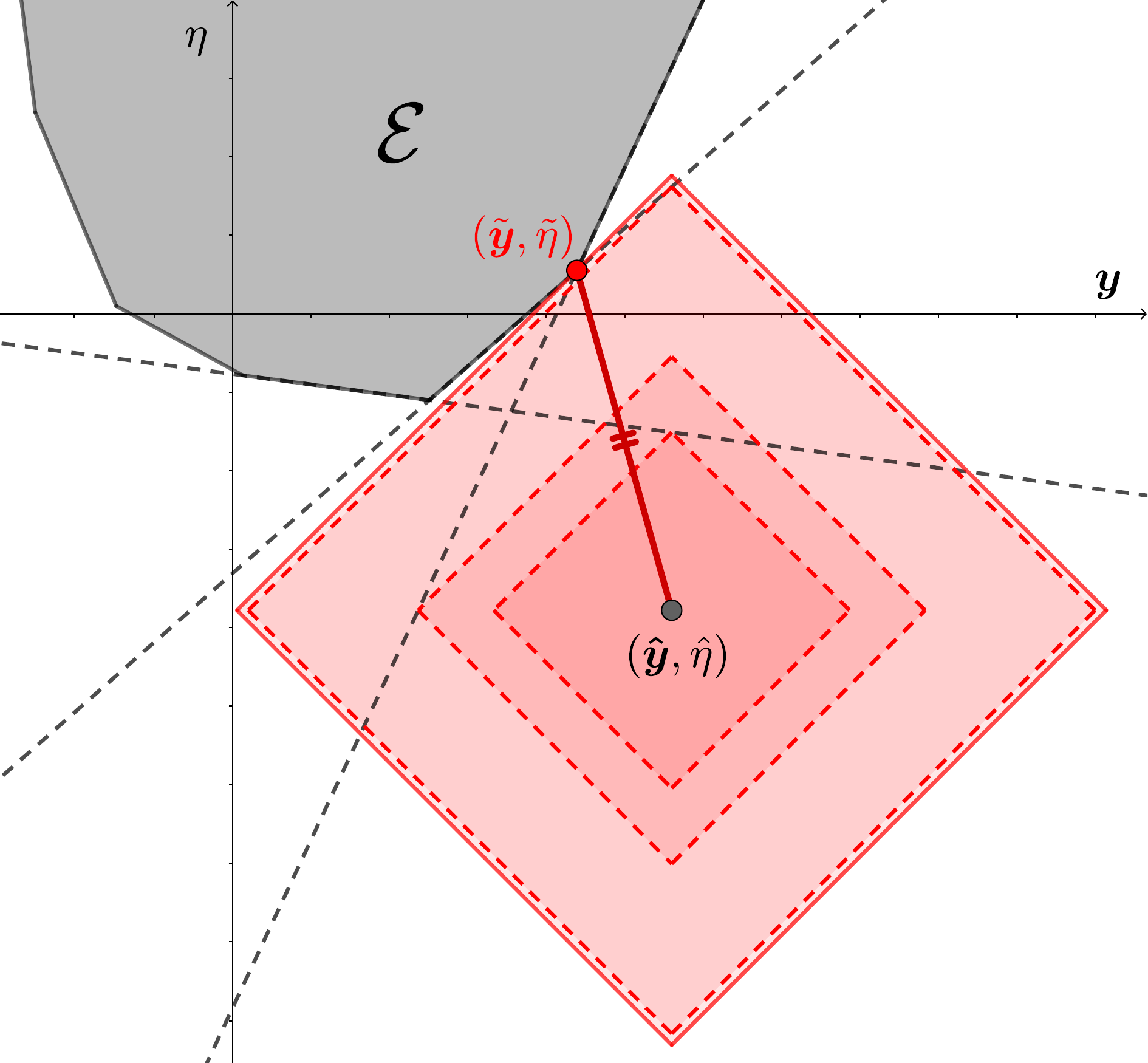}\label{fig:deep_cuts_primal_norm-linf}}\hspace{0.2cm}
    \caption{Primal-dual perspectives of the separation problem. 
    Separating $(\hat{\yy},\hat{\eta})\notin \epi$ with an $\ell_p$-deepest cut accounts for finding $(\tilde{\yy},\tilde{\eta})\in\epi$ with minimum $\ell_q$-distance to $(\hat{\yy},\hat{\eta})$.
    The red lines illustrate the contour lines of the objective value of SSP, which also correspond to $\ell_q$-balls around $(\hat{\yy},\hat{\eta})$.}\label{fig:deep_cuts_primal_norm}
\end{figure*}

\begin{theorem}\label{thrm:deepest-cuts-primal}
Separation problem \eqref{separation-subproblem-Lp} is equivalent to the following Lagrangian dual problem
\begin{equation}\label{primal-ssp}
\begin{split}
     [\text{Primal SSP}]\quad \min\quad & \|({\yy}-\hat{\yy}, \eta-\hat{\eta})\|_{q}\\
    \text{s.t.} \quad & \eta\ge {\cc}^{\top}{\xx}+{\ff}^{\top}{\yy}\\
    & A{\xx}\ge {\bb}-B{\yy}\\
    & {\xx} \ge {\vzero}, \eta\ge \hat{\eta},
\end{split}
\end{equation}
in which $({\yy},{\xx},\eta)$ are the variables and $\ell_q$ is the dual norm of $\ell_p$.
\end{theorem}
The following result follows from strong duality and the definition of $\epi$.
\begin{corollary}\label{col-deepest-cut-projection}
    $d^*_{\ell_p}(\hat{\yy},\hat{\eta})$ measures the $\ell_q$ distance of $(\hat{\yy},\hat{\eta})$ from $\epi$. That is,
    \begin{align}
        d^*_{\ell_p}(\hat{\yy},\hat{\eta})=\min_{({\yy},\eta)\in \epi: \eta\ge \hat{\eta}}\|{\yy}-\hat{\yy}, \eta-\hat{\eta}\|_q.\label{eq:deepest-cut-projection}
    \end{align}
\end{corollary}
Let $(\tilde{\yy},\tilde{\eta})$ be the optimal solution of the Primal SSP.
In convex analysis, the solution of \eqref{eq:deepest-cut-projection} for $q=2$ is known as the \textit{projection} of $(\hat{\yy},\hat{\eta})$ onto $\epi$. Thus, we refer to $(\tilde{\yy},\tilde{\eta})$ henceforth as the $\ell_q$-projection of $(\hat{\yy},\hat{\eta})$, and refer to \eqref{eq:deepest-cut-projection} as the \textit{projection subproblem}. 
Figure~\ref{fig:deep_cuts_primal_norm} illustrates these projections for different values of $q$. Observe that the $\ell_q$-projection or the $\ell_p$-deepest cut might not be unique for $p=1$ ($q=\infty$) or $p=\infty$  ($q=1$).
The following proposition states that deepest cuts support $\epi$ at $\ell_q$-projections, even when the projection or the cut are not unique.

\begin{proposition}\label{prop:deepest-cuts-support}
    Let $(\tilde{\yy},\tilde{\eta})\in\epi$ be an $\ell_q$-projection of  $(\hat{\yy},\hat{\eta})$ onto $\epi$. 
    Then, any $\ell_p$-deepest cut separating  $(\hat{\yy},\hat{\eta})$ from $\epi$ supports $\epi$ at $(\tilde{\yy},\tilde{\eta})$.
\end{proposition}

From the duality between separation and projection established in the above theorem and proposition, we derive the following important technical results. 

\textbf{Deepest cuts minimally resolve infeasibility in FSP.} 
By Theorem~\ref{thrm:deepest-cuts-primal} and as illustrated in Figure~\ref{fig:deep_cuts_primal_norm}, producing an $\ell_p$-deepest cut amounts to finding the point $(\tilde{\yy},\tilde{\eta})$ of least $\ell_q$-distance from $(\hat{\yy},\hat{\eta})$ for which a feasible solution ${\xx}$ exists that satisfies the system
\begin{align*}
    \left\{{\cc}^{\top}{\xx} \le \tilde{\eta}-{\ff}^{\top}\tilde{\yy};\quad A{\xx}\ge {\bb}-B\tilde{\yy};\quad {\xx} \ge {\vzero}\right\}.
\end{align*}
Hence, producing a deepest cut can be viewed as resolving infeasibility of FSP \eqref{generic-milp-fsp} through minimal deviation from $(\hat{\yy},\hat{\eta})$ with respect to the $\ell_q$-norm.
If FSP is feasible for $(\hat{\yy},\hat{\eta})$ (i.e., if  ${\|\tilde{\yy}-\hat{\yy},\tilde{\eta}-\hat{\eta}\|_q=0}$), then $(\hat{\yy},\hat{\eta})$ is optimal for MP. Effectively, $d^*_{\ell_p}(\hat{\yy},\hat{\eta})$ measures how far $(\hat{\yy},\hat{\eta})$ is from being optimal by measuring the minimal deviation in $(\hat{\yy},\hat{\eta})$ that renders FSP feasible. Thus, producing a deepest cut assesses how inaccurate our current guess of the optimal solution is.

\textbf{Sparsity, density and flatness of deepest cuts.}
We have empirically observed that the deepest cuts generated at the early stages of the BD algorithm tend to be relatively flat. That is, the coefficients of the $\yy$-variables in the cut are mostly zero, and in some cases the cut is completely flat, i.e., all $\yy$ coefficients are zero. Here, we provide a justification for this observation and discuss its implications. Along this vein, we first note the following property of $\ell_1$-deepest cuts.

\begin{proposition}\label{prop:deepst_cuts_flat}
For sufficiently small $\hat{\eta}$, the $\ell_1$-deepest cut separating $(\hat{\yy},\hat{\eta})$ from $\epi$ is the flat cut $\eta\ge Q^*$, where $Q^*=\min_{\yy}Q(\yy)$ is the optimal value of $Q$ for unrestricted $\yy$.
\end{proposition}

Proposition~\ref{prop:deepst_cuts_flat} implies that, at early iterations of the BD algorithm, an $\ell_1$-deepest cut can provide a lower bound of at least $Q^*$. Since $Q^*$ is obtained by relaxing $\yy\in Y$, $Q^*$ is at most equal to the optimal value of the LP relaxation of OP, thus the bound is expected to be effective when the integrality gap is low (e.g., in facility location problems). 

More generally, for small $p$ (i.e., large $q$) and relatively small $\hat{\eta}$, we may approximate ${\|({\yy}-\hat{\yy}, \eta-\hat{\eta})\|_{q}\approx \eta-\hat{\eta}}$. Therefore, in line with Proposition~\ref{prop:deepst_cuts_flat}, we can expect that the coefficients of the $\yy$-variables in the $\ell_p$-deepest cuts (i.e., $\hat\pi_0{\ff}^{\top}-\hat{\BFpi}^{\top}B$) will be close to zero, which means that the deepest cuts are relatively \textit{sparse}, which is also in line with using Lasso or $\ell_1$-regularization in statistics for producing sparse solutions \citep[see e.g.,][]{tibshirani1996regression}. By the same token, large values of $p$ (e.g., $p=\infty$) induce \textit{dense} cuts, in that the coefficients of the $\yy$ variables are mostly non-zero.

\textbf{Deepest cuts are more likely to be optimality cuts than feasibility cuts.} 
In our experiments with deepest cuts, we found that they are likely to be optimality cuts, even when classical Benders produces a feasibility cut (i.e., even when $\hat{\yy}\notin \dom(Q)$). Intuitively, since at each iteration of BD, $\eta^{\itt}$ is an under-estimator of the convex piece-wise linear function $Q$ and deepest cuts support $\epi$, the coefficient of $\eta$ in a deepest cut (namely, $\pi_0$) is most likely non-zero (i.e., the cut is an optimality cut).
Proposition~\ref{prop:deepst_cuts_optimality} below guarantees that deepest cuts are optimality cuts when $\hat{\eta}$ is sufficiently small. Proposition~\ref{prop:deepst_cuts_optimality2} below further suggests that, even if $\hat{\eta}$ is not very small, deepest cuts are more likely to be optimality cuts (note that the $\ell_q$-projection is always unique for $1<q<\infty$).
This is particularly appealing from a practical standpoint, as the contribution of Benders optimality cuts to closing the gap is usually more pronounced than that of feasibility cuts \citep[cf.,][]{saharidis2010improving,de2013improved}. 

\begin{proposition}\label{prop:deepst_cuts_optimality}
For $p>1$, provided that $\hat{\eta} < Q^*\coloneqq\min\limits_{\yy}Q(\yy)$, the $\ell_p$-deepest cut(s) separating $(\hat{\yy},\hat{\eta})$ are optimality cuts for any arbitrary $\hat{\yy}$ (i.e., even if $\hat{\yy}\notin \dom(Q)$).
\end{proposition}

\begin{proposition}\label{prop:deepst_cuts_optimality2}
Provided that the $\ell_q$-projection of $(\hat{\yy},\hat{\eta})$ onto $\epi$ is the unique point $(\tilde{\yy},\tilde{\eta})$ and $\hat{\eta} < \tilde{\eta}$, the $\ell_p$-deepest cuts separating $(\hat{\yy},\hat{\eta})$ are optimality cuts even if $\hat{\yy}\notin \dom(Q)$.
\end{proposition}

\section{General Benders Distance Functions }\label{sec:general-distance-functions}

We now introduce a general distance function based on duality theory that we call a \textit{Benders distance function}. Such a function $d$ must (a) identify (based on its sign) whether the incumbent point $(\hat{\yy},\hat \eta)$ is inside, outside, or on the boundary of the epigraph $\epi$, and (b) if outside, convey how ``far" the incumbent point is from the boundary. Crucially, we do not explicitly define the metric on which a Benders distance function $d$ is based; this is by design, and interestingly isn't needed. We provide a formal definition of Benders distance functions in Appendix \ref{app:monotonicity} and show that so long as a monotonicity property linked to convexity holds, a sufficient notion of distance exists.

Following our methodology, cutting off the incumbent point $(\hat{\yy},\hat{\eta})$ at maximum distance from the epigraph can be achieved by maximizing the Benders distance function $d$ with respect to all feasible dual solutions, yielding in the process the optimal value function $d^*$ below, which we call the \textit{epigraph distance function}.
\begin{align}
    [\mbox{BSP}]\;d^*(\hat{\yy},\hat{\eta})=\max_{({\BFpi},\pi_0)\in\Pi} \quad d(\hat{\yy},\hat{\eta}|{\BFpi},\pi_0). \label{eq:separation_general}
\end{align}

We start in Section \ref{sec:normalized-distance-functions} by introducing an important special class of Benders distance function, normalized distance functions. In Section \ref{sec:projective}, we introduce projective distance functions and show that they admit a simple characterization of the distance to candidate hyperplanes through the gradient of the normalization function. In Section \ref{sec:linear_pseudonorms} we connect deepest cuts with other types of Benders cuts from the literature through  linear normalization functions.
Finally, in Section \ref{sec:modified-bd-algorithm} we study the convergence behavior of our distance-based BD algorithm.

\subsection{Normalized Distance Functions and Normalization Functions}\label{sec:normalized-distance-functions}

Here we introduce an important special class of Benders distance functions, which we call \textit{normalized Benders distance functions}. 
The distance functions in this class generalize the geometric $\ell_p$-norm-based distance function that we introduced in Section \ref{sec:lp-norm-deepest-cuts} by replacing the denominator of  $d_{\ell_p}$ \eqref{distance-Lp} with a general \textit{normalization function} which is only required to be a positive homogeneous function of the dual variables.

\begin{definition}[Normalized distance function]
Let $d_h(\hat{\yy},\hat \eta|{\BFpi},\pi_0) = \frac{{\BFpi}^{\top} ({\bb}-B\hat {\yy})+\pi_0({\ff}^{\top}\hat{\yy}-\hat\eta)}{h({\BFpi},\pi_0)}$ where $h$ is a positive homogeneous function (i.e., $h(\alpha{\BFpi},\alpha\pi_0)=\alpha h({\BFpi},\pi_0)$ for any $\alpha\ge0$). We call $d_h$ a \textbf{normalized distance function}, and refer to $h$ as its \textbf{normalization function}.
\end{definition} 

As noted in Appendix \ref{app:monotonicity}, normalized distance functions are well-defined Benders distance functions in that they produce cuts that correctly separate the infeasible solutions from $\epi$ and their value monotonically increases as we move away from the boundary of $\epi$.
As we will show, many types of Benders cuts proposed in the literature, such as Minimum Infeasible Subsystems (MIS) cuts \citep{fischetti2010note}, can be viewed as being generated using normalized distance functions.

The normalization function $h$ governs the behavior of the distance function, and quantifies our perception of the quality of the cut it produces.
Homogeneity of $h$ is critical. Indeed, with constant (i.e., non-homogeneous) $h({\BFpi},\pi_0)=1$, \eqref{eq:separation_general} is equivalent to the na\"ive CGSP \eqref{generic-milp-cgsp}, for which $d^*(\hat{\yy},\hat\eta)\in\{0,+\infty\}$. In this case, $d^*$ is simply a characteristic function of $\epi$, which, regardless of the quality of the cut, merely provides a binary answer to whether or not $(\hat{\yy},\hat\eta)$ is the optimal solution, without any further indication of how far $(\hat{\yy},\hat\eta)$ is from being optimal.

Our next result sheds light on the desirable properties of the normalization function $h$, and paves the way for reformulating the separation problem. 
Note when $h$ is linear, the reformulation reduces to the Charnes-Cooper transformation \citep{charnes1962programming} for linear-fractional programs.
\begin{proposition}\label{prop:subproblem-transformation}
Let $h$ be a normalization function, $d_h(\hat{\yy},\hat{\eta}|{\BFpi},\pi_0) = \frac{{\BFpi}^{\top} ({\bb}-B\hat{\yy})+\pi_0({\ff}^{\top}\hat{\yy}-\hat{\eta})}{h({\BFpi},\pi_0)}$ the distance function induced by $h$, and $\Pi_h=\{({\BFpi},\pi_0)\in\Pi:h({\BFpi},\pi_0)\le 1\}$ the cone $\Pi$ truncated by the constraint $h({\BFpi},\pi_0)\le 1$.
Then, the separation problem \eqref{eq:separation_general} is equivalent to the normalized separation problem (NSP) defined below. Furthermore, $h({\BFpi},\pi_0)\le 1$ is binding at optimality.
\begin{align}
     [\text{NSP}]\quad d^*_h(\hat{\yy},\hat{\eta})=\max_{({\BFpi},\pi_0)\in\Pi_h} \;&{\BFpi}^{\top} ({\bb}-B\hat{\yy})+\pi_0({\ff}^{\top}\hat{\yy}-\hat{\eta}) \label{normalized-separation-problem}
\end{align}
\end{proposition}

From a polyhedral perspective, Proposition~\ref{prop:subproblem-transformation} shows the equivalence between choosing a distance function and truncating the cone $\Pi$ with a specific normalization function. 
Note that all general norms (including the $\ell_p$-norms introduced in Section~\ref{sec:lp-norm-deepest-cuts}) are normalization functions. 
Figure~\ref{fig:norms-truncate} in Appendix~\ref{app:supplementary_figures} illustrates the effect of $\ell_1$-, $\ell_2$- and $\ell_{\infty}$-norms on truncating the cone of dual solutions.

Proposition~\ref{prop:subproblem-transformation} also shows that BSP can be converted into a problem of optimizing a linear function over a convex set $\Pi_h$ when $h$ is convex, which also allows us to leverage the re-optimization capabilities of the solver whenever possible. For instance, a certificate produced at iteration $t$ of the BD algorithm can be used for warm-starting the separation subproblem at iteration $t+1$ if $\Pi_h$ does not depend on $({\yy}^{\itt},\eta^{\itt})$.
Moreover, a convex piece-wise linear function $h({\BFpi},\pi_0)$ amounts to solving linear programs with different objective function coefficients at each iteration (see Section~\ref{sec:linear_pseudonorms} and Appendix~\ref{app:reformulations-lp}), thus can be re-optimized using a primal simplex method (see Appendix~\ref{app:reoptimization}).

We end this section by providing a primal view of normalized distance functions, which  shows that producing a cut with respect to $d_h$ amounts to resolving infeasibility in FSP \eqref{generic-milp-fsp} through minimal perturbation in its right-hand-side values along $\nabla_{(\BFpi,\pi_0)}h(\BFpi^*,\pi^*_0)$. The proof uses Euler's homogeneous function theorem, which implies that $h(\BFpi,\pi_0)=\BFpi^{\top}\nabla_{\BFpi}h+\pi_0\nabla_{\pi_0}h$.

\begin{theorem}\label{thrm:homogeneous-functions-primal}
    Let $h(\BFpi,\pi_0)$ be a convex differentiable positive homogeneous function. Assuming that NSP \eqref{normalized-separation-problem} admits a bounded optimal solution $(\BFpi^*, \pi_0^*)$, NSP is equivalent to the following LP
    \begin{equation}\label{primal-homogeneous}
    \begin{split}
        d^*_h(\hat{\yy},\hat{\eta})=\min\quad & z\\
        \text{s.t.} \quad & A\xx \ge \bb-B\hat{\yy}-z\nabla_{\BFpi}h(\BFpi^*,\pi^*_0)\\
        & \cc^{\top}\xx\le \hat{\eta}-\ff^{\top}\hat{\yy}+z\nabla_{\pi_0}h(\BFpi^*,\pi^*_0)\\
        & \xx\ge \vzero.
    \end{split}
    \end{equation}
\end{theorem}

\subsection{Projective Normalization Functions and Projective Distance Functions.}\label{sec:projective}

We now introduce a class of normalized distance functions that possess several properties of geometric distance functions, including defining a metric for measuring distance between the incumbent point and the hyperplanes/epigraph and establishing a duality between separation and projection.

\begin{definition}[Projective normalization function]
Let $\BFtau=B^{\top}\BFpi-\pi_0 \ff$ for $(\BFpi,\pi_0)\in \Pi$. We call 
a positive homogeneous function $g(\BFtau,\pi_0)$ a \textit{projective normalization function} iff $g(\BFtau,\pi_0) > 0$ for $(\BFtau,\pi_0)\ne \vzero$, and refer to the distance function induced by $g$ as a \textit{projective distance function}.
\end{definition}

We start by showing that projective distance functions admit a simple characterization of the distance between  incumbent point and  candidate hyperplane, which is simply the $1$-dimensional distance along the gradient of $g$ from the incumbent point to the hyperplane.

\begin{proposition}\label{prop-projective-functions-distance-hyperplane}
    Let $g$ be a projective normalization function. For any $(\hat{\yy}, \hat{\eta})$ and any $(\BFpi,\pi_0)\in\Pi$ such that $(\BFtau, \pi_0)\ne\vzero$ with $\BFtau=B^{\top}\BFpi-\pi_0 \ff$, the line $(\yy, \eta) = (\hat{\yy},\hat{\eta})+ z\nabla_{(\BFtau, \pi_0)}g(\BFtau, \pi_0)$ intersects the hyperplane $\hp(\BFpi,\pi_0)$ at $(\tilde{\yy}, \tilde{\eta})$ for $z=d_{g}(\hat{\yy},\hat{\eta}|\BFpi,\pi_0)$, thus implicitly defining $d_g$.
\end{proposition}

Some examples of projective normalization functions include $\|\BFtau,\pi_0\|_p$ as in $d_{\ell_p}$ distance functions, as well as some linear normalization functions of the form $g(\BFtau,\pi_0)=\BFomega^{\top}\BFtau+\omega_0\pi_0$. More specifically, substituting $\BFtau=B^{\top}\BFpi-\pi_0 \ff$, we obtain $h(\BFpi,\pi_0)=\BFpi^{\top}B\BFomega+\pi_0(\omega_0-\ff^{\top}\BFomega)$. Therefore, for $h(\BFpi,\pi_0)=\BFpi^{\top}\ww+w_0\pi_0$ to be projective it is necessary that $\ww\in \text{Im}(B)$ (i.e., there exists $\BFomega\in\mathbb{R}^n$ such that $\ww=B\BFomega$.) 
We provide two examples of such linear functions in Sections \ref{sec:relaxed-l1-norm} and \ref{sec:cw-norm}. 

Our next result establishes a relationship between $d_{\ell_p}$ distance functions and general projective distance functions by using the fact that  $d_{\ell_p}$ measures the $\ell_q$-distance to the candidate hyperplanes.

\begin{proposition}\label{prop:projective-functions-bound_lp}
    Let $g$ be a projective normalization function. 
    For any $p\ge 1$ and its dual $q$ (i.e., $\frac{1}{p}+\frac{1}{q}=1$), any $(\hat{\yy}, \hat{\eta})$ and any $(\BFpi,\pi_0)\in\Pi$ with $\BFtau=B^{\top}\BFpi-\pi_0 \ff$ such that $(\hat{\yy}, \hat{\eta})\notin \cut(\BFpi, \pi_0)$ we have
    \begin{align*}
        d_{\ell_p}(\hat{\yy},\hat{\eta}|\BFpi,\pi_0)\le d_{g}(\hat{\yy},\hat{\eta}|\BFpi,\pi_0)\; \|\nabla_{(\BFtau, \pi_0)} g(\BFtau, \pi_0)\|_q.
    \end{align*}
\end{proposition}

\begin{figure*}[h]
    \centering
    \includegraphics[clip,width=0.43\textwidth]{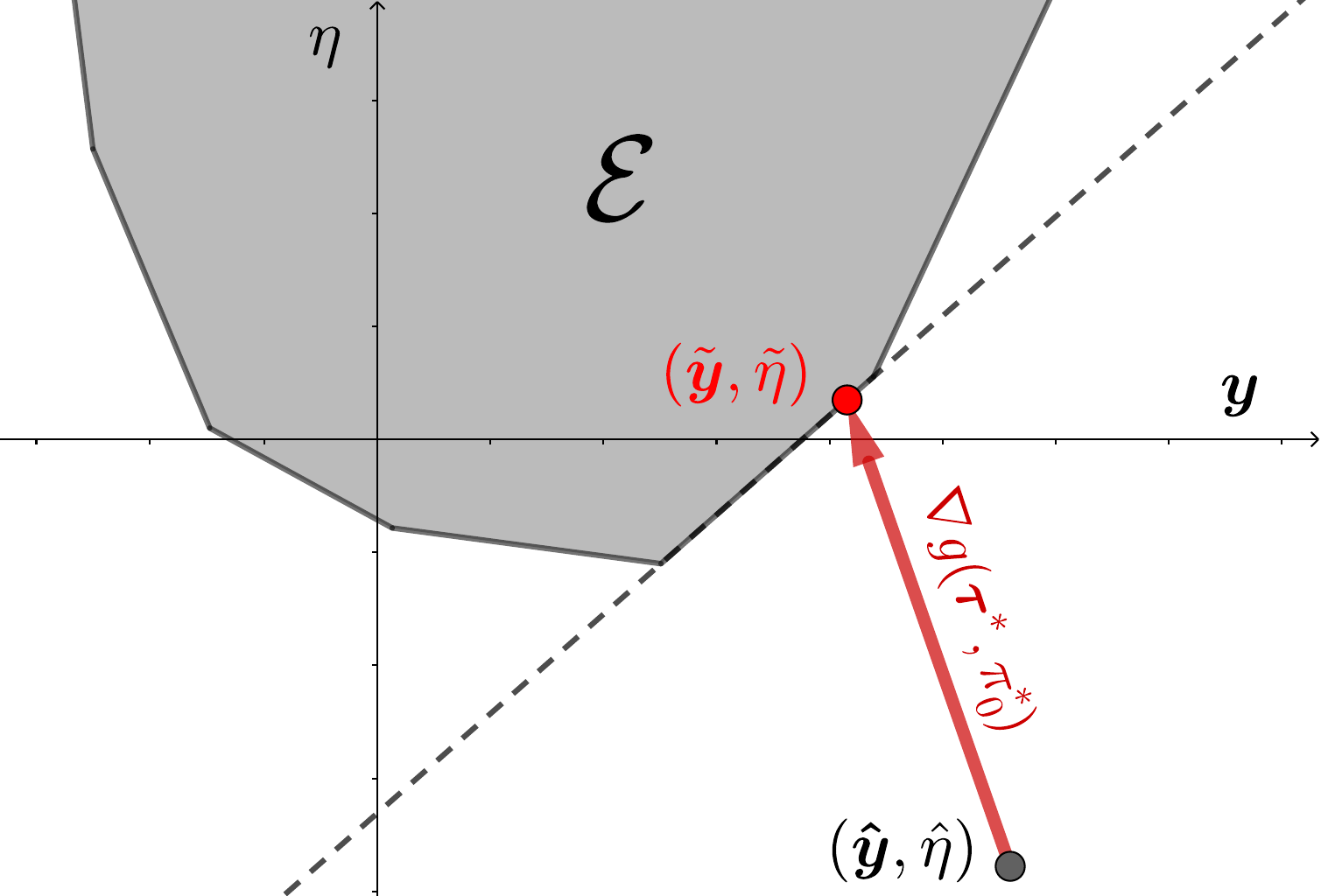}
    \caption{Primal view of projective normalization function $g$, and $(\tilde{\yy},\tilde{\eta})\in \epi$, the $g$-projection of $(\hat{\yy},\hat{\eta})\notin \epi$.}
    \label{fig:projective}
\end{figure*}

Next, akin to Theorem~\ref{thrm:deepest-cuts-primal}, Proposition~\ref{prop:projective-functions-primal} below  establishes a duality between separation and projection for a projective normalization function $g$ by explicitly yielding a point $(\tilde{\yy}, \tilde{\eta})\in\epi$ (the optimizer of \eqref{primal-projective} below) which we call the $g$-projection. 
This result together with Proposition \ref{prop-projective-functions-distance-hyperplane} imply that the cuts selected according to $d_g$ support $\epi$ at the $g$-projection (see Figure~\ref{fig:projective}.)

\begin{proposition}\label{prop:projective-functions-primal}
   Let $g$ be a projective normalization function and $(\BFpi^*, \pi_0^*)$ an optimal solution to the separation problem NSP \eqref{normalized-separation-problem} with $h(\BFpi,\pi_0)=g(\BFtau,\pi_0)$ where $\BFtau=B^{\top}\BFpi-\pi_0 \ff$. Then
    \begin{equation}\label{primal-projective}
    \begin{split}
         d^*_g(\hat{\yy},\hat{\eta})=\min \; \{z: (\yy,\eta)\in \epi,\quad  (\yy, \eta) = (\hat{\yy},\hat{\eta})+ z\nabla_{(\BFtau, \pi_0)}g(\BFtau^*, \pi_0^*)\}.
    \end{split}
    \end{equation}
\end{proposition}

The separation-projection duality also implies that we can extract a deep Benders cut by solving a classical Benders separation problem \eqref{eq:projective-primal-dsp} when the $g$-projection $(\tilde{\yy}, \tilde{\eta})$ is known; we use this property in Appendix~\ref{app:separable_subproblem} to develop an algorithm for solving separable subproblems.
\begin{align}
    \max_{\uu\in \mathcal{U}} \quad \uu^{\top}(\bb-B\tilde{\yy}).\label{eq:projective-primal-dsp}
\end{align}

Another remarkable property of a projective $g$ is that $g(\BFtau,\pi_0)\le 1$ directly truncates the cone $\Gamma$ defined in \eqref{eq:Gamma}. Given the one-to-one correspondence between $\Gamma$ and the F-cone of $\epi$ \citep[cf.][for the definition of F-cone of a convex set]{conforti2014integer}, when the optimal solution $(\BFtau^*,\pi^*_0)$ is an extreme ray of $\Gamma$, the resulting cut almost surely exposes a facet of $\epi$ \citep{conforti2019facet}. 
\begin{align}
  \Gamma = \{(\BFtau,\pi_0): \exists (\BFpi,\pi_0)\in \Pi \text{ s.t. } \BFtau=B^{\top}\BFpi-\pi_0 \ff\}  \label{eq:Gamma}
\end{align}

\subsection{Linear Pseudonorms}\label{sec:linear_pseudonorms}
Consider the class of normalization functions defined by choosing parameters $({\ww},w_0)$ such that $h({\BFpi},\pi_0)={\BFpi}^{\top}{\ww}+\pi_0w_0\ge 0$ for all $({\BFpi},\pi_0)\in\Pi$. Here, we study how different values of $({\ww},w_0)$ impact the resulting normalized Benders distance function, as well as how the cuts produced relate to other cut selection strategies in the literature. Note that a linear function $h$ of this form satisfies most axioms of a norm; that is, $h$ is subadditive (i.e., $h({\uu}+{\vv})\le h({\uu})+h({\vv})$), homogeneous (i.e., $h(\alpha {\uu})=\alpha h({\uu})$ for any $\alpha\ge 0$), and positive over $\Pi$, but not necessarily positive definite (i.e., $h({\BFpi},\pi_0)=0$ does not necessarily imply $({\BFpi},\pi_0)={\vzero}$). Hence, we call $h$ a \textit{linear pseudonorm} over $\Pi$.

With $h({\BFpi},\pi_0)={\BFpi}^{\top}{\ww}+\pi_0w_0 \ge 0$ for all $(\BFpi,\pi_0)\in \Pi$, 
using Proposition~\ref{prop:subproblem-transformation}, we can re-state the separation problem as the following LP, which contains only one additional variable and one additional constraint compared to DSP \eqref{generic-milp-projection-dsp} in the classical BD algorithm.
\begin{equation}
    \begin{split}\label{sp-linear-norm-fractional-linear}
    \max\quad & {\BFpi}^{\top}({\bb}-B\hat{\yy})+\pi_0({\ff}^{\top}\hat{\yy}-\hat{\eta})\\
    \text{s.t.} \quad & {\BFpi}^{\top}A\le \pi_0 {\cc}^{\top}\\
    & {\BFpi}^{\top}{\ww}+\pi_0w_0 \le 1\\
    & {\BFpi} \ge {\vzero}, \pi_0\ge 0.
\end{split}
\end{equation} 

Separation problem \eqref{sp-linear-norm-fractional-linear} is the \textit{minimal infeasible subsystems} (MIS) subproblem proposed by \cite{fischetti2010note}, which is derived by treating the separation problem as approximating the minimal source of infeasibility of FSP \eqref{generic-milp-fsp} by minimizing a positive linear function ${\BFpi}^{\top}{\ww}+\pi_0 w_0$ over the alternative polyhedron of $\Pi$ (i.e., $\Pi$ truncated by constraint ${\BFpi}^{\top}({\bb}-B\hat{\yy})+\pi_0({\ff}^{\top}\hat{\yy}-\hat{\eta})=1$). 
Therefore, the MIS subproblem can be viewed as a special type of the Benders separation subproblem \eqref{normalized-separation-problem} in which the normalization function $h$ takes the form of $h({\BFpi},\pi_0)={\BFpi}^{\top}{\ww}+\pi_0w_0$.

As noted by \cite{fischetti2010note}, the choice of parameters $({\ww}, w_0)$ can profoundly impact the effectiveness of MIS cuts. In their implementation, the authors set $w_0=1$ and initially set $w_i=1$ for all $i=1,\dots,m$. They further suggest that setting $w_i=0$ for the null rows of $B$ (i.e., row $i$ such that $B_{ij}=0$ for all $j$) may lead to a substantial improvement in the convergence of the BD algorithm. 
Below, we propose four ways to choose parameters $({\ww}, w_0)$ based on the parameters of the problem instance and discuss their implications and connections to other cut selection strategies.

\subsubsection{Classical Benders pseudonorm.}\label{sec:bd_pseudonorm}
A trivial choice is to set $\ww=\vzero$ and $w_0=1$, 
which causes \eqref{sp-linear-norm-fractional-linear} to reduce to DSP \eqref{generic-milp-projection-dsp} in the classical BD algorithm (here, $\pi_0^* = 1$). In other words, defining
$d_{\text{CB}}(\hat{\yy},\hat{\eta}|{\BFpi},\pi_0) = \frac{{\BFpi}^{\top} ({\bb}-B\hat{\yy})+\pi_0({\ff}^{\top}\hat{\yy}-\hat{\eta})}{\pi_0}$ as the classical Benders (CB) distance function, then
$$d^*_{\text{CB}}(\hat{\yy},\hat{\eta})=\max_{({\BFpi},\pi_0)\in\Pi}d_{\text{CB}}(\hat{\yy},\hat{\eta}|{\BFpi},\pi_0)=\tilde{Q}(\hat{\yy})+{\ff}^{\top}\hat{\yy}-\hat{\eta}=Q(\hat{\yy})-\hat{\eta}.$$ 
Observe that, at iteration $t$ of BD, $d^*_{\text{CB}}({\yy}^{\itt},\eta^{\itt})=Q({\yy}^{\itt})-\eta^{\itt}$ overestimates the optimality gap. Thus, $d^*_{\text{CB}}({\yy}^{\itt},\eta^{\itt})=0$ means $({\yy}^{\itt},\eta^{\itt})$ is an optimal solution to \eqref{generic-milp-projection-reform}, which is precisely the stopping criterion used in the classical BD algorithm. Moreover, while $h_{\text{CB}}(\BFpi,\pi_0)\coloneqq\pi_0$ seems to satisfy the requirements of a projective normalization function, $h_{\text{CB}}$ is not projective in general because DSP \eqref{generic-milp-projection-dsp} need not be bounded. However, when $d^*_{\text{CB}}(\hat{\yy},\hat{\eta}) < \infty$ (i.e., $\hat{\yy}$ is feasible), the projection direction is $(\vzero, 1)$. Therefore, as illustrated in Figure~\ref{fig:cut-selection-deep-classical}, $d^*_{\text{CB}}(\hat{\yy},\hat{\eta})$ can be geometrically interpreted as the distance from the point $(\hat{\yy},\hat{\eta})$ to the boundary of $\epi$ along the $\eta$-axis.  

\subsubsection{Relaxed \texorpdfstring{$\ell_1$}{L1} pseudonorm.}\label{sec:relaxed-l1-norm}
Here, we introduce a linear normalization function based on $\ell_1$ which not only assigns meaningful values to parameters $({\ww}, w_0)$, but also leads to a geometric interpretation of the MIS subproblem.
Expanding $\ell_1$ and using the triangle inequality, we obtain 
$$\|({\BFpi}^{\top} B - \pi_0{\ff}^{\top},\pi_0)\|_1 =\pi_0+\sum\nolimits_{j=1}^{n}|\pi_0f_j-\sum\nolimits_{i=1}^{m}\pi_i B_{ij}|\le \pi_0(1+\sum\nolimits_{j=1}^{n}|f_j|)+\sum\nolimits_{i=1}^{m}\pi_i\sum\nolimits_{j=1}^{n}|B_{ij}|.$$
Hence, we refer to $h_{R\ell_1}({\BFpi},\pi_0)={\BFpi}^{\top}{\ww}+\pi_0w_0$ with $w_0=1+\sum_{j=1}^{n}|f_j|$ and $w_i=\sum_{j=1}^{n}|B_{ij}|$ as the \textit{relaxed $\ell_1$ pseudonorm} ($R\ell_1$ for short), and denote $d_{R\ell_1}(\hat{\yy},\hat{\eta}|{\BFpi},\pi_0)= \frac{{\BFpi}^{\top} ({\bb}-B\hat{\yy})+\pi_0({\ff}^{\top}\hat{\yy}-\hat{\eta})}{h_{R\ell_1}({\BFpi},\pi_0)}$.
Note that $w_i=\sum_{j=1}^{n}|B_{ij}|=0$ for the null rows of $B$, which is in line with the recipe of \cite{fischetti2010note}. 
Proposition \ref{prop:depth-gap-lp} below also shows that $d_{\ell_p}$ distance functions are bounded by $d_{CB}$ and $d_{R\ell_1}$.

\begin{proposition}\label{prop:depth-gap-lp}
The following relationship  holds between $d_{\text{CB}}$, $d_{\ell_p}$, and $d_{R\ell_1}$ for any $(\BFpi,\pi_0)\in \Pi$:
\begin{align*}
    d_{\text{CB}}(\hat{\yy},\hat{\eta}|\BFpi,\pi_0)\ge d_{\ell_{\infty}}(\hat{\yy},\hat{\eta}|\BFpi,\pi_0)\ge\dots\ge d_{\ell_p}(\hat{\yy},\hat{\eta}|\BFpi,\pi_0)\ge \dots\ge d_{\ell_1}(\hat{\yy},\hat{\eta}|\BFpi,\pi_0)\ge d_{R\ell_1}(\hat{\yy},\hat{\eta}|\BFpi,\pi_0).
\end{align*}
\end{proposition}

A stronger connection between $d_{\ell_p}$ and $d_{R\ell_1}$ exists when $B\ge 0$ (i.e., $B_{ij}\ge 0$ for all $i$ and $j$) or $B\le 0$ (i.e., $B_{ij}\le 0$ for all $i$ and $j$), which is the case in many structured MILPs, such as facility location, network design and network interdiction problems. When $B\ge 0$, $\sum_{i=1}^{m}\pi_i\sum_{j=1}^{n}|B_{ij}|=\BFpi^{\top}B\ee$, where $\ee$ is the vector of all ones in $\mathbb{R}^n$. 
This allows us to restate $h_{R\ell_1}$ as $g(\BFtau,\pi_0)=\BFtau^{\top}\ee+\pi_0(1+\ff^{\top}\ee+\sum_j |f_j|)$, where  $\BFtau=B^{\top}\BFpi-\pi_0\ff$.
Similarly, when $B\le 0$, $g(\BFtau,\pi_0)=-\BFtau^{\top}\ee+\pi_0(1-\ff^{\top}\ee+\sum_j|f_j|)$.
Consequently, $R\ell_1$ is projective, and using Propositions \ref{prop:projective-functions-bound_lp} and \ref{prop:depth-gap-lp}, we derive the following result.
As a result, when the factor $\frac{1}{\|(\BFdelta, 1+\ff^{\top}\BFdelta+\sum_j |f_j|)\|_q}$ is reasonably close to 1, we can replace $d_{\ell_p}$ with $d_{R\ell_1}$ and still produce cuts that are  geometrically deep  but by solving a simpler LP.
\begin{corollary}\label{corollary:lp-rl1}
Provided that $B\ge 0$ or $B\le 0$, for any $p\ge 1$ and any $(\hat{\yy},\hat{\eta})$ we have:
$$d_{\ell_p}(\hat{\yy},\hat{\eta}|\BFpi,\pi_0)\ge d_{R\ell_1}(\hat{\yy},\hat{\eta}|\BFpi,\pi_0) \ge \frac{1}{\|(\BFdelta, 1+\ff^{\top}\BFdelta+\sum_j |f_j|)\|_q}d_{\ell_p}(\hat{\yy},\hat{\eta}|\BFpi,\pi_0), \qquad \forall (\BFpi,\pi_0)\in \Pi$$
where $\BFdelta=\ee$ when $B\ge 0$ and $\BFdelta=-\ee$ when $B\le 0$,  and $\ell_q$ is the dual norm of $\ell_p$.
\end{corollary}

\subsubsection{Magnanti-Wong-Papadakos pseudonorm.}\label{sec:mw-norm}
The Magnanti-Wong procedure for producing a Pareto-optimal cut using a given core point $\bar{\yy}$ (i.e., $\bar{\yy}\in \mbox{relint}(Y)$) involves solving the following subproblem \citep{magnanti1981accelerating}:
\begin{align}
    \max_{\uu\in \mathcal{U}}\;\left\{ {\uu}^{\top}({\bb}-B\bar{\yy}): {\uu}^{\top}({\bb}-B \hat{\yy})=\tilde Q(\hat{\yy})\right\},\label{sp-magnanti-wong}
\end{align}
where $\mathcal{U}=\{{\uu}\ge {\vzero}:{\uu}^{\top} A\le {\cc}^{\top}\}$ and $\tilde Q(\hat{\yy})$ is obtained by solving DSP \eqref{generic-milp-projection-dsp}. The constraint ${\uu}^{\top}({\bb}-B \hat{\yy})=\tilde Q(\hat{\yy})$ in \eqref{sp-magnanti-wong} is imposed to guarantee that the dual solution ${\uu}$ is one of the alternative optimal solutions of the DSP induced by $\hat{\yy}$. However, as noted by \cite{papadakos2008practical}, one can still produce a Pareto-optimal cut by suppressing this constraint and instead solving
\begin{align}
    \tilde Q(\bar{\yy})=\max_{\uu\in \mathcal{U}}\;{\uu}^{\top}({\bb}-B\bar{\yy}).\label{sp-papadakos}
\end{align}

Note that $\tilde Q(\bar{\yy})-{\uu}^{\top}({\bb}-B\bar{\yy})\ge 0$ for any ${\uu}\in \mathcal{U}$, and problem \eqref{sp-papadakos} is equivalent to minimizing $\tilde Q(\bar{\yy})-{\uu}^{\top}({\bb}-B\bar{\yy})$. Additionally, $\frac{\BFpi}{\pi_0}\in \mathcal{U}$ for any $({\BFpi},\pi_0)\in\Pi$ such that $\pi_0>0$. Consequently, one can approximate a Pareto-optimal cut when cutting off the point $(\hat{\yy},\hat{\eta})$ by employing
\begin{align*}
    d_{\text{MWP}}(\hat{\yy},\hat{\eta}|{\BFpi},\pi_0)=\frac{{\BFpi}^{\top} ({\bb}-B\hat{\yy})+\pi_0({\ff}^{\top}\hat{\yy}-\hat{\eta})}{\pi_0 \bar{\mu}-{\BFpi}^{\top}({\bb}-B\bar{\yy})}\label{distance-mw}
\end{align*}
as the distance function, in which $\bar{\yy}$ is a core point and $\bar{\mu}>\tilde Q(\bar{\yy})$ ensures that the denominator is positive.
This is equivalent to using the linear normalization function
$h_{\text{MWP}}({\BFpi},\pi_0)={\BFpi}^{\top}{\ww}+\pi_0w_0$ with $({\ww}, w_0)=(B\bar{\yy}-\bb, \bar{\mu})$. We refer to $h_{\text{MWP}}$ as the \textit{Magnanti-Wong-Papadakos (MWP) pseudonorm}, which connects our distance functions to this well-known cut selection strategy.

\begin{remark}
    $\tilde Q(\bar{\yy})$ needs to be computed only once in the course of the BD algorithm (for details, see Algorithm~\ref{pseudo-code-bd-modified} in Section~\ref{sec:modified-bd-algorithm}).
\end{remark}

\begin{remark}
  While the Magnanti-Wong-Papadakos procedure is applicable to optimality cuts, the MWP pseudonorm introduced here is applicable to both optimality and feasibility cuts.
\end{remark}

\subsubsection{Conforti-Wolsey pseudonorm.}\label{sec:cw-norm}
Recently, \cite{conforti2019facet} proposed an interesting procedure for producing a facet-defining cut for separating a point from a convex set. Given a core point $\bar{\yy}$ and $\bar{\eta}>Q(\bar{\yy})=\tilde{Q}(\bar{\yy})+\ff^{\top}\bar{\yy}$ so that $(\bar{\yy},\bar{\eta})$ is in the interior of $\epi$, the geometric interpretation of this idea in our context is to find the closest point to $(\hat{\yy},\hat{\eta})$ on the line segment between $(\bar{\yy},\bar{\eta})$ and $(\hat{\yy},\hat{\eta})$ that renders FSP \eqref{generic-milp-fsp} feasible, which translates to solving
\begin{equation}
    \begin{split}\label{separtion-cw}
        \min\quad & \lambda\\
    \text{s.t.} \quad & -{\cc}^{\top}{\xx} + \lambda\left(\bar{\eta}-\hat{\eta}-\ff^{\top}(\bar{\yy}-\hat{\yy})\right) \ge -\hat{\eta} + {\ff}^{\top}\hat{\yy}\\
    & A{\xx}+\lambda B(\bar{\yy}-\hat{\yy})\ge \bb- B\hat{\yy}\\
    & {\xx} \ge {\vzero}, 1\ge \lambda \ge 0.
    \end{split}
\end{equation}
First, note that we may suppress $\lambda\le 1$ since $(\bar{\yy},\bar{\eta})\in \epi$. Next, assigning dual variable $\pi_0$ to the first constraint and $\BFpi$ to the second set of constraints, we may express LP \eqref{separtion-cw} in its dual form as
\begin{equation*}
    \begin{split}
        \max \quad&{\BFpi}^{\top} ({\bb}-B\hat{\yy})+\pi_0({\ff}^{\top}\hat{\yy}-\hat{\eta})\\
        \text{s.t.}\quad& \BFpi^{\top}B(\bar{\yy}-\hat{\yy}) +\pi_0\left(\bar{\eta}-\hat{\eta}-\ff^{\top}(\bar{\yy}-\hat{\yy})\right)\le 1\\
        &({\BFpi},\pi_0)\in\Pi,
    \end{split}
\end{equation*}
which is equivalent to employing $h_{\text{CW}}({\BFpi},\pi_0)={\BFpi}^{\top}{\ww}+\pi_0w_0$ in  \eqref{normalized-separation-problem} with $w_0= \bar{\eta}-\hat{\eta}-\ff^{\top}(\bar{\yy}-\hat{\yy})$ and $\ww=B(\bar{\yy}-\hat{\yy})$. We refer to $h_{\text{CW}}$ as the \textit{Conforti-Wolsey (CW) pseudonorm}. 
Note that the coefficients of $h_{\text{CW}}$ change as $(\hat{\yy},\hat{\eta})$ changes; thus, unlike other normalization functions presented so far, one should update the normalization constraint for each new point being separated. 
Also, it is not difficult to see that the CW pseudonorm is projective with direction $(\bar{\yy}-\hat{\yy}, \bar{\eta}-\hat{\eta})$, which is in line with the natural geometric interpretation of the CW separation problem \eqref{separtion-cw}. 

\begin{remark}
    From a geometric perspective, CW fixes a target point $(\bar{\yy},\bar{\eta})$ in the interior of $\epi$ and gets as close to this point as possible. In contrast, relaxed $\ell_1$ fixes a direction towards $\epi$ and moves along this direction as much as possible before reaching  $\epi$. Conversely, $d_{\ell_p}$ neither fixes a target point nor a direction; it computes a direction that provides a shortest path to $\epi$ according to $\ell_q$.
\end{remark}

\subsection{Distance-Based Benders Decomposition Algorithm}\label{sec:modified-bd-algorithm}
We present an overview of our proposed Benders decomposition algorithm based on general Benders distance functions in Algorithm~\ref{pseudo-code-bd-modified}. 
Theorem \ref{thrm:complete-Benders-distance-measure}, below, establishes finite convergence of this algorithm for a specific practical class of Benders distance function.

\begin{algorithm}[h]
    \SingleSpacedXI
	\caption{Distance-Based Benders Decomposition Algorithm}
	\label{pseudo-code-bd-modified}
	\begin{algorithmic}[1]
	    \State Select a Benders distance function $d$.
		\State $t \leftarrow 1$, $\hat\Pi_t \leftarrow \emptyset$
		\State Solve MP with $\hat \Pi_t$ in place of $\Pi$ and obtain master solution $({\yy}^{\itt},\eta^{\itt})$.
		\State Solve BSP \eqref{eq:separation_general} to obtain $d^*({\yy}^{\itt},\eta^{\itt})$ and the optimal solution $(\hat{\BFpi},\hat{\pi}_0)$.
		\If{$d^*({\yy}^{\itt},\eta^{\itt})>0$}
		    \State Add a new cut to MP: Set $\hat \Pi_{t+1}\leftarrow \hat \Pi_t\cup\{(\hat{\BFpi},\hat{\pi}_0)\}$, $t\leftarrow t+1$ and loop to Step 3.
		\Else
		    \State Stop. $({\yy}^{\itt},\eta^{\itt})$ is an optimal solution for MP with optimal value $\eta^{\itt}$.
		\EndIf
	\end{algorithmic}
\end{algorithm}

\begin{theorem}\label{thrm:complete-Benders-distance-measure}
Let $d_h$ be a Benders normalized distance function with a convex piece-wise linear normalization function $h$. Then BD Algorithm~\ref{pseudo-code-bd-modified} converges to an optimal solution or asserts infeasibility of MP in a finite number of iterations.
\end{theorem}

Algorithm~\ref{pseudo-code-bd-modified} is finitely convergent when $h$ is a linear function of $({\BFpi},\pi_0)$ (see Section~\ref{sec:linear_pseudonorms} for such linear functions) or when $\ell_1$- or $\ell_{\infty}$-deepest cuts are produced. For other cases (e.g., Euclidean deepest cuts), one may choose to employ Euclidean deepest cuts in conjunction with known finitely convergent separation routines (e.g., classical BD cuts) to guarantee convergence while continuing to benefit from the desirable properties of deepest cuts. 

We present Algorithm~\ref{pseudo-code-bd-modified} for theoretical completeness, but remark that modern implementations of BD add Benders cuts to the cut pool of branch-and-cut using callbacks \citep{fortz2009improved,maher2021implementing}. The resulting algorithm, known as branch-and-Benders-cut, allows for solving the integer master problem in a single run, thus potentially saving computation time by avoiding solving multiple integer master problems. In addition, it allows for separating both integer and fractional master solutions, which may prove useful for producing effective Benders cuts.

\section{Solution Methods for Structured Problems}\label{sec:method_structured}
As shown in Appendix \ref{app:reformulations-lp} and Section \ref{sec:linear_pseudonorms}, generating a distance-based cut with a $\ell_p$ or linear normalization function is no harder than solving a LP/QP. 
Nevertheless, BD is most efficient when there exists an oracle which can exploit the combinatorial properties of the problem instances for solving the subproblems.
For instance, the primal Benders subproblem in several classes of structured problems exhibits useful properties such as generalized upper bounds, which means for a fixed $\yy$, many constraints can simply be treated as variable bounds on the continuous variables. 
Here, we show how we can exploit such combinatorial structures to produce these cuts efficiently. 

\subsection{Directed Depth-Maximizing Algorithm (DDMA) for Linear Normalization Functions}\label{sec:directed-depth-maximizing}

We begin by considering the linear normalization function case (i.e., the normalized separation problem is the MIS subproblem). Note that the normalized separation problem is an LP, which is structurally similar to a classical Benders subproblem, but the addition of the normalization constraint $h(\BFpi, \pi_0)\le 1$ and presence of the variable $\pi_0$ hinder exploiting the combinatorial structures that otherwise could be used in a classical Benders subproblem.
In light of the primal view of the normalized separation problems (Theorem~\ref{thrm:homogeneous-functions-primal}), and noting that $\nabla_{\BFpi}h = \ww$ and $\nabla_{\pi_0}h = w_0$ for $h(\BFpi,\pi_0)=\ww^{\top} \BFpi + w_0\pi_0$, we can progressively increase the cut depth in an iterative manner. Note that fixing $z$ in \eqref{primal-homogeneous} produces a feasibility problem in the space of $\xx$ variables, which shares the constraints $\{A\xx\ge \bb-B\hat{\yy}-z\ww, \xx\ge \vzero\}$ with the classical Benders primal subproblem, and additionally includes the constraint $\cc^{\top}\xx\le \hat{\eta}+z w_0 - \ff^{\top}\hat{\yy}$. Hence, instead of solving this feasibility problem, we can maximize $\cc^{\top}\xx$ over $\{A\xx\ge \bb-B\hat{\yy}-z\ww, \xx\ge \vzero\}$ using the oracle, and check if the constraint $\cc^{\top}\xx\le \hat{\eta}+z w_0 - \ff^{\top}\hat{\yy}$ is satisfied. Therefore, we can consider two cases. 
\begin{itemize}
    \item If $\{A\xx\ge \bb-B\hat{\yy}-z\ww, \xx\ge \vzero\}$ is infeasible, then there exists a Farkas certificate $\bar{\vv}$ such that $\bar{\vv}^{\top}(\bb-B\hat{\yy}-z\ww) > 0$. Therefore, a feasible $z$ must satisfy $z\ge \frac{\bar{\vv}^{\top}(\bb-B\hat{\yy})}{\ww^{\top}\bar{\vv}}$.
    \item Otherwise, there exists an optimal primal-dual solution $(\bar{\xx},\bar{\uu})$ such that $\cc^{\top}\bar{\xx}=\bar{\uu}^{\top}(\bb-B\hat{\yy}-z\ww)$. Therefore, to satisfy the constraint $\cc^{\top}\xx\le \hat{\eta}+z w_0 - \ff^{\top}\hat{\yy}$ we should have $\bar{\uu}^{\top}(\bb-B\hat{\yy}-z\ww)\le \hat{\eta}+z w_0 - \ff^{\top}\hat{\yy}$, which yields $z\ge \frac{\bar{\uu}^{\top}(\bb-B\hat{\yy})+(\ff^{\top}\hat{\yy}-\hat{\eta})}{\ww^{\top}\bar{\uu}+w_0}$.
\end{itemize}
Note that in both cases we compare $z$ to the depth of the current cut with respect to $h$. As we describe in Algorithm \ref{alg:directed-depth-maximizing}, this means that starting with $z=0$, we can iteratively improve the depth of the cut by solving primal subproblems until $z$ becomes feasible (and equals the maximum depth). 
\begin{remark}
    Algorithm \ref{alg:directed-depth-maximizing} produces a sequence of dual solutions in an increasing order of depth with respect to $h$, all of which cut off the incumbent solution $(\hat{\yy},\hat{\eta})$. Consequently, terminating Algorithm \ref{alg:directed-depth-maximizing} before convergence still guarantees convergence of the BD algorithm. 
\end{remark}

\begin{remark}
    For projective linear pseudonorms (i.e., when there exists $\BFomega$ such that $\ww=B\BFomega$) as in relaxed $\ell_1$ and Conforti-Wolsey, Algorithm \ref{alg:directed-depth-maximizing} also produces the $g$-projection of $(\hat{\yy},\hat{\eta})$ onto $\epi$.
\end{remark}

\begin{remark}
    With CW pseudonorm, Algorithm \ref{alg:directed-depth-maximizing} computes the exact step size in the \textit{in-out search} \citep{ben2007acceleration}, as opposed to heuristic methods \citep[e.g.,][]{fischetti2010out}.
\end{remark}

\begin{algorithm}[t]
    \SingleSpacedXI
	\caption{Directed Depth-Maximizing Algorithm (DDMA)}
	\label{alg:directed-depth-maximizing}
	\begin{algorithmic}[1]
        \State \textbf{STEP 0:} Set $z\leftarrow 0$
        \While{not converged}
        \State \textbf{STEP 1 (Cut Generation):} Attempt solving the following primal subproblem:
        \begin{align*}
  [\text{PSP}]\quad\min\left\{\cc^{\top}\xx: A\xx\ge \bb-B\hat{\yy}-z\ww, \xx \ge {\vzero}\right\}.
        \end{align*}
        \State \textbf{if} PSP is feasible \textbf{then} set $(\bar{\BFpi}, \bar{\pi}_0)=(\bar{\uu},1)$ where $\bar{\uu}$ is an optimal dual solution. 
         \State \textbf{else} set $(\bar{\BFpi}, \bar{\pi}_0)=(\bar{\vv},0)$, where $\bar{\vv}$ is a Farkas certificate.
        \State \textbf{STEP 2 (Depth Maximization):} Set $\underline{z}=\frac{\bar{\BFpi}^{\top}(\bb-B\hat{\yy}) +\bar{\pi}_0(\ff^{\top}\hat{\yy}-\hat{\eta})}{\ww^{\top}\bar{\BFpi} + w_0\bar{\pi}_0}$.
        \If{$z = \underline{z}$}
            \State \textbf{Stop}. $(\bar{\BFpi}, \bar{\pi}_0)$ is optimal with maximum depth $d^*_h(\hat{\yy}, \hat{\eta})=z$. 
            \State If there exists $\BFomega$ such that $\ww=B\BFomega$, then output $\tilde{\yy}=\hat{\yy}+z\BFomega$ as the projection of $\hat{\yy}$.
        \EndIf
        
        \State Update $z\leftarrow  \underline{z}$.
        \EndWhile
	\end{algorithmic}
\end{algorithm}

\subsection{Guided Projections Algorithm for Producing \texorpdfstring{$\ell_p$}{Lp}-deepest Cuts}\label{sec:guided_projections_algorithm}

We now describe a procedure for deriving $\ell_p$-deepest cuts for general $p$.
Recall from Theorem~\ref{thrm:deepest-cuts-primal} that producing an $\ell_p$-deepest cut is equivalent to finding the $\ell_q$-projection of the incumbent point $(\hat{\yy},\hat{\eta})$ onto the epigraph $\epi$. Instead of finding this projection directly, we can iteratively \textit{guide} the projection, as illustrated in Figure~\ref{fig:guided-projections}, by moving from the incumbent point to its projection on the epigraph by successively identifying constraints of $\epi$ by calling the oracle; thus, we call this iterative procedure the \textit{Guided Projections Algorithm} (GPA), see Algorithm~\ref{alg:guided-projections} for pseudocode.

\begin{algorithm}[t!]
    \SingleSpacedXI
	\caption{Guided Projections Algorithm (GPA)}
	\label{alg:guided-projections}
	\begin{algorithmic}[1]
	    \State \textbf{STEP 0:} Initialize
	    $t\leftarrow 0$, $\calC^{(0)}\leftarrow\mathbb{R}^{n+1}$, $(\tilde{\yy}^{(0)},\tilde{\eta}^{(0)})\leftarrow (\hat{\yy},\hat{\eta})$.
	   
	    \While{not converged}
  \State \textbf{STEP 1 (Cut Generation):} Attempt solving the following primal subproblem:
  \vspace{-0.2cm}
  \begin{align*}
      [\text{PSP}]\quad\tilde{Q}(\tilde{\yy}^\itt)=\min\left\{\cc^{\top}\xx: A\xx\ge \bb-B\tilde{\yy}^\itt, \xx \ge {\vzero}\right\}.
  \end{align*}
  \vspace{-0.2cm}
   \State \textbf{if} PSP is feasible \textbf{then} set $({\BFpi}^\itt,{\pi}^\itt_0)=(\bar{\uu},1)$ where $\bar{\uu}$ is an optimal dual solution. 
   \State \textbf{else} set $({\BFpi}^\itt,{\pi}^\itt_0)=(\bar{\vv},0)$, where $\bar{\vv}$ is a Farkas certificate.
  \If{$(\tilde{\yy}^\itt,\tilde{\eta}^\itt)\in \cut({\BFpi}^\itt,{\pi}^\itt_0)$}
  \State \textbf{Stop}. $(\tilde{\yy}^\itt,\tilde{\eta}^\itt)$ is the $\ell_q$-projection of $(\hat{\yy},\hat{\eta})$ onto $\epi$
  \EndIf
  \State \textbf{STEP 2 (Projection):} Update $\calC^{(t+1)}\leftarrow\calC^\itt\cap \cut({\BFpi}^\itt,{\pi}^\itt_0)$
  \State Find the $\ell_q$-projection of $(\hat{\yy},\hat{\eta})$ onto $\calC^{(t+1)}$,
  and let  $(\tilde{\yy}^{(t+1)},\tilde{\eta}^{(t+1)})$ be this projection.
  \vspace{-0.2cm}
  \begin{align}
      (\tilde{\yy}^{(t+1)},\tilde{\eta}^{(t+1)})&\leftarrow\agmin_{({\yy},\eta)\in \calC^{(t+1)}}\quad \|({\yy}-\hat{\yy}, \eta-\hat{\eta})\|_{q} \label{gpa-projection}
  \end{align}
  \vspace{-0.2cm}
    \State $t\leftarrow t+1$
	    \EndWhile
	\end{algorithmic}
\end{algorithm}

GPA separates the projection problem \eqref{primal-ssp} into two simpler problems in a row generation manner: (a) tightening the approximation of $\epi$ using classical Benders subproblems, and (b) projecting the incumbent point $(\hat{\yy},\hat{\eta})$ onto the approximation of $\epi$. 
Starting with $(\tilde{\yy}^{(0)},\tilde{\eta}^{(0)})=(\hat{\yy},\hat{\eta})$, GPA produces a classical cut by calling the oracle at $\tilde{\yy}^{(0)}$, and adds the cut to $\calC^{(0)}$ to obtain $\calC^{(1)}$.
GPA then iterates by updating the intermediate projection and calling the oracle.
Given that $d^*_{\ell_p}(\hat{\yy},\hat{\eta})$ is the $\ell_q$ distance from $(\hat{\yy},\hat{\eta})$ to $\epi$, the intermediate projection $(\tilde{\yy}^\itt,\tilde{\eta}^\itt)$ provides both lower- and upper-bounds on $d^*_{\ell_p}(\hat{\yy},\hat{\eta})$ through
$\|(\hat{\yy},\hat{\eta})-(\tilde{\yy}^\itt,\tilde{\eta}^\itt)\|_q$ and $\|(\hat{\yy},\hat{\eta})-(\tilde{\yy}^\itt,Q(\tilde{\yy}^\itt))\|_q$, respectively.
As the algorithm iterates and $\calC^\itt$ becomes a tighter approximation of $\epi$, these bounds converge; thus, the intermediate projections converge to a projection of $(\hat{\yy},\hat{\eta})$ onto $\epi$. 

As GPA iterates, we obtain a sequence of dual solutions, and by construction, each one of them can be used for separating $(\hat{\yy},\hat{\eta})$ from $\epi$ and guaranteeing convergence of BD. Note that the dual solution associated with the deepest cut might not be one of these solutions, but a convex combination of them. Therefore, one can choose to add one or more of the cuts produced by GPA to the BD master problem. 
Figure~\ref{fig:gpa_bounds} in Appendix \ref{app:supplementary_figures} illustrates iterations of GPA and how GPA can replace a classical feasibility cut with a deep optimality cut.

\begin{remark}
    $\calC^{(0)}$ can be initialized with a few simple constraints such as box constraints or the constraints that define $Y$, which in many cases, also correspond to feasibility cuts. 
\end{remark}
\begin{remark}
    Subproblem \eqref{gpa-projection} is an LP for $q\in\{1,\infty\}$ and a QP for $q>1$ and integer, which is expected to be easier than the corresponding separation problem when $\calC^\itt$ contains few constraints.
\end{remark}

\section{Computational Experiments}\label{sec:computational-experiments}

In this section, we illustrate the benefits of deepest cuts and other variants of Benders cuts.
We start by quantifying the benefits of our specialized algorithms for deriving the cuts. We then compare performance of distance-based Benders cuts on several benchmark instances.

\subsection{Experimental Setup}
\subsubsection{Benchmark instances.}
We used four classes of optimization problems from the literature to test our methods. In what follows, we provide a summary of these problems and the benchmark instances for each; detailed descriptions and formulations are in Appendix \ref{app:instances}. 

Our first set of benchmark instances come from the \textit{capacitated facility location problem} (CFLP), whose structure is known to be well-suited for BD \citep{fischetti2016benders}. 
We used two sets of benchmark instances from the literature. The \texttt{CAP} data set from the OR-Library \citep{orlib} and \texttt{CST} instances first introduced in \citep{cornuejols1991comparison}, also known as the \texttt{GK} instances \citep{gortz2012simple}. We considered CFLP under two settings: deterministic demands and stochastic demands. In the former, the subproblem is non-separable, while in the latter, a number of scenarios describe the stochastic demands, resulting in separable subproblems. For the deterministic case, we considered instances with up to $n=1000$ facilities. For the stochastic case, following the recipe of \citet{bodur2016strengthened}, we generated $|K|\in\{254, 512, 1024\}$ demand scenarios based on 16 \texttt{CAP} instances as detailed in Appendix~\ref{app:cflp}, resulting in a total of 48 instances.

We used instances from the \textit{uncapacitated facility location problem} (UFLP) for our second set, for which classical Benders cuts can be derived by solving knapsack problems \citep{fischetti2017redesigning}.
We used 21 instances with up to $n=1000$ facilities from the \texttt{M*} set \citep{kratica2001solving} and 40 instances with up to $n=500$ facilities from the \texttt{KG} set \citep{ghosh2003neighborhood}. See Appendix~\ref{app:uflp} for details.

For the third set, we used instances from the \textit{multicommodity capacitated network design problem} (MCNDP), which widely appears in applications such as telecommunications, transportation, and logistics, and is shown to be computationally challenging \citep{crainic2001bundle}. Given that this problem requires producing both optimality and feasibility cuts, we can evaluate the quality of different Benders cuts in the face of infeasible subproblems. We considered 10 classes of instances from the \texttt{R} set (\texttt{R01}--\texttt{R10}), each consisting of nine instances with different capacities and fixed costs. We also considered a stochastic version of this problem by generating $|K|\in\{16, 32, 64\}$ demand scenarios from the nominal demands as described in Appendix~\ref{app:mcndp}.

Our fourth set of benchmark instances are related to the \textit{stochastic network interdiction problem} (SNIP) proposed by \cite{pan2008minimizing}, and further used in \citep{bodur2016strengthened,boland2016proximity}. As detailed in Appendix~\ref{app:snip}, the set contains five network structures, each with $n=320$ binary first-stage decision variables and $|K|=456$ scenarios.
We considered parameter settings \texttt{snipno=3} and \texttt{snipno=4} with budget ranging in $\{30,40,50\}$ resulting in a total of 30 instances.

\subsubsection{Summary of cut selection routines.}
We implemented Benders decomposition with eight cut selection strategies based on different choices of the normalization function:

\noindent $\boldsymbol{\ell_p}$ \textbf{norms}: Distance functions with $h(\BFpi,\pi_0)=\|\BFpi^{\top}B-\pi_0\ff^{\top}, \pi_0\|_p$ with $p\in\{1,2,\infty\}$. Separation problems are either solved as LP (for $p\in\{1,\infty\}$) or QP (for $p=2$) (see  Appendix~\ref{app:reformulations-lp}), or solved iteratively using the Guided Projections Algorithm (GPA, Algorithm~\ref{alg:guided-projections}).
    
\noindent \textbf{Linear pseudonorms:}  Distance function with linear normalization function $h(\BFpi,\pi_0)=\ww^{\top}\BFpi+w_0\pi_0$. Separation problems are either solved as LP according to the MIS subproblem \eqref{sp-linear-norm-fractional-linear}, or with the Directed Depth-Maximizing Algorithm (DDMA, Algorithm~\ref{alg:directed-depth-maximizing}) iteratively. Special cases are:
\begin{itemize}
    \item \textbf{MISD}: The default choice of $(\ww,w_0)$ in the MIS subproblem as suggested by \cite{fischetti2010note}, that is $w_0=1$, $w_i=0$ if the $i$'th row of $B$ is all zeros and $w_i=1$ otherwise.
    \item $\boldsymbol{R\ell_1}$: $(\ww,w_0)$ chosen according to relaxed $\ell_1$ as described in Section~\ref{sec:relaxed-l1-norm}.
    \item \textbf{MWP}: $(\ww,w_0)$ chosen according to MWP as described in Section~\ref{sec:mw-norm}.
    \item \textbf{CW}: $(\ww,w_0)$ chosen according to CW as described in Section~\ref{sec:cw-norm}.
    \item \textbf{CB}: Classical Benders cuts which correspond to setting $h(\BFpi,\pi_0)=\pi_0$. 
\end{itemize}

\subsubsection{Implementation details.}
For all variants of BD, we employed the following stabilization techniques to ensure an effective implementation. We restricted the separation of fractional solutions to nodes with depth at most 5, and generated Benders cuts after the solver's internal cuts were added (i.e., \texttt{IsAfterCutLoop()==true}). For a fractional solution $(\hat{\yy},\hat{\eta})$, we implemented a simple in-out procedure and separated instead the convex combination of $(\hat{\yy},\hat{\eta})$ and a core point $(\bar{\yy},Q(\bar{\yy}))$ with $0.1$ weight assigned to $(\bar{\yy},Q(\bar{\yy}))$. We also implemented an early termination of the root-node processing as suggested in \cite{bodur2016strengthened} by aborting the cut loop when less than a $0.05\%$ reduction in optimality gap was observed within the last 5 iterations. We considered a cut  violated if its slack was at least $10^{-6}$. As a reference, we also solved the instances as compact MILP models using \texttt{Cplex} with a single thread and other parameters set at their default values. Details regarding coefficient scaling and reoptimization can be found in Appendices~\ref{app:coefficient-scaling} and \ref{app:reoptimization}.

For producing an $\ell_p$-deepest cut, we ran GPA for a maximum of 10 iterations and picked the cut with highest depth according to $\ell_p$. We ran DDMA for a maximum of 10 iterations, or until $|z-\underline{z}|<10^{-6}$, and picked the last cut (which, by construction, is best according to the respective distance function). Note that producing a \texttt{CB} cut corresponds to running DDMA for one iteration.

We computed the gap at the root node via $\frac{U^*-L_0}{U^*}$, where $L_0$ is the lower bound at the root node, and $U^*$ is the best upper bound identified for the problem instance across all methods. 
We used \textit{shifted geometric mean} (SGM) for aggregating the performance results, where $\text{SGM}\{x_i\}_{i=1}^n=\left(\prod_{i=1}^n(x_i+s)\right)^{\frac{1}{n}}-s$, with $s=1$ for gaps and ratios, and $s=10$ for times, number of cuts and branch-and-bound nodes.

\subsection{Comparison of Separation Routines}
As is well-known, BD is most effective when the underlying structural properties of the instances can be exploited \citep{rahmaniani2017benders}. For instance, CFLP and MCNDP exhibit combinatorial structures that can be exploited when solving classical Benders primal subproblems (see Appendices \ref{app:cflp-cut} and \ref{app:mcndp-cut}). In a nutshell, a large number of constraints in both problems can be turned into bounds on the $\xx$ variables when $\yy$ is fixed. Additionally, a cut lifting step can be applied to strengthen the produced cuts by solving a series of continuous knapsack problems. Given the inherent degeneracy of these problems, the latter benefits all cut selection strategies when the solution methods rely on producing classical Benders cuts (including \texttt{CB} itself). 
Therefore, here we report on the effectiveness of different cut selection strategies when the separation problems are (i) solved as general LP/QP's and (ii) solved using our specialized algorithms DDMA (Algorithm~\ref{alg:directed-depth-maximizing}) and GPA (Algorithm~\ref{alg:guided-projections}), which exploit the combinatorial properties of the instances.

Tables \ref{tab:deterministic_solver_rootgap}--\ref{tab:deterministic_solver_cuts} in Appendix~\ref{app:results} present detailed results for moderately-sized instances of CFLP and MCNDP when the separation problems are solved using a solver
(termination criteria: optimality gap of 0.1\% or time limit of 500 seconds). 
Considering the root node gaps and total number of cuts, we note that all distance-based cuts are predominantly more effective compared to \texttt{CB}. For instance, for MCNDP instances, $\ell_1$ achieves a root node gap of 1.75\% compared to 23.71\% for \texttt{CB}. Interestingly, $R\ell_1$ (2\%), \texttt{MWP} (2.59\%) and \texttt{CW} (3.12\%) achieve competitive gaps while outperforming \texttt{MISD} (6.57\%), which highlights the effectiveness of these choices of coefficients in the linear normalization functions compared to the default coefficients (\texttt{MISD}). It is also interesting to note the percentage of cuts that are of feasibility type produced by each method when solving instances of MCNDP (other models only produce optimality cuts).
In particular, $\ell_p$-deepest cuts result in significantly fewer feasibility cuts (1.66\%--4.29\%) compared to \texttt{CB} (28.57\%) and \texttt{MISD} (12.69\%) which also confirms our theoretical insights based on Propositions~\ref{prop:deepst_cuts_optimality} and \ref{prop:deepst_cuts_optimality2}.

Next we demonstrate the effectiveness of our specialized algorithms GPA (for $\ell_p$ deepest cuts) and DDMA (for \texttt{CW}, \texttt{MWP}, $R\ell_1$, and \texttt{CB}). 
As detailed in in
Tables \ref{tab:deterministic_iterative_rootgap}--\ref{tab:deterministic_iterative_cuts} in Appendix~\ref{app:results}, we observe similar trends in terms of relative effectiveness of different choices of cuts, but with a much lower computation cost (around 80\% reduction for distance-based cuts and 50\% reduction for \texttt{CB}). 
Of note, we observe more than a 90\% time reduction for deriving the $\ell_2$-deepest cuts, despite the projection subproblems in GPA being QP. As a result, we observe a significant increase in the number of instances solved to optimality: BD with $\ell_1$-deepest cuts solved all 169 instances, followed by $\ell_{\infty}$ and $\ell_2$ (167 instances), \texttt{CW} (166), $R\ell_1$ and \texttt{MWP} (165) and \texttt{CB} (164).

\subsection{Experiments on CFLP and UFLP}\label{sec:results-nonseparable}

In our second experiment, we focus on large-scale instances of CFLP and UFLP. For CFLP, we considered the 12 largest instances from the \texttt{CAP} set (i.e., \texttt{capa1-4}, \texttt{capb1-4}, and \texttt{capc1-4}) as well as 96 instances from the \texttt{CST} set (with a combination of $\{100, 200, 500, 1000\}$ facilities, $\{500, 1000\}$ customers). We set a time limit of 500 seconds and an optimality gap threshold of $0.1\%$ for all variants of BD and when solving the instances as compact MILP using \texttt{Cplex}.

\begin{figure*}[ht]
    \centering
    \subfigure[CFLP: Gap at the root node.]{
 		\includegraphics[clip,width=0.29\textwidth]{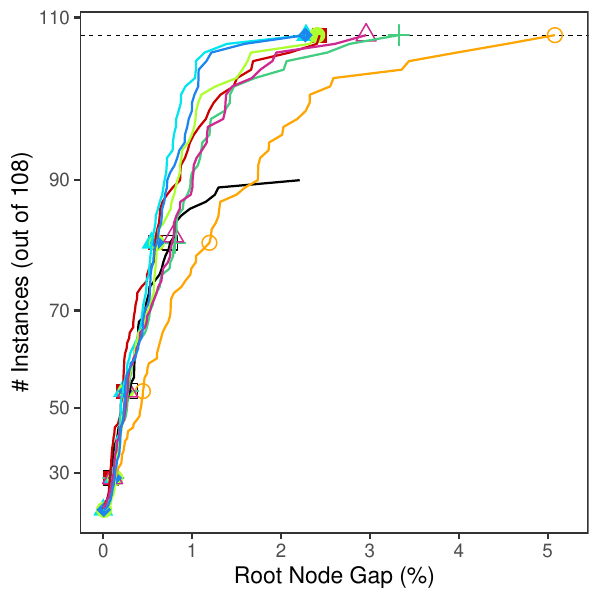}\label{fig:pp-cflp-large-rngap}}
    \subfigure[CFLP: Total number of cuts.]{
 		\includegraphics[clip,width=0.29\textwidth]{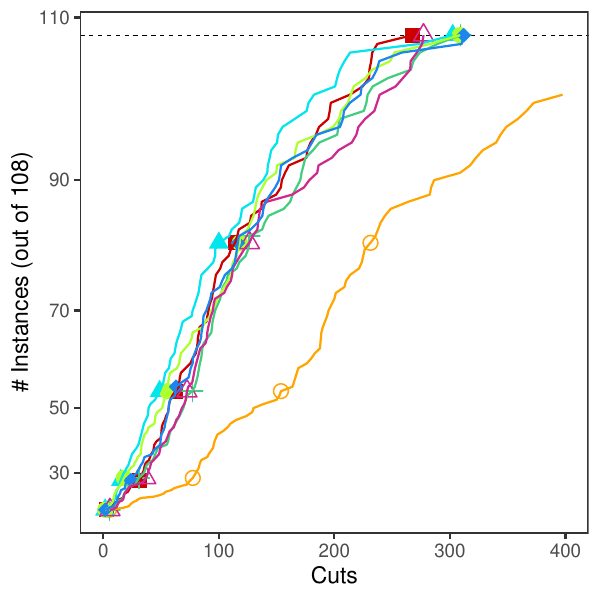}\label{fig:pp-cflp-large-cut}}
    \subfigure[CFLP: Total computing time.]{
 	    \includegraphics[clip,width=0.365\textwidth]{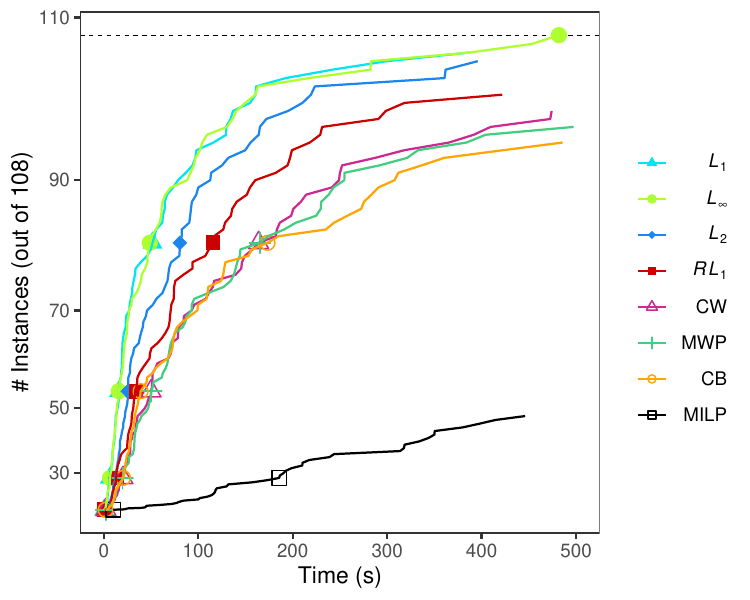}\label{fig:pp-cflp-large-time}}
    \subfigure[UFLP: Gap at the root node.]{
 		\includegraphics[clip,width=0.29\textwidth]{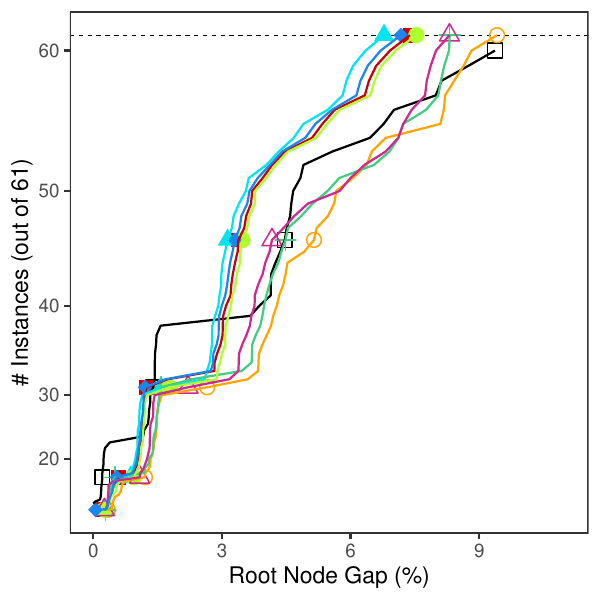}\label{fig:pp-uflp-rngap}}
    \subfigure[UFLP: Total number of cuts.]{
 		\includegraphics[clip,width=0.29\textwidth]{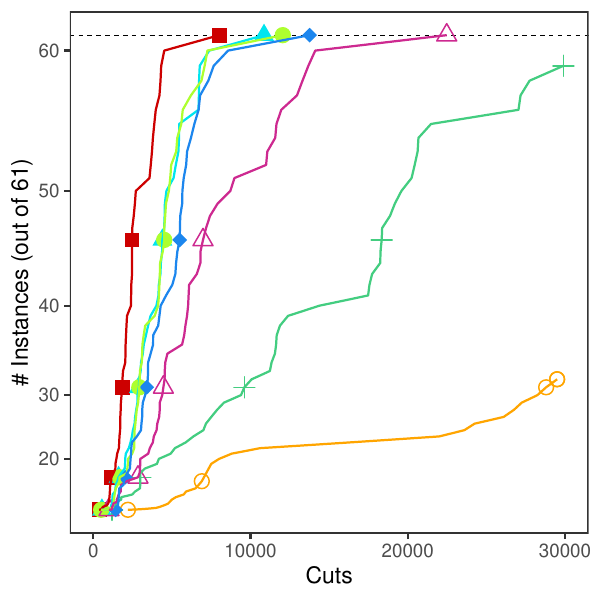}\label{fig:pp-uflp-cut}}
    \subfigure[UFLP: Total computing time.]{
 	    \includegraphics[clip,width=0.365\textwidth]{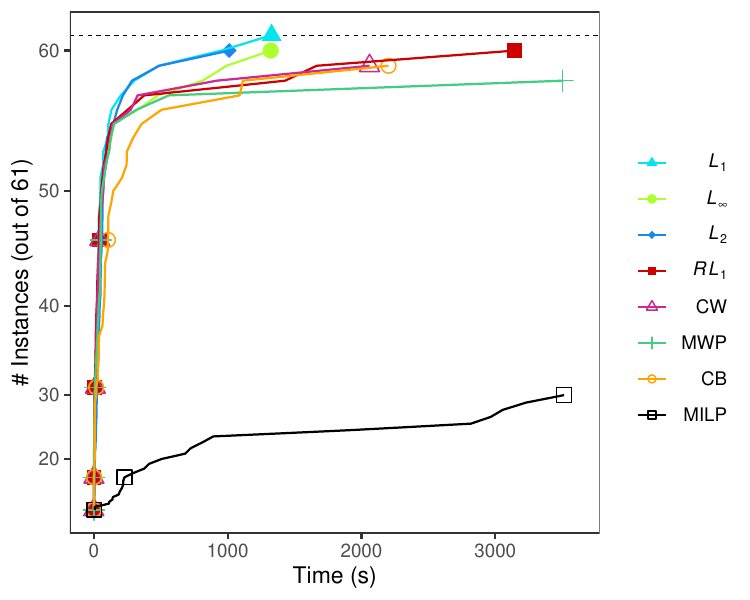}\label{fig:pp-uflp-time}}
    \caption{Comparing the performance of Benders with different cut selection strategies and \texttt{Cplex} (\texttt{MILP}) on large instances of CFLP and UFLP.}\label{fig:pp-flp-large}
\end{figure*}

For ULFP, we considered 21 instances from the \texttt{M*} set with number of facilities/customers ranging in $\{100, 200, 300, 500, 1000\}$ as well as 40 instances from the \texttt{KG} set with $\{250, 500\}$ facilities/customers. We set a time limit of 3600 seconds and, as before an optimality gap threshold of $0.1\%$.
Performance profiles of the models over large-scale instances of CFLP and UFLP are illustrated in Figure~\ref{fig:pp-flp-large}. Note that a line closer to the top left corner indicates better performance.

First, we observe that the distance-based variants of BD are significantly more efficient than \texttt{CB} and the baseline solver (denoted \texttt{MILP}). This is evident from both the computing time of the models and the number of instances solved to optimality. 
More specifically, $\ell_{\infty}$ solved all 108 instances of CFLP, followed closely by $\ell_1$ (106) and $\ell_2$ (105), while $R\ell_1$, \texttt{CW} and \texttt{MWP} solved 101, 99, and 97 instances, respectively. This is while \texttt{CB} solved 95 instances and \texttt{MILP} solved only 48 out of 108 instances (see Figure \ref{fig:pp-cflp-large-time}). 
Similarly, $\ell_{1}$ solved all 61 instances of UFLP, followed by $\ell_{\infty}$, $\ell_2$ and $R\ell_{1}$ with 60 instances each, whereas \texttt{CW}, \texttt{CB} and \texttt{MWP} respectively solved 59, 59, and 58 instances, while \texttt{MILP} solved only 30 out of 61 instances (see Figure \ref{fig:pp-uflp-time}). 
These results highlight effectiveness of GPA and DDMA in deriving the cuts despite the large scale of the separation problems. 
 
Generally, we observe that distance-based cuts, particularly $\ell_1$, require producing fewer cuts and outperform \texttt{CB} and the solver in providing tighter bounds at the root node, which further confirms the effectiveness of deepest cuts in closing the optimality gap, despite the classical cuts in CFLP and UFLP (when derived as in Appendices \ref{app:cflp-cut} and \ref{app:uflp-cut}) often being deemed sufficiently effective in the literature \citep{fischetti2016benders,fischetti2017redesigning}. Specifically, with respect to closing the optimality gap at the root node, we can see that deepest cuts (particularly with $\ell_1$) have a clear edge over both the solver and other variants of BD including \texttt{CW} and \texttt{MWP} (Figures \ref{fig:pp-cflp-large-rngap} and \ref{fig:pp-uflp-rngap}). 
Moreover, the distance-based cuts clearly outperform \texttt{CB} in terms of cut quality (Figures~\ref{fig:pp-cflp-large-cut} and \ref{fig:pp-uflp-cut}).

\subsection{Experiments on Two-Stage Stochastic Instances}

One common application of BD is in tackling two-stage stochastic programs with integer first-stage decision variables.  
Therefore, in our final set of experiments, we study two-stage stochastic programs based on instances from CFLP, MCNDP and SNIP. As discussed earlier, the block-diagonal structure in these problems allows for separating the subproblems upon fixing the first-stage decision variables. We introduce a Stochastic Projections Algorithm (SPA, Algorithm~\ref{alg:sampled-projections} in Appendix~\ref{app:separable_subproblem}) to exploit the repetitive structure of these problems while simultaneously leveraging the geometric properties of distance functions. We use GPA and DDMA as subroutines of SPA for producing the cuts and the projection points.

\begin{figure*}[t]
    \centering
    \subfigure[CFLP: Root node gap.]{
 		\includegraphics[clip,width=0.29\textwidth]{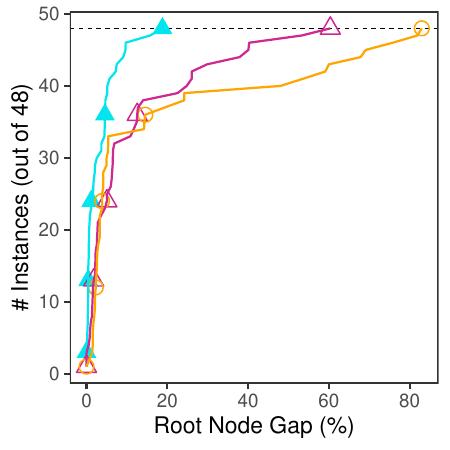}\label{fig:pp-scflp-rngap}}
    \subfigure[CFLP: Total cuts.]{
 		\includegraphics[clip,width=0.29\textwidth]{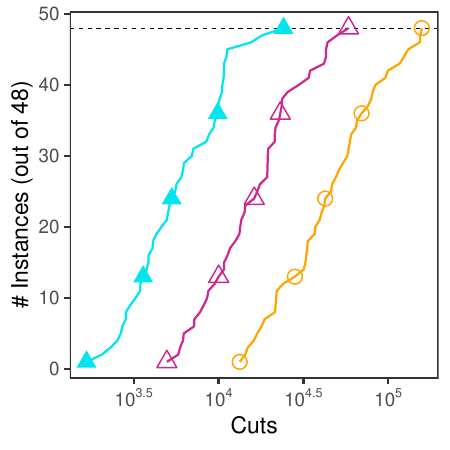}\label{fig:pp-scflp-cut}}
 	\subfigure[CFLP: Total computing time.]{
 	    \includegraphics[clip,width=0.365\textwidth]{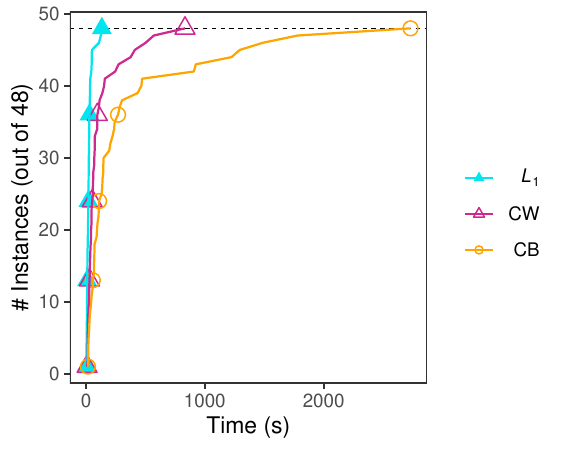}\label{fig:pp-scflp-time}}
      
      \subfigure[MCNDP: Root node gap.]{
 		\includegraphics[clip,width=0.29\textwidth]{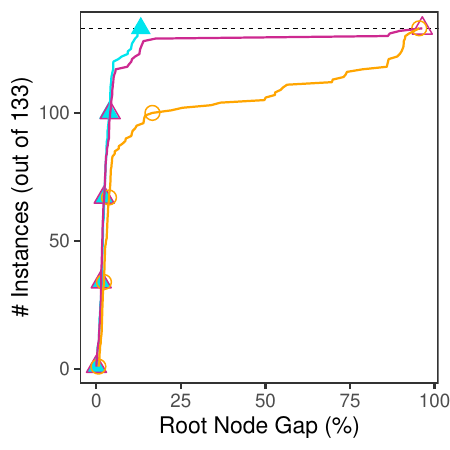}\label{fig:pp-sndp-rngap}}
    \subfigure[MCNDP: Total cuts.]{
 		\includegraphics[clip,width=0.29\textwidth]{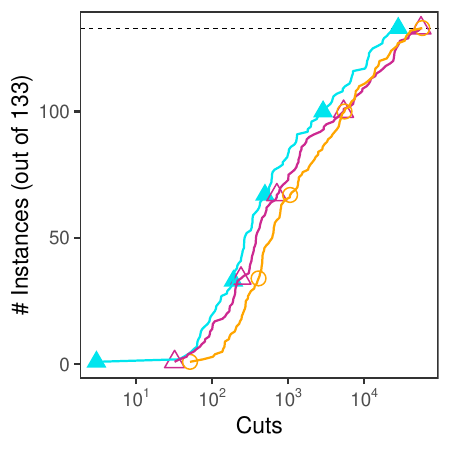}\label{fig:pp-sndp-cut}}
 	\subfigure[MCNDP: Total computing time.]{
 	    \includegraphics[clip,width=0.365\textwidth]{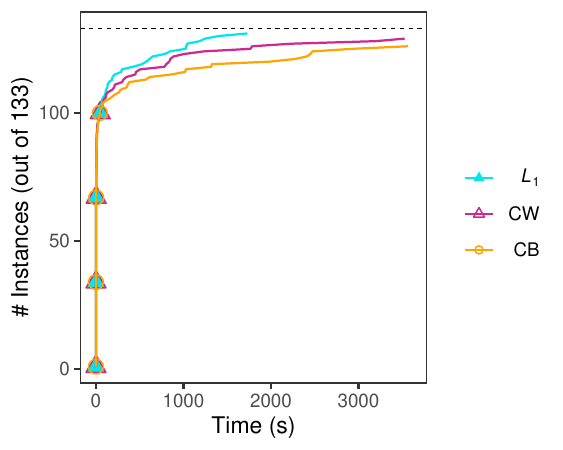}\label{fig:pp-sndp-time}}\hspace{0.2cm}

      \subfigure[SNIP: Root node gap.]{
 		\includegraphics[clip,width=0.29\textwidth]{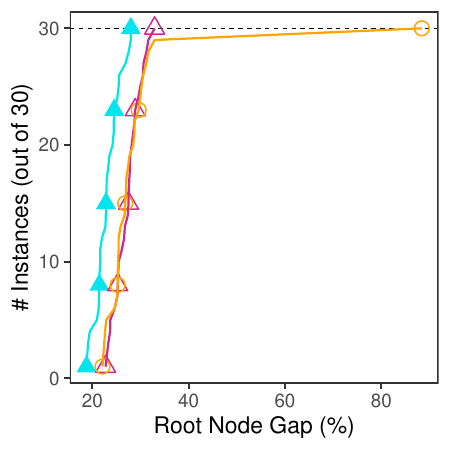}\label{fig:pp-snip-rngap}}
    \subfigure[SNIP: Total cuts.]{
 		\includegraphics[clip,width=0.29\textwidth]{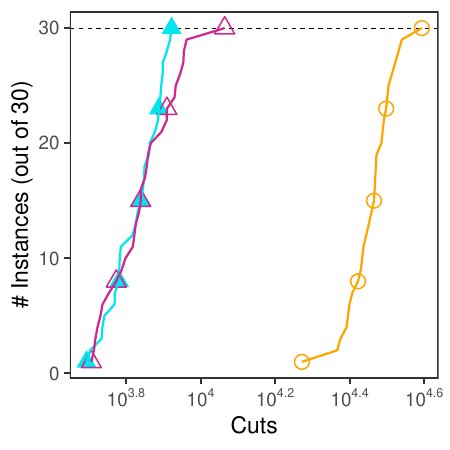}\label{fig:pp-snip-cut}}
 	\subfigure[SNIP: Total computing time.]{
 	    \includegraphics[clip,width=0.365\textwidth]{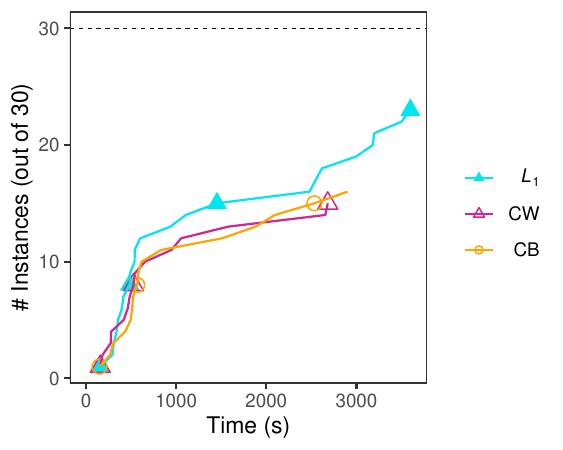}\label{fig:pp-snip-time}}\hspace{0.2cm}
    \caption{Comparing performance of $\ell_1$-deepest cuts, \texttt{CW} cuts and \texttt{CB} cuts on instances of stochastic CFLP (top), stochastic MCNDP (middle), and SNIP (bottom).}\label{fig:pp-stochastic}
\end{figure*}

We report the gaps and computing times at the root node as well as the total computing time and number of B\&B nodes of different variants of BD together with those of \texttt{Cplex} on solving the extended MILPs in Tables \ref{tab:stochastic_rgaptime} and \ref{tab:stochastic_timenodes}, respectively, where time limit and optimality gap are set to 3600 seconds and 0.1\%. We also report the number of cuts in Table \ref{tab:stochastic_cuts}. 
For a more visual comparison, we plot the performance profiles of $\ell_1$ (the top-performing $\ell_p$ deepest cut), \texttt{CW} (a representative of linear pseudonorms) and baseline \texttt{CB} in Figures \ref{fig:pp-stochastic} and \ref{fig:pp-stochastic-supp}. 

All variants of BD efficiently exploit the separability of the instances, outperforming the solver by a significant margin. This margin becomes more significant when the size of the instance increases, as after a point \texttt{MILP} fails to even process the root node. As observed in Table \ref{tab:stochastic_rgaptime}, $\ell_1$-deepest cuts are predominantly more effective in closing the gap at the root node, while requiring essentially the same time compared to other cuts. Effectiveness of these cuts is particularly pronounced in the CFLP instances, where BD with $\ell_1$-deepest cuts is more than three times faster than other variants of BD and orders of magnitude faster than \texttt{MILP}.
We also demonstrate the behavior of $\ell_1$, \texttt{CW} and \texttt{CB} cuts in improving the lower bound as a function of time in Figure~\ref{fig:lowerbounds_plot}, and observe that $\ell_1$ is able to improve the lower bound at a much faster rate for instances of stochastic CFLP.
The SNIP instances turn out to be the most challenging for all variants of BD and \texttt{MILP}, yet $\ell_1$ deepest cuts remain comparatively more effective. In summary, $\ell_1$ leads by solving 202 out of 211 instances to optimality within the allotted time, followed by $\ell_{\infty}$ and $R\ell_1$ (195 each), \texttt{CW} and \texttt{MWP} (192 each), and \texttt{CB} (190), with \texttt{MILP} solving only 148 instances.

\section{Conclusions}\label{sec:conclusions}

In this paper, we proposed and analyzed theoretically and computationally a new method for selecting Benders cuts, aimed at improving the effectiveness of the cuts in closing the gap and reducing the running time of the BD algorithm. Our technique is based on generating Benders cuts that explicitly take cut depth into account. As a measure of cut depth, we first considered Euclidean distance from the master solution to the candidate cuts, and then extended this measure to general $\ell_p$-norms. We provided a comprehensive study of deepest cuts and described 
their properties from a primal perspective. We showed that producing an $\ell_p$-deepest cut is equivalent to finding an $\ell_q$-projection of the point being separated onto the epigraph of the original problem. We also showed how the separation problems can be solved as linear or quadratic programs.

From a theoretical perspective, we generalized our notion of distance by defining what we call a Benders distance function, and developed a notion of monotonicity which allows these functions to be treated as distance functions despite the fact that they do not necessarily satisfy the axioms of metrics.
As an important family of Benders distance functions, we introduced normalized distance functions, and illustrated their connection to some well-known cut selection strategies. Specifically, we established the connection to MIS cuts, and provided three novel ways of choosing the normalization coefficients in the MIS subproblem, that connect our distance functions to the Magnanti-Wong-Papadakos and Conforti-Wolsey procedures.

We showed that a specific class of normalized distance functions, which includes $\ell_p$, relaxed $\ell_1$ and CW, admit a simple characterization of the distance between incumbent point and candidate hyperplane, which is the 1-dimensional distance along the gradient of the normalization function.
In particular, we showed that an intuitive geometric connection exists between CW, our $R\ell_1$ and $\ell_p$ normalization functions: CW fixes a target point inside the epigraph and gets as close to this point as possible. $R\ell_1$ instead fixes a direction toward the epigraph and moves along this direction until a point on the boundary of the epigraph is reached. With $\ell_p$-deepest cuts, we neither fix a target point nor a direction; we find the shortest path to the epigraph.

From a practical perspective, leveraging the duality between $\ell_p$-deepest cuts and $\ell_q$-projections, we introduced our Guided Projections Algorithm for producing $\ell_p$-deepest cuts in a way that can exploit the combinatorial structure of problem instances. By the same token, we introduced our Directed Depth-Maximizing Algorithm for deriving linearly-normalized cuts. Moreover, we showed how the separability and stochastic nature of two-stage stochastic programs can be exploited through our Stochastic Projections Algorithm for deriving distance-based cuts.
    
Our computational experiments on various benchmark problems showed the benefits of deepest cuts and other distance-based cuts, particularly when generated using our tailored algorithms GPA and DDMA, in decreasing the number of cuts as well as the computation time of the BD algorithm. Besides the theoretical insights, our results showed that, distance-based cuts, particularly deepest cuts, are effective in closing the gap at the root node and speeding up the convergence of BD.
In summary, our results suggest that the $\ell_1$ norm is a good choice to generate deep cuts and often works best, while carefully chosen linear alternatives, such as relaxed $\ell_1$ and CW, are often competitive and should also be considered for improving the convergence of BD.

\bibliographystyle{informs2014.bst}
\bibliography{Refs.bib}

\ECSwitch

\begin{APPENDICES}

\section{Proofs of Statements}\label{app:proofs}
In this appendix, we provide the proof of the propositions and theorems given in the body of the paper. For convenience, we formally restate the propositions and theorems as well.

\subsection{Proof of Proposition~\ref{prop:Lq-distance}}\label{app-proof:prop:Lq-distance}

\begin{repeattheorem}[Proposition~\ref{prop:Lq-distance}]
Given $q\ge 1$ and $p\ge 1$ such that $\ell_p$ is the dual norm of $\ell_q$ (i.e., $\frac{1}{p}+\frac{1}{q}=1$), the $\ell_q$-distance from the point $\hat{\zz}\in\mathbb{R}^{n+1}$ to hyperplane ${\BFalpha}^{\top}{\zz}+\beta=0$ is
\begin{align*}
    \min\limits_{\zz: {\BFalpha}^{\top}{\zz}+\beta=0}\|\zz-\hat{\zz}\|_q=\frac{|{\BFalpha}^{\top}\hat{\zz}+\beta|}{\|{\BFalpha}\|_p}.
\end{align*}
\end{repeattheorem}
\prf{For generality, we prove the proposition for general norms using the definition of dual norms; proof for $\ell_p$ norms follows directly.
By definition of dual norms, we have
$$\|\BFalpha\|_*=\max_{\xx}\left\{\frac{|\BFalpha^{\top}\xx|}{\|\xx\|}\right\}.$$
Replacing $\xx=\zz-\hat{\zz}$, we get
\begin{align}
    \|\BFalpha\|_*=\max_{\zz}\left\{\frac{|\BFalpha^{\top}(\zz-\hat{\zz})|}{\|\zz-\hat{\zz}\|}\right\}.\label{eq:general-norm-dist}
\end{align} 
For $\zz\in\mathbb{R}^{n+1}\setminus\{\hat{\zz}\}$, define $\tilde{\zz}(\zz)$ to be the intersection of hyperplane ${\BFalpha}^{\top}{\zz}+\beta=0$ and the line that crosses points $(\zz,\hat{\zz})$. Note that the intersection point for any optimal $\zz$ exists, since the line crossing $(\zz,\hat{\zz})$ cannot be parallel to the hyperplane ${\BFalpha}^{\top}{\zz}+\beta=0$ for optimal $\zz$. This is because a parallel line crossing $(\zz,\hat{\zz})$ and hyperplane ${\BFalpha}^{\top}{\zz}+\beta=0$ would imply that $\BFalpha^{\top}(\zz-\hat{\zz})=0$, which cannot be optimal, since $\|\BFalpha\|_*>0$.  
Now, since $\hat{\zz}$ does not belong to the hyperplane ${\BFalpha}^{\top}{\zz}+\beta=0$, there exists $\theta(\zz)\ne 0$ such that $\zz-\hat{\zz}=\theta(\zz)\times(\tilde{\zz}(\zz)-\hat{\zz})$. We can therefore rewrite \eqref{eq:general-norm-dist} as
\begin{align}
    \|\BFalpha\|_*=\max_{\zz}\left\{\frac{|\theta(\zz)||\BFalpha^{\top}(\tilde{\zz}(\zz)-\hat{\zz})|}{\|\theta(\zz)\times(\tilde{\zz}(\zz)-\hat{\zz})\|}\right\}=\max_{\zz}\left\{\frac{|\BFalpha^{\top}(\tilde{\zz}(\zz)-\hat{\zz})|}{\|(\tilde{\zz}(\zz)-\hat{\zz})\|}\right\},\label{eq:general-norm-dist2}
\end{align}
where the last equality holds since norms are homogeneous. Consequently, without loss of generality we may restrict $\zz$ to the points on the hyperplane ${\BFalpha}^{\top}{\zz}+\beta=0$, that is
\begin{align}
    \|\BFalpha\|_*=\max_{\zz:{\BFalpha}^{\top}{\zz}+\beta=0}\left\{\frac{|\BFalpha^{\top}(\zz-\hat{\zz})|}{\|\zz-\hat{\zz}\|}\right\}=\max_{\zz:{\BFalpha}^{\top}{\zz}+\beta=0}\left\{\frac{|\BFalpha^{\top}\hat{\zz}+\beta|}{\|\zz-\hat{\zz}\|}\right\},\label{eq:general-norm-dist3}
\end{align}
where we have used $\beta=-{\BFalpha}^{\top}{\zz}$. But $|\BFalpha^{\top}\hat{\zz}+\beta|$ is constant, therefore we may rewrite \eqref{eq:general-norm-dist3} as
\begin{align}
    \|\BFalpha\|_*=|\BFalpha^{\top}\hat{\zz}+\beta|\max_{\zz:{\BFalpha}^{\top}{\zz}+\beta=0}\left\{\frac{1}{\|\zz-\hat{\zz}\|}\right\}=\frac{|\BFalpha^{\top}\hat{\zz}+\beta|}{\min\limits_{\zz:{\BFalpha}^{\top}{\zz}+\beta=0}\|\zz-\hat{\zz}\|},\label{eq:general-norm-dist4}
\end{align}
which completes the proof for general norm. The proof for $\ell_q$ follows by replacing $\|\cdot\|=\|\cdot\|_q$ and $\|\cdot\|_*=\|\cdot\|_p$.
}

\subsection{Proof of Theorem~\ref{thrm:deepest-cuts-primal}}\label{app-proof:thrm:deepest-cuts-primal}

\begin{repeattheorem}[Theorem~\ref{thrm:deepest-cuts-primal}.]
Separation problem \eqref{separation-subproblem-Lp} is equivalent to the following Lagrangian dual problem.
\begin{equation}
\begin{split}\label{primal-ssp-proof}
    [\text{Primal SSP}]\quad \min\quad & \|({\yy}-\hat{\yy}, \eta-\hat{\eta})\|_{q}\\
    \text{s.t.} \quad & \eta\ge {\cc}^{\top}{\xx}+{\ff}^{\top}{\yy}\\
    & A{\xx}+B{\yy}\ge {\bb}\\
    & \xx\ge \vzero, \eta \ge \hat{\eta}
\end{split}
\end{equation}
in which $({\yy},{\xx},\eta)$ are the variables and $\ell_q$ is the dual norm of $\ell_p$.
\end{repeattheorem}
\prf{SSP \eqref{separation-subproblem-Lp} can be equivalently stated as (see Proposition~\ref{prop:subproblem-transformation}):
\begin{equation}
    \begin{split}\label{eq:separation-subproblem-Lp-1}
        \max_{({\BFpi},\pi_0)\in\Pi}\quad &{\BFpi}^{\top} ({\bb}-B\hat{\yy})+\pi_0({\ff}^{\top}\hat{\yy}-\hat{\eta})\\
        \text{s.t.}\quad &\|({\BFpi}^{\top} B-\pi_0 {\ff}^{\top},\pi_0)\|_p\le 1.
    \end{split}
\end{equation}
In the following, we prove the statement for $1<p<\infty$, since for $p\in\{1,\infty\}$ the dual can be directly derived using LP duality by reformulating \eqref{eq:separation-subproblem-Lp-1} as an LP (see Appendix~\ref{app:reformulations-lp}). Note that $\|({\BFpi}^{\top} B-\pi_0 {\ff}^{\top},\pi_0)\|_p\le 1$ is equivalent to introducing auxiliary variables $\BFtau$ and setting $\|(\BFtau,\pi_0)\|_p\le 1$ and $-\BFtau\le {\BFpi}^{\top} B-\pi_0 {\ff}^{\top} \le \BFtau$. Hence, we may restate \eqref{eq:separation-subproblem-Lp-1} as
\begin{align}
    \max\quad &{\BFpi}^{\top} ({\bb}-B\hat{\yy})+\pi_0({\ff}^{\top}\hat{\yy}-\hat{\eta})\label{ssp-Lp-obj}\\
    \text{s.t.}\quad 
    & {\BFpi}^{\top}A-\pi_0 {\cc}^{\top}\le 0\label{ssp-Lp-x}\\
    & {\BFpi}^{\top} B-\pi_0 {\ff}^{\top}\le \BFtau^{\top}\label{ssp-Lp-yplus}\\
    & \pi_0 {\ff}^{\top}-{\BFpi}^{\top} B\le \BFtau^{\top}\label{ssp-Lp-yminus}\\
    & \|(\BFtau,\pi_0)\|_p\le 1\label{ssp-Lp-norm}\\
    & {\BFpi}\ge {\vzero}, \pi_0\ge 0. \label{ssp-Lp-sign}
\end{align}
Assigning non-negative Lagrange multipliers ${\xx}$, ${\yy}^+$, ${\yy}^-$, and $z$ respectively to constraints \eqref{ssp-Lp-x}-\eqref{ssp-Lp-norm}, and multipliers $\BFalpha$ and $\gamma$ to the sign constraints \eqref{ssp-Lp-sign}, the Lagrangian function can be stated as
\begin{align}
    &\calL(\xx, \yy^+, \yy^-, z, \BFalpha, \gamma, \BFpi, \pi_0, \BFtau) = {\BFpi}^{\top} ({\bb}-B\hat{\yy})+\pi_0({\ff}^{\top}\hat{\yy}-\hat{\eta}) - ({\BFpi}^{\top}A-\pi_0 {\cc}^{\top})\xx\nonumber\\ 
    & - ({\BFpi}^{\top} B - \pi_0 {\ff}^{\top} - \BFtau^{\top})\yy^+ - (\pi_0 {\ff}^{\top}-{\BFpi}^{\top} B - \BFtau^{\top})\yy^- - (\|(\BFtau,\pi_0)\|_p - 1)z + \BFpi^{\top}\BFalpha + \pi_0\gamma. \label{ssp-Lagrangian}
\end{align}
Therefore, the Lagrangian dual function becomes
\begin{align}
    \calL(\xx, \yy^+, \yy^-, z, \BFalpha, \gamma) = \max_{\BFpi, \pi_0, \BFtau} \calL(\xx, \yy^+, \yy^-, z, \BFalpha, \gamma, \BFpi, \pi_0, \BFtau). \label{ssp-Lagrangian-dual-function}
\end{align}
Taking the derivatives with respect to $\BFpi$, $\pi_0$ and $\BFtau$:
\begin{align}
    \nabla_{\BFpi} \calL = \vzero \Rightarrow & A\xx +B(\hat{\yy} + \yy^+ - \yy^-) = {\bb}+\BFalpha \label{ssp-Lagrangian-derivative-pi}\\
    \frac{\partial\calL}{\partial\pi_0} = 0 \Rightarrow & {\cc}^{\top}\xx + \ff^{\top}(\hat{\yy}+\yy^+ - \yy^-) + \gamma - \hat{\eta}  = z \frac{\pi_0^{p-1}}{\|(\BFtau,\pi_0)\|_p^{p-1}} \label{ssp-Lagrangian-derivative-pi0}\\
     \frac{\partial\calL}{\partial\tau_j} = 0 \Rightarrow & y^+_j+y^-_j = z \frac{\tau_j^{p-1}}{\|(\BFtau,\pi_0)\|_p^{p-1}} \label{ssp-Lagrangian-derivative-tau}
\end{align}
Substituting \eqref{ssp-Lagrangian-derivative-pi}, \eqref{ssp-Lagrangian-derivative-pi0} and \eqref{ssp-Lagrangian-derivative-tau} into \eqref{ssp-Lagrangian-dual-function} we obtain
$$\calL(\xx, z, \yy^+, \yy^-, \BFalpha, \gamma) = z\left(1-\|\BFtau,\pi_0)\|_p+\frac{\pi_0^p}{\|(\BFtau,\pi_0)\|_p^{p-1}}+\sum_{j}\frac{\tau_j^p}{\|(\BFtau,\pi_0)\|_p^{p-1}}\right)=z.$$
Let us define $\eta = {\cc}^{\top}\xx + \ff^{\top}(\hat{\yy}+\yy^+ - \yy^-) + \gamma$. From \eqref{ssp-Lagrangian-derivative-pi0} and \eqref{ssp-Lagrangian-derivative-tau} we obtain:
\begin{align}
    \frac{\eta - \hat{\eta}}{z} = \frac{\pi_0^{p-1}}{\|(\BFtau,\pi_0)\|_p^{p-1}} &\Rightarrow \left(\frac{\eta - \hat{\eta}}{z}\right)^\frac{p}{p-1} = \frac{\pi_0^p}{\|(\BFtau,\pi_0)\|_p^p}\label{ssp-Lagrangian-derivative-pi02}\\
     \frac{y^+_j+y^-_j}{z} = \frac{\tau_j^{p-1}}{\|(\BFtau,\pi_0)\|_p^{p-1}} &\Rightarrow \left(\frac{y^+_j+y^-_j}{z}\right)^\frac{p}{p-1} = \frac{\tau_j^p}{\|(\BFtau,\pi_0)\|_p^p}\label{ssp-Lagrangian-derivative-tau2}
\end{align}
Adding up \eqref{ssp-Lagrangian-derivative-pi02} and \eqref{ssp-Lagrangian-derivative-tau2} and setting $q=\frac{p}{p-1}$ yields
$$\frac{(\eta-\hat{\eta})^q+\sum_{j}(y^+_j+y^-_j)^q}{z^q}=1 \Rightarrow z = \|\yy^+ + \yy^-, \eta-\hat{\eta}\|_q.$$
 Let $\yy = \hat{\yy}+\yy^+ - \yy^-$. Equation \eqref{ssp-Lagrangian-derivative-pi0} implies that $\gamma$ should be sufficiently large so that $\eta-\hat{\eta}\ge 0$. Hence, we may replace $\gamma\ge 0$ with $\eta \ge {\cc}^{\top}\xx + \ff^{\top}\yy$ and $\eta\ge \hat{\eta}$. We can also treat $\BFalpha\ge \vzero$ as simple slack variables and remove them to convert \eqref{ssp-Lagrangian-derivative-pi} to inequality. Minimizing $z$ implies that at the optimal solution, $y^+_j+y^-_j=|y_j-\hat{y}_j|$. 
 Hence, we can state the dual problem as
\begin{equation*}
\begin{split}
    \min\quad & \|({\yy}-\hat{\yy}, \eta-\hat{\eta})\|_{q}\\
    \text{s.t.} \quad & \eta\ge {\cc}^{\top}{\xx}+{\ff}^{\top}{\yy}\\
    & A{\xx}+B{\yy}\ge {\bb}\\
    & \xx\ge \vzero, \eta \ge \hat{\eta}
\end{split}
\end{equation*}
}

\subsection{Proof of Proposition~\ref{prop:deepest-cuts-support}}\label{app-proof:prop:deepest-cuts-support}

\begin{repeattheorem}[Proposition~\ref{prop:deepest-cuts-support}.]
Let $(\tilde{\yy},\tilde{\eta})\in\epi$ be an $\ell_q$-projection of  $(\hat{\yy},\hat{\eta})$ onto $\epi$. Then, any $\ell_p$-deepest cut separating  $(\hat{\yy},\hat{\eta})$ from $\epi$ supports $\epi$ at $(\tilde{\yy},\tilde{\eta})$.    
\end{repeattheorem}
\prf{Let $(\hat{\BFpi},\hat{\pi}_0)$ be the solution associated with the $\ell_p$-deepest cut. By Theorem~\ref{thrm:deepest-cuts-primal} we have 
\begin{align}
    \|(\tilde{\yy}-\hat{\yy},\tilde{\eta}-\hat{\eta})\|_q =\frac{\hat{\BFpi}^{\top}({\bb}-B\hat{\yy})+\hat{\pi}_0({\ff}^{\top}\hat{\yy}-\hat{\eta})}{\|(\hat{\BFpi}^{\top}B-\hat{\pi}_0 {\ff}^{\top},\hat{\pi}_0)\|_p}.\label{eq:deepest-cuts-touch-1}
\end{align}
On the other hand, $(\tilde{\yy},\tilde{\eta})\in\epi$ implies $\hat{\BFpi}^{\top}({\bb}-B\tilde{\yy})+\hat{\pi}_0({\ff}^{\top}\tilde{\yy}-\tilde{\eta})\le 0$. To the contrary, assume that $(\tilde{\yy},\tilde{\eta})$ is not on the hyperplane. Then, $\hat{\BFpi}^{\top}({\bb}-B\tilde{\yy})+\hat{\pi}_0({\ff}^{\top}\tilde{\yy}-\tilde{\eta})$ must be negative, implying
\begin{align}
    0 < -\frac{\hat{\BFpi}^{\top}({\bb}-B\tilde{\yy})+\hat{\pi}_0({\ff}^{\top}\tilde{\yy}-\tilde{\eta})}{\|(\hat{\BFpi}^{\top}B-\hat{\pi}_0 {\ff}^{\top},\hat{\pi}_0)\|_p}.\label{eq:deepest-cuts-touch-2}
\end{align}
Adding \eqref{eq:deepest-cuts-touch-1} and \eqref{eq:deepest-cuts-touch-2} we get
\begin{align*}
    \|(\tilde{\yy}-\hat{\yy},\tilde{\eta}-\hat{\eta})\|_q <\frac{(\hat{\BFpi}^{\top}B-\hat{\pi}_0 {\ff}^{\top}) (\tilde{\yy}-\hat{\yy})+\hat{\pi}_0(\tilde{\eta}-\hat{\eta})}{\|(\hat{\BFpi}^{\top}B-\hat{\pi}_0 {\ff}^{\top},\hat{\pi}_0)\|_p}.
\end{align*}
But this contradicts with H\"{o}lder's inequality since $\ell_p$ and $\ell_q$ are dual norms.}

\subsection{Proof of Proposition~\ref{prop:deepst_cuts_flat}}\label{app-proof:prop:deepst_cuts_flat}
\begin{repeattheorem}[Proposition~\ref{prop:deepst_cuts_flat}.]
For sufficiently small $\hat{\eta}$, the $\ell_1$-deepest cut separating $(\hat{\yy},\hat{\eta})$ from $\epi$ is the flat cut $\eta\ge Q^*$, where $Q^*=\min_{\yy}Q(\yy)$ is the optimal value of $Q$ for unrestricted $\yy$.
\end{repeattheorem}
\prf{Since the dual norm of $\ell_1$ is $\ell_{\infty}$, the objective function in Primal SSP \eqref{primal-ssp} is to minimize the component with largest absolute value in $(\yy-\hat{\yy},\eta-\hat{\eta})$, which,
for sufficiently small $\hat{\eta}$, is $\eta-\hat{\eta}$. Thus, we can restate Primal SSP as the following LP
\begin{equation}
\begin{split}\label{primal-ssp-max}
     -\hat{\eta}+\min\quad & \eta\\
    \text{s.t.} \quad & \eta\ge {\cc}^{\top}{\xx}+{\ff}^{\top}{\yy}\\
    & A{\xx}\ge {\bb}-B{\yy}\\
    & {\xx} \ge {\vzero}.
\end{split}
\end{equation}
Let $(\tilde{\eta},\tilde{\yy},\tilde{\xx})$ be the optimal solution of \eqref{primal-ssp-max}. Observe that $\tilde{\eta}=Q(\tilde{\yy})=\min_{\yy}Q(\yy)$, that is $(\tilde{\eta},\tilde{\yy})$ is an optimal corner point of $\epi$. Further, let $\pi_0$ and $\BFpi$ be the dual multipliers. The dual LP is
\begin{equation}
\begin{split}\label{primal-ssp-max_dual}
     -\hat{\eta}+\max\quad & \BFpi^{\top}\bb\\
    \text{s.t.} \quad & \BFpi^{\top}A\le \pi_0 \cc\\
    &\BFpi^{\top}B=\pi_0 \ff\\
    &\pi_0=1\\
    & {\BFpi} \ge {\vzero}.
\end{split}
\end{equation}
Let $(\hat{\BFpi},\hat{\pi}_0)$ be the optimal solution to \eqref{primal-ssp-max_dual}. The $\ell_1$-deepest cut is $\hat{\BFpi}^{\top}({\bb}-B{\yy})+\hat\pi_0({\ff}^{\top}{\yy}-\eta)=\hat{\BFpi}^{\top}{\bb}-\eta\le 0$. By strong duality, $\hat{\BFpi}^{\top}{\bb}=\tilde{\eta}=Q^*$, hence the deepest cut is the flat cut $\eta\ge Q^*$.
}

\subsection{Proof of Proposition~\ref{prop:deepst_cuts_optimality}}\label{app-proof:prop:deepst_cuts_optimality}
\begin{repeattheorem}[Proposition~\ref{prop:deepst_cuts_optimality}.]
For $p>1$, provided that $\hat{\eta} < Q^*\coloneqq\min\limits_{\yy}Q(\yy)$, the $\ell_p$-deepest cut(s) separating $(\hat{\yy},\hat{\eta})$ are optimality cuts for any arbitrary $\hat{\yy}$ (i.e., even if $\hat{\yy}\notin \dom(Q)$).
\end{repeattheorem}
\prf{Since $\hat{\eta} < Q^*$, we can separate $(\hat{\yy},\hat{\eta})$ from $\epi$ using the flat cut $\bar{\cut}=\{(\yy,\eta):\eta \ge Q^*\}$. Let $(\hat{\BFpi},\hat{\pi}_0)\in \Pi$ be the dual solution associated with the deepest cut, and assume to the contrary that the deepest cut is vertical, that is $\hat{\pi}_0=0$.

Let $(\tilde{\yy}^{H},\tilde{\eta}^{H})$ and $(\tilde{\yy}^{V},\tilde{\eta}^{V})$ be the $\ell_q$-projections of $(\hat{\yy},\hat{\eta})$ onto $\bar{\cut}$ and $\cut(\hat{\BFpi},\hat{\pi}_0)$, respectively. Observe that $(\tilde{\yy}^{H},\tilde{\eta}^{H})=(\hat{\yy}, Q^*)$ and $\tilde{\eta}^{V}=\hat{\eta}$, and that the $\ell_q$-projection of $(\hat{\yy},\hat{\eta})$ onto $\bar{\cut}\cap\cut(\hat{\BFpi},\hat{\pi}_0)$ is $(\tilde{\yy}^{V},\tilde{\eta}^{H})$.
Let $\bar d$ be the $\ell_q$-distance of $(\hat{\yy},\hat{\eta})$ from $\bar{\cut}\cap \cut(\hat{\BFpi},\hat{\pi}_0)$. Note that
\begin{align*}
    \bar d = & \|(\hat{\yy},\hat{\eta})-(\tilde{\yy}^{V},\tilde{\eta}^{H})\|_q = \left(\|\hat{\yy}-\tilde{\yy}^{V}\|_q^q+\|\hat{\eta}-\tilde{\eta}^{H}\|_q^q\right)^{\frac{1}{q}}\\
    =&\left(\|\hat{\yy}-\tilde{\yy}^{V}\|_q^q+\|\hat{\eta}-\tilde{\eta}^{V}\|_q^q+\|\hat{\yy}-\tilde{\yy}^{H}\|_q^q+\|\hat{\eta}-\tilde{\eta}^{H}\|_q^q\right)^{\frac{1}{q}}\\
    =&\left(\|(\hat{\yy},\hat{\eta})-(\tilde{\yy}^{V},\tilde{\eta}^{V})\|_q^q+\|(\hat{\yy},\hat{\eta})-(\tilde{\yy}^{H},\tilde{\eta}^{H})\|_q^q\right)^{\frac{1}{q}}\\
    =&\left((d^*)^q+(Q^*-\hat{\eta})^q\right)^{\frac{1}{q}},
\end{align*}
where $d^*=\|(\hat{\yy},\hat{\eta})-(\tilde{\yy}^{V},\tilde{\eta}^{V})\|_q$. This implies that $\bar d>d^*$ since $q<\infty$ and $Q^*>\hat{\eta}$.
However, given that both $\bar{\cut}$ and $\cut(\hat{\BFpi},\hat{\pi}_0)$ support $\epi$, the $\ell_q$-distance from $(\hat{\yy},\hat{\eta})$ to $\epi$ (i.e., $d^*$) must be at least equal to $\bar d$, that is $d^*\ge \bar d$, which is a contradiction.
}

\subsection{Proof of Proposition~\ref{prop:deepst_cuts_optimality2}}\label{app-proof:prop:deepst_cuts_optimality2}
\begin{repeattheorem}[Proposition~\ref{prop:deepst_cuts_optimality2}.]
Provided that the $\ell_q$-projection of $(\hat{\yy},\hat{\eta})$ onto $\epi$ is the unique point $(\tilde{\yy},\tilde{\eta})$ and $\hat{\eta} < \tilde{\eta}$, the $\ell_p$-deepest cuts separating $(\hat{\yy},\hat{\eta})$ are optimality cuts even if $\hat{\yy}\notin \dom(Q)$.
\end{repeattheorem}
\prf{
Let $(\hat{\BFpi},\hat{\pi}_0)\in \Pi$ be the dual solution associated with the deepest cut.
Assume to the contrary that the deepest cut is a feasibility cut, that is $\hat{\pi}_0=0$. Since the projection is unique, the $\ell_q$-projection of $(\hat{\yy},\hat{\eta})$ onto the vertical cut $\cut(\hat{\BFpi},\hat{\pi}_0)$ must be $(\tilde{\yy},\hat{\eta})$, which contradicts with the assumption that $\tilde{\eta}>\hat{\eta}$. 
}

\subsection{Proof of Proposition~\ref{prop:subproblem-transformation}}\label{app-proof:prop:subproblem-transformation}

\begin{repeattheorem}[Proposition~\ref{prop:subproblem-transformation}.]
Let $h$ be a normalization function, $d_h(\hat{\yy},\hat{\eta}|{\BFpi},\pi_0) = \frac{{\BFpi}^{\top} ({\bb}-B\hat{\yy})+\pi_0({\ff}^{\top}\hat{\yy}-\hat{\eta})}{h({\BFpi},\pi_0)}$ the distance function induced by $h$, and $\Pi_h=\{({\BFpi},\pi_0)\in\Pi:h({\BFpi},\pi_0)\le 1\}$ the cone $\Pi$ truncated by the constraint $h({\BFpi},\pi_0)\le 1$.
Then, the separation problem \eqref{eq:separation_general} is equivalent to the normalized separation problem (NSP) defined below, and $h({\BFpi},\pi_0)\le 1$ is binding at optimality.
\begin{align}
     [\text{NSP}]\quad d^*_h(\hat{\yy},\hat{\eta})=\max_{({\BFpi},\pi_0)\in\Pi_h} \;&{\BFpi}^{\top} ({\bb}-B\hat{\yy})+\pi_0({\ff}^{\top}\hat{\yy}-\hat{\eta}) \label{app:normalized-separation-problem}
\end{align}
\end{repeattheorem}
\prf{The separation problem \eqref{eq:separation_general} can be equivalently expressed as
\begin{align}
   \max_{q>0}\left\{\max_{(\tilde{\BFpi},\tilde{\pi}_0)\in\Pi: h(\tilde{\BFpi},\tilde{\pi}_0)=q}\frac{\tilde{\BFpi}^{\top} ({\bb}-B\hat{\yy})+\tilde{\pi}_0({\ff}^{\top}\hat{\yy}-\hat{\eta})}{q}\right\}.\label{generic-subproblem-transformation-obj}
\end{align}
Define $({\BFpi},\pi_0)= \frac{1}{q}(\tilde{\BFpi},\tilde{\pi}_0)$. Since $\Pi$ is a cone it follows that $({\BFpi},\pi_0)\in\Pi$. Additionally, since $h$ is homogeneous, we have $h({\BFpi},\pi_0)=\frac{1}{q}h(\tilde{\BFpi},\tilde{\pi}_0)=1$. Therefore, the inner maximization in \eqref{generic-subproblem-transformation-obj} can be restated as
\begin{align}
   \max_{({\BFpi},\pi_0)\in\Pi: h({\BFpi},\pi_0)=1}{\BFpi}^{\top} ({\bb}-B\hat{\yy})+\pi_0({\ff}^{\top}\hat{\yy}-\hat{\eta}),\label{generic-subproblem-transformation-obj2}
\end{align}
which is constant with respect to $q$. Therefore, \eqref{generic-subproblem-transformation-obj} itself is equivalent to \eqref{generic-subproblem-transformation-obj2}. 
We next show that \eqref{generic-subproblem-transformation-obj2} is equivalent to \eqref{app:normalized-separation-problem}, that is $h({\BFpi},\pi_0)=1$ can be replaced with $h({\BFpi},\pi_0)\le 1$. Let $(\hat{\BFpi},\hat{\pi}_0)\in\Pi_h$ be an arbitrary solution to \eqref{app:normalized-separation-problem} with $\hat\alpha=h(\hat{\BFpi},\hat{\pi}_0)<1$. Note that $(\bar{\BFpi},\bar{\pi}_0)=\frac{1}{\hat\alpha}(\hat{\BFpi},\hat{\pi}_0)\in \Pi_h$ with $h(\bar{\BFpi},\bar{\pi}_0)=1$. Additionally, we have
$$\bar{\BFpi}^{\top} ({\bb}-B\hat{\yy})+\bar{\pi}_0({\ff}^{\top}\hat{\yy}-\hat{\eta})=\frac{\hat{\BFpi}^{\top} ({\bb}-B\hat{\yy})+\hat{\pi}_0({\ff}^{\top}\hat{\yy}-\hat{\eta})}{\hat{\alpha}}\ge \hat{\BFpi}^{\top} ({\bb}-B\hat{\yy})+\hat{\pi}_0({\ff}^{\top}\hat{\yy}-\hat{\eta}),$$
which is strict if $\hat{\BFpi}^{\top} ({\bb}-B\hat{\yy})+\hat{\pi}_0({\ff}^{\top}\hat{\yy}-\hat{\eta})>0$. Thus, at optimality $h({\BFpi},\pi_0)\le 1$ is binding.}

\subsection{Proof of Theorem~\ref{thrm:homogeneous-functions-primal}}\label{app-proof:thrm:homogeneous-functions-primal}

\begin{repeattheorem}[Theorem~\ref{thrm:homogeneous-functions-primal}.]
Let $h(\BFpi,\pi_0)$ be a convex differentiable positive homogeneous function. Assuming that NSP \eqref{normalized-separation-problem} admits a bounded optimal solution $(\BFpi^*, \pi_0^*)$, NSP is equivalent to the following LP
\begin{equation}\label{prf-primal-homogeneous}
\begin{split}
     d^*_h(\hat{\yy},\hat{\eta})=\min\quad & z\\
    \text{s.t.} \quad & A\xx \ge \bb-B\hat{\yy}-z\nabla_{\BFpi}h(\BFpi^*,\pi^*_0)\\
    & \cc^{\top}\xx\le \hat{\eta}-\ff^{\top}\hat{\yy}+z\nabla_{\pi_0}h(\BFpi^*,\pi^*_0)\\
    & \xx\ge \vzero
\end{split}
\end{equation}
\end{repeattheorem}
\prf{Since $h$ is homogeneous of degree 1, by Euler's homogeneous function theorem, we have
\begin{align*}
h(\BFpi,\pi_0)=\BFpi^{\top}\nabla_{\BFpi}h(\BFpi,\pi_0)+\pi_0\nabla_{\pi_0}h(\BFpi,\pi_0).
\end{align*}
Define $\tilde{h}(\BFpi,\pi_0)=\BFpi^{\top}\nabla_{\BFpi}h(\BFpi^*,\pi^*_0)+\pi_0\nabla_{\pi_0}h(\BFpi^*,\pi^*_0))$. Note that replacing $h$ with $\tilde{h}$ in NSP \eqref{normalized-separation-problem} results in the same objective value $d^*_h(\hat{\yy},\hat{\eta})$.
Consequently, we can restate NSP \eqref{normalized-separation-problem} as 
\begin{equation*}
\begin{split}
    \max \quad &{\BFpi}^{\top} ({\bb}-B\hat{\yy})+\pi_0({\ff}^{\top}\hat{\yy}-\hat{\eta})\\
    \text{s.t.}\quad & \BFpi^{\top} A - \pi_0 \cc \le \vzero\\
    & \BFpi^{\top}\nabla_{\BFpi}h(\BFpi^*,\pi^*_0)+\pi_0\nabla_{\pi_0}h(\BFpi^*,\pi^*_0))\le 1\\
    &(\BFpi,\pi_0)\ge \vzero
\end{split}    
\end{equation*}
which is an LP. Taking the dual of this LP yields the result.
}

\subsection{Proof of Proposition~\ref{prop-projective-functions-distance-hyperplane}}\label{app-proof:prop-projective-functions-distance-hyperplane}

\begin{repeattheorem}[Proposition~\ref{prop-projective-functions-distance-hyperplane}.]
    Let $g$ be a projective normalization function. For any $(\hat{\yy}, \hat{\eta})$ and any $(\BFpi,\pi_0)\in\Pi$ such that $(\BFtau, \pi_0)\ne\vzero$ with $\BFtau=B^{\top}\BFpi-\pi_0 \ff$, the line $(\yy, \eta) = (\hat{\yy},\hat{\eta})+ z\nabla_{(\BFtau, \pi_0)}g(\BFtau, \pi_0)$ intersects the hyperplane $\hp(\BFpi,\pi_0)$ at $(\tilde{\yy}, \tilde{\eta})$ for $z=d_{g}(\hat{\yy},\hat{\eta}|\BFpi,\pi_0)$, thus implicitly defining $d_g$.
\end{repeattheorem}
    
\prf{
Consider the line $(\yy,\eta)=(\hat{\yy},\hat{\eta})+z\nabla g(\BFtau, \pi_0)$. For this line to intersect the hyperplane $\BFpi^{\top}(\bb-B{\yy})+\pi_0(\ff^{\top}{\yy}-{\eta})=0$ at a point $(\tilde{\yy},\tilde{\eta})$ we must have
\begin{align*}
\BFpi^{\top}(\bb-B\tilde{\yy})+\pi_0(\ff^{\top}\tilde{\yy}-\tilde{\eta})=0 \text{ and } (\tilde{\yy},\tilde{\eta})=(\hat{\yy},\hat{\eta})+z\nabla g(\BFtau, \pi_0).
\end{align*}
Therefore, $z$ must satisfy
\begin{align*}
z \left((\BFpi^{\top}B-\pi_0 \ff^{\top})\nabla_{\BFtau} g(\BFtau, \pi_0)+\pi_0 \nabla_{\pi_0} g(\BFtau, \pi_0)\right) = \BFpi^{\top}(\bb-B\hat{\yy})+\pi_0(\ff^{\top}\hat{\yy}-\hat{\eta})
\end{align*}
Replacing $\BFtau=B^{\top}\BFpi-\pi_0 \ff$ and using Euler's homogeneous function theorem, we obtain
\begin{align*}
    z \left((\BFpi^{\top}B-\pi_0 \ff^{\top})\nabla_{\BFtau} g(\BFtau, \pi_0)+\pi_0 \nabla_{\pi_0} g(\BFtau, \pi_0)\right) = z\left(\BFtau^{\top}\nabla_{\BFtau} g(\BFtau, \pi_0)+\pi_0 \nabla_{\pi_0} g(\BFtau, \pi_0)\right) = zg(\BFtau,\pi_0)    
\end{align*}
However, by definition of projective normalization function, $g(\BFtau,\pi_0)$ must be strictly positive since $(\BFtau, \pi_0)\ne\vzero$. Consequently, the intersection point $(\tilde{\yy},\tilde{\eta})$ exists and we obtain
\begin{align*}
z = \frac{\BFpi^{\top}(\bb-B\hat{\yy})+\pi_0(\ff^{\top}\hat{\yy}-\hat{\eta})}{g(\BFtau,\pi_0)} = d_{g}(\hat{\yy},\hat{\eta}|\BFpi,\pi_0).
\end{align*}
}

\subsection{Proof of Proposition~\ref{prop:projective-functions-bound_lp}}\label{app-proof:prop:projective-functions-bound_lp}

\begin{repeattheorem}[Proposition~\ref{prop:projective-functions-bound_lp}.]
Let $g$ be a projective normalization function. 
For any $p\ge 1$ and its dual $q$ (i.e., $\frac{1}{p}+\frac{1}{q}=1$), any $(\hat{\yy}, \hat{\eta})$ and any $(\BFpi,\pi_0)\in\Pi$ with $\BFtau=B^{\top}\BFpi-\pi_0 \ff$ such that $(\hat{\yy}, \hat{\eta})\notin \cut(\BFpi, \pi_0)$ we have
\begin{align*}
    d_{\ell_p}(\hat{\yy},\hat{\eta}|\BFpi,\pi_0)\le d_{g}(\hat{\yy},\hat{\eta}|\BFpi,\pi_0)\; \|\nabla g(\BFtau, \pi_0)\|_q.
\end{align*}
\end{repeattheorem}
\prf{The proof follows similarly as in the proof of Proposition~\ref{prop-projective-functions-distance-hyperplane}. Note that when the line $(\yy,\eta)=(\hat{\yy},\hat{\eta})+z\nabla g(\BFtau, \pi_0)$ intersects the hyperplane $\BFpi^{\top}(\bb-B{\yy})+\pi_0(\ff^{\top}{\yy}-{\eta})=0$ at a point $(\tilde{\yy},\tilde{\eta})$, the $\ell_q$ distance from $(\hat{\yy},\hat{\eta})$ to $(\tilde{\yy},\tilde{\eta})$ is
\begin{align*}
    \|(\hat{\yy},\hat{\eta})-(\tilde{\yy},\tilde{\eta})\|_q = z\|\nabla g(\BFtau, \pi_0)\| = d_{g}(\hat{\yy},\hat{\eta}|\BFpi,\pi_0)\; \|\nabla g(\BFtau, \pi_0)\|_q.
\end{align*}
Consequently, given that $d_{\ell_p}(\hat{\yy},\hat{\eta}|\BFpi,\pi_0)$ measures the $\ell_q$ distance from  $(\hat{\yy},\hat{\eta})$ to the hyperplane $\BFpi^{\top}(\bb-B{\yy})+\pi_0(\ff^{\top}{\yy}-{\eta})=0$, and $(\tilde{\yy},\tilde{\eta})$ lies on this hyperplane, we obtain
\begin{align*}
    d_{\ell_p}(\hat{\yy},\hat{\eta}|\BFpi,\pi_0)\le d_{g}(\hat{\yy},\hat{\eta}|\BFpi,\pi_0)\; \|\nabla g(\BFtau, \pi_0)\|_q.
\end{align*}
When the line $(\yy,\eta)=(\hat{\yy},\hat{\eta})+z\nabla g(\BFtau, \pi_0)$ and the hyperplane $\BFpi^{\top}(\bb-B{\yy})+\pi_0(\ff^{\top}{\yy}-{\eta})=0$ do not intersect, we have $d_{g}(\hat{\yy},\hat{\eta}|\BFpi,\pi_0) = \infty$ and the claim follows since $d_{\ell_p}(\hat{\yy},\hat{\eta}|\BFpi,\pi_0) < \infty$ and $\|\nabla g(\BFtau, \pi_0)\|_q > 0$.
}

\subsection{Proof of Proposition~\ref{prop:projective-functions-primal}}\label{app-proof:prop:projective-functions-primal}
\begin{repeattheorem}[Proposition~\ref{prop:projective-functions-primal}.]
Let $g$ be a projective normalization function and $(\BFpi^*, \pi_0^*)$ an optimal solution to the separation problem NSP \eqref{normalized-separation-problem} with $h(\BFpi,\pi_0)=g(\BFtau,\pi_0)$ where $\BFtau=B^{\top}\BFpi-\pi_0 \ff$. Then
    \begin{equation}\label{app-primal-projective}
    \begin{split}
         d^*_g(\hat{\yy},\hat{\eta})=\min \; \{z: (\yy,\eta)\in \epi,\quad  (\yy, \eta) = (\hat{\yy},\hat{\eta})+ z\nabla_{(\BFtau, \pi_0)}g(\BFtau^*, \pi_0^*)\}.
    \end{split}
    \end{equation}
\end{repeattheorem}
\prf{Given $h(\BFpi,\pi_0)=g(\BFtau,\pi_0) = g(\BFpi^{\top}B-\pi_0 \ff, \pi_0)$, by chain rule, we have
\begin{align*}
    \nabla_{\BFpi}h &= B\nabla_{\BFtau}g\\
    \nabla_{\pi_0}h &= -\ff^{\top}\nabla_{\BFtau}g+\nabla_{\pi_0}g
\end{align*}
Consequently, by evaluating Theorem~\ref{thrm:homogeneous-functions-primal} at $h$ we deduce
\begin{equation*}
\begin{split}
     d^*_h(\hat{\yy},\hat{\eta})=\min\quad & z\\
    \text{s.t.} \quad & A\xx \ge \bb-B\hat{\yy}-zB\nabla_{\BFtau}g(\BFtau^*,\pi^*_0)\\
    & \cc^{\top}\xx\le \hat{\eta}-\ff^{\top}\hat{\yy}-z\ff^{\top}\nabla_{\BFtau}g(\BFtau^*,\pi^*_0)+z\nabla_{\pi_0}g(\BFtau^*,\pi^*_0)\\
    & \xx\ge \vzero
\end{split}
\end{equation*}
Setting $\yy=\hat{\yy}+z\nabla_{\BFtau}g(\BFtau^*,\pi^*_0)$ and $\eta=\hat{\eta}+z\nabla_{\pi_0}g(\BFtau^*,\pi^*_0)$ we obtain
\begin{equation*}
\begin{split}
     d^*_g(\hat{\yy},\hat{\eta})=\min\quad & z\\
    \text{s.t.} \quad & A\xx + B \yy\ge b\\
    & \cc^{\top}\xx +\ff^{\top}\yy\le \eta\\
    & \xx\ge \vzero\\
    & (\yy, \eta) = (\hat{\yy},\hat{\eta})+ z\nabla_{(\BFtau, \pi_0)}g(\BFtau^*, \pi_0^*)
\end{split}
\end{equation*}
Replacing the first three sets of constraints with $(\yy,\eta)\in \epi$ gives the result.
}

\subsection{Proof of Proposition~\ref{prop:depth-gap-lp}}\label{app-proof:prop:depth-gap-lp}

\begin{repeattheorem}[Proposition~\ref{prop:depth-gap-lp}.]
The following relationship  holds between $d_{\text{CB}}$, $d_{\ell_p}$, and $d_{R\ell_1}$ for any $(\BFpi,\pi_0)\in \Pi$:
\begin{align*}
    d_{\text{CB}}(\hat{\yy},\hat{\eta}|\BFpi,\pi_0)\ge d_{\ell_{\infty}}(\hat{\yy},\hat{\eta}|\BFpi,\pi_0)\ge\dots\ge d_{\ell_p}(\hat{\yy},\hat{\eta}|\BFpi,\pi_0)\ge \dots\ge d_{\ell_1}(\hat{\yy},\hat{\eta}|\BFpi,\pi_0)\ge d_{R\ell_1}(\hat{\yy},\hat{\eta}|\BFpi,\pi_0).
\end{align*}
\end{repeattheorem}
\prf{Recall that the $d_{\ell_p}$, $d_{R\ell1}$, and $d_{\text{CB}}$ distance functions are defined as
\begin{align*}
    d_{\ell_p}(\hat{\yy},\hat{\eta}|{\BFpi},\pi_0)=&\frac{{\BFpi}^{\top} ({\bb}-B\hat{\yy})+\pi_0({\ff}^{\top}\hat{\yy}-\hat{\eta})}{\|({\BFpi}^{\top} B-\pi_0 {\ff}^{\top},\pi_0)\|_p}\\
    d_{R\ell1}(\hat{\yy},\hat{\eta}|{\BFpi},\pi_0)=&\frac{{\BFpi}^{\top} ({\bb}-B\hat{\yy})+\pi_0({\ff}^{\top}\hat{\yy}-\hat{\eta})}{\sum_{i=1}^m\pi_i\sum_{j=1}^n|B_{ij}|+(1+\sum_{j=1}^n|f_j|)\pi_0}\\
    d_{\text{CB}}(\hat{\yy},\hat{\eta}|{\BFpi},\pi_0)=&\frac{{\BFpi}^{\top} ({\bb}-B\hat{\yy})+\pi_0({\ff}^{\top}\hat{\yy}-\hat{\eta})}{\pi_0}
\end{align*}
The proof follows by noting the following facts:
\begin{enumerate}[label=(\roman*)]
    \item $\|\cdot\|_p\le \|\cdot\|_{p'}$ for any $1\le p'<p$. This implies that
    $$d_{\ell_{\infty}}(\hat{\yy},\hat{\eta}|{\BFpi},\pi_0)\ge\dots\ge d_{\ell_p}(\hat{\yy},\hat{\eta}|{\BFpi},\pi_0)\ge \dots\ge d_{\ell_1}(\hat{\yy},\hat{\eta}|{\BFpi},\pi_0)\qquad \forall (\BFpi,\pi_0)\in \Pi.$$
    
    \item $\|({\BFpi}^{\top} B-\pi_0 {\ff}^{\top},\pi_0)\|_p\ge \pi_0$ for any $p\ge 1$. This implies that
    $$d_{\text{CB}}(\hat{\yy},\hat{\eta}|{\BFpi},\pi_0)\ge d_{\ell_{\infty}}(\hat{\yy},\hat{\eta}|{\BFpi},\pi_0)\qquad \forall (\BFpi,\pi_0)\in \Pi.$$
    
    \item $\|({\BFpi}^{\top} B-\pi_0 {\ff}^{\top},\pi_0)\|_1\le \sum_{i=1}^m\pi_i\sum_{j=1}^n|B_{ij}|+(1+\sum_{j=1}^n|f_j|)\pi_0$. This implies that
    $$d_{\ell_1}(\hat{\yy},\hat{\eta}|{\BFpi},\pi_0)\ge d_{R\ell_1}(\hat{\yy},\hat{\eta}|{\BFpi},\pi_0)\qquad \forall (\BFpi,\pi_0)\in \Pi.$$
\end{enumerate}
}

\subsection{Proof of Theorem~\ref{thrm:complete-Benders-distance-measure}}\label{app-proof:thrm:complete-Benders-distance-measure}

\begin{repeattheorem}[Theorem~\ref{thrm:complete-Benders-distance-measure}.]
Let $d_h$ be a Benders normalized distance function with $h$ a convex piece-wise linear function. Then BD Algorithm~\ref{pseudo-code-bd-modified} converges to an optimal solution or asserts infeasibility of MP in a finite number of iterations.
\end{repeattheorem}
\prf{First, we show that the BD algorithm does not stagnate in a degenerate loop. Let $\hat\Pi_t$ be the set of dual solutions obtained before iteration $t$ of the BD algorithm. Let $\mbox{MP}^{\itt}$ be the current approximation of MP with $({\yy}^{\itt},\eta^{\itt})$ its optimal solution and let $(\bar{\BFpi},\bar\pi_0)$ be the dual solution obtained from BSP \eqref{eq:separation_general} for separating $({\yy}^{\itt},\eta^{\itt})$.
If $\bar{\BFpi}^{\top} ({\bb}-B{\yy}^{\itt})+\bar\pi_0({\ff}^{\top}{\yy}^{\itt}-\eta^{\itt})=0$, then $({\yy}^{\itt},\eta^{\itt})$ is optimal for MP since $\eta^{\itt}$ is a lower bound on the optimal value of MP. Hence, assume that $\bar{\BFpi}^{\top} ({\bb}-B{\yy}^{\itt})+\bar\pi_0({\ff}^{\top}{\yy}^{\itt}-\eta^{\itt})>0$. Since $({\yy}^{\itt},\eta^{\itt})$ is feasible for $\mbox{MP}^{\itt}$, it follows that $\hat{\BFpi}^{\top} ({\bb}-B{\yy}^{\itt})+\hat{\pi}_0({\ff}^{\top}{\yy}^{\itt}-\eta^{\itt})\le 0$ for each $(\hat{\BFpi},\hat{\pi}_0)\in\hat\Pi_t$; hence, $(\bar{\BFpi},\bar\pi_0)$ cannot be a conical (i.e., scaling or a convex) combination of the solutions contained in $\hat\Pi_t$, meaning that, at each iteration, \eqref{eq:separation_general} will produce a cut that is not implied by the cuts hitherto obtained.

Finally, since $h$ is positive homogeneous and $\Pi$ is a cone, by Proposition~\ref{prop:subproblem-transformation} we can restate the separation subproblem \eqref{eq:separation_general} as
\begin{align}
    \max_{({\BFpi},\pi_0)\in\Pi_h} \;&{\BFpi}^{\top} ({\bb}-B{\yy}^{\itt})+\pi_0({\ff}^{\top}{\yy}^{\itt}-\eta^{\itt}),\label{generic-complete-concave-transformation}
\end{align}
where $\Pi_h=\{({\BFpi},\pi_0)\in\Pi:h({\BFpi},\pi_0)\le 1\}$.  Since $h$ is a convex piece-wise linear function, $\Pi_h$ is a polyhedron. Let $\Pi^v_h$ and $\Pi^r_h$ be the set of extreme points and rays of $\Pi_h$, respectively. Note that $\Pi^v_h\subset \Pi$ and $\Pi^r_h\subset \Pi$, and that they do not depend on $({\yy}^{\itt},\eta^{\itt})$. If \eqref{generic-complete-concave-transformation} is bounded, then its optimal solution is attained at one of the points in $\Pi^v_h$, otherwise an extreme ray belonging to $\Pi^r_h$ causes unboundedness. Either way, the produced extreme point/ray of $\Pi_h$ serves as the certificate. Therefore, the number of iterations is bounded by $|\Pi^v_h|+|\Pi^r_h|$.}

\section{Reformulation of $\ell_p$ Separation Problems as Linear/Quadratic Programs}\label{app:reformulations-lp}
Here we show how SSP \eqref{separation-subproblem-Lp} can be cast as an LP/QP using standard reformulation techniques. By homogeneity of $\ell_p$-norms, using Proposition~\ref{prop:subproblem-transformation} we can rewrite
SSP \eqref{separation-subproblem-Lp} as
\begin{equation}
    \begin{split}\label{separation-problem-lp}
        \max_{({\BFpi},\pi_0)\in\Pi}\quad&{\BFpi}^{\top} ({\bb}-B\hat{\yy})+\pi_0({\ff}^{\top}\hat{\yy}-\hat{\eta})\\
        \text{s.t.}\quad&\|({\BFpi}^{\top} B-\pi_0{\ff}^{\top},\pi_0)\|_p\le 1.
\end{split}
\end{equation}
We may express the constraint $\|({\BFpi}^{\top} B-\pi_0{\ff}^{\top},\pi_0)\|_p\le 1$ as a set of linear/quadratic constraints depending on the choice of $p$ as follows.
\paragraph{For $p=\infty$:} We have $\|({\BFpi}^{\top} B - \pi_0{\ff}^{\top},\pi_0)\|_{\infty}=\max\left\{\pi_0, \max\limits_{j=1,\dots,n}\{|{\BFpi}^{\top} B_{.j}-\pi_0f_j|\}\right\}$, where $B_{.j}$ is the $j$'th column of matrix $B$. Therefore, $\|({\BFpi}^{\top} B - \pi_0{\ff}^{\top},\pi_0)\|_{\infty}\le 1$ can be represented by the $2n$ linear constraints $-1\le {\BFpi}^{\top} B_{.j}-\pi_0f_j\le 1$ for each $j$, and a bound constraint $\pi_0\le 1$.
    
\paragraph{For $p=1$:} We may rewrite $\|({\BFpi}^{\top} B - \pi_0{\ff}^{\top},\pi_0)\|_1=\pi_0+\sum_{j=1}^{n}|{\BFpi}^{\top} B_{.j}-\pi_0f_j|\le 1$ as $\pi_0+\sum_{j=1}^n \tau_j\le 1$ by introducing $n$ new variables $\BFtau\in \mathbb{R}^n_+$ and $2n$ constraints $-\BFtau\le {\BFpi}^{\top} B-\pi_0 \ff\le \BFtau$.
    
\paragraph{For $p=2$:} One only needs to rewrite $\|({\BFpi}^{\top} B - \pi_0{\ff}^{\top},\pi_0)\|_2\le 1$ as ${\pi}^2_0+\sum_{j=1}^{n}({\BFpi}^{\top} B_{.j}-\pi_0f_j)^2\le 1$ to cast \eqref{separation-problem-lp} as a convex quadratically constrained linear program.
    
\paragraph{For $p>2$ and integer:} Note that $\|({\BFpi}^{\top} B - \pi_0{\ff}^{\top},\pi_0)\|_p\le 1$ is equivalent to ${\pi}^p_0+\sum_{j=1}^{n}\tau_j^p\le 1$, where $-\tau_j\le{\BFpi}^{\top} B_{.j}-\pi_0f_j\le \tau_j$. The constraint ${\pi}^p_0+\sum_{j=1}^{n}\tau_j^p\le 1$ can be expressed as quadratic constraints using a series of transformations.  
For instance, with $p=4$, it is not difficult to see that $\pi_0^4+\sum_{j=1}^{n}\tau_j^4\le 1$ may be expressed using auxiliary variables $\{\beta_j\}_{j=0}^n$ as the following set of second-order constraints \eqref{eq:transform-second-order}; similar transformations may be used for other values of $p$.
\begin{align}
    \beta_0^2+\sum_{j=1}^{n}\beta_j^2\le 1,\quad \pi_0^2\le \beta_0,\quad \tau_j^2\le \beta_j\quad \forall j. \label{eq:transform-second-order}
\end{align}

\section{A Monotonicity Property of Distance Functions}\label{app:monotonicity}
In this appendix, we formally define Benders distance functions and present a notion of mononicity for them.
Formally, we define a Benders distance function, which is a generalization of the geometric distance functions induced by $\ell_p$-norms presented earlier, as follows.

\begin{definition}[Benders distance function]\label{def:Benders-distance}
Function $d(\hat{\yy},\hat{\eta}|{\BFpi},\pi_0): \mathbb{R}^{n+1}\times\Pi\to \mathbb{R}$ is a Benders distance function if (i) it certifies $d(\hat{\yy},\hat{\eta}|{\BFpi},\pi_0)> 0$ iff $(\hat{\yy},\hat{\eta})$ is exterior to $\cut(\BFpi,\pi_0)$, $d(\hat{\yy},\hat{\eta}|{\BFpi},\pi_0)= 0$ iff $(\hat{\yy},\hat{\eta})$ is on the boundary of $\cut(\BFpi,\pi_0)$, and $d(\hat{\yy},\hat{\eta}|{\BFpi},\pi_0)< 0$ iff $(\hat{\yy},\hat{\eta})$ is in the interior of $\cut(\BFpi,\pi_0)$, and (ii) $d^*(\hat{\yy},\hat{\eta})$ defined as the objective value of the Benders separation problem (BSP) \eqref{eq:separation_general} is convex.
\end{definition}

\begin{definition}[Epigraph distance function]\label{def:Benders-distance3}
For a given Benders distance function $d$, we call $d^*$ as defined in \eqref{eq:separation_general} the epigraph distance function induced by $d$.
\end{definition}

The sign of $d^*(\hat{\yy},\hat{\eta})$ determines if $(\hat{\yy},\hat{\eta})$ is in the exterior ($\text{sign}=+1$), interior ($\text{sign}=-1$) or on the boundary of $\epi$ ($\text{sign}=0$). 
Moreover, the following definition  provides a weak characterization of how we expect distance functions to behave. That is, given any line segment between a point on the epigraph and a point outside the epigraph, $d^*$ monotonically increases as we move farther away from the epigraph (see Figure~\ref{fig:monotonic_distances}).

\begin{figure*}[h]
    \centering
    \includegraphics[clip,width=0.45\textwidth]{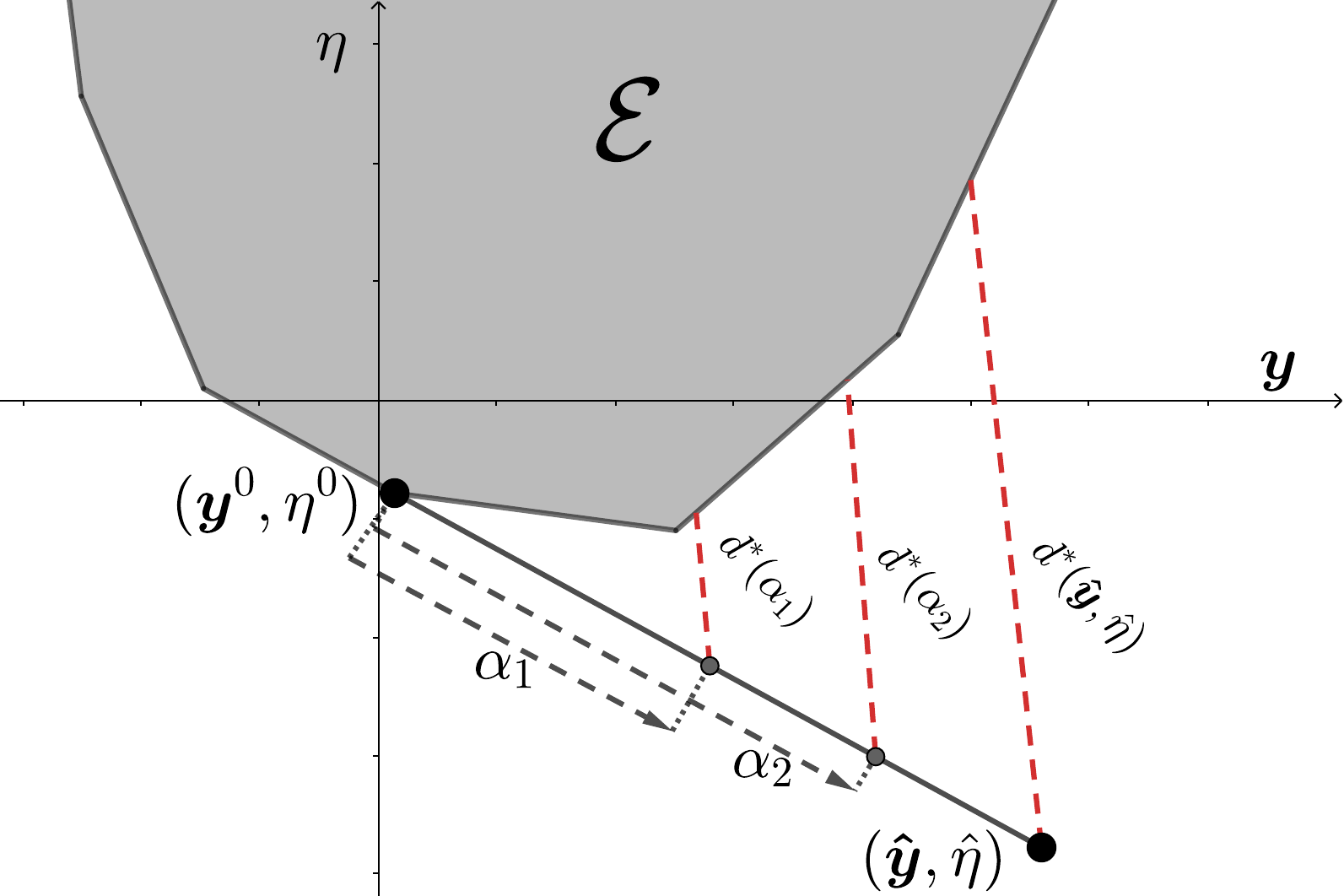}
    \caption{Epigraph distance functions are monotonic. As we move away from the boundary of $\epi$, $d^*$ gets larger.}
   \label{fig:monotonic_distances}
\end{figure*}

\begin{definition}[Monotonicity of epigraph distance function] For arbitrary $(\hat{\yy},\hat{\eta})\notin \epi$ and $({\yy}^0,{\eta}^0)\in \partial \epi$ (boundary of $\epi$) such that the open line segment between $(\hat{\yy},\hat{\eta})$ and $({\yy}^0,{\eta}^0)$ lies in the exterior of $\epi$, define $d^*(\alpha)=d^*((1-\alpha)({\yy}^0,{\eta}^0)+\alpha(\hat{\yy},\hat{\eta}))$.
We say the epigraph distance function $d^*$ is \textbf{monotonic} if $d^*(\alpha_1)\le d^*(\alpha_2)$ for any $0\le\alpha_1<\alpha_2\le 1$. We say $d^*$ is \textbf{strictly monotonic} if $d^*(\alpha_1)<d^*(\alpha_2)$ for any $0\le\alpha_1 < \alpha_2\le 1$.
\end{definition}

Due to the convexity of $d^*$, we show in Theorem~\ref{thrm:monotonic_convex} that any Benders distance function induces an epigraph distance function that is a measure of how far $(\hat{\yy},\hat{\eta})$ is from the boundary of $\epi$.

\begin{theorem}\label{thrm:monotonic_convex}
    Epigraph distance functions are monotonic.
\end{theorem}
\prf{Let $\alpha_i\in[0,1]$ for $i=1,2$ and assume that $\alpha_2>\alpha_1$. Define $(\bar{\yy}^{(i)},\bar{\eta}^{(i)})=(1-\alpha_i)({\yy}^0,{\eta}^0)+\alpha_i(\hat{\yy},\hat{\eta})$ for $i=1,2$. 
Since $0\le\alpha_1<\alpha_2\le 1$, we may state $(\bar{\yy}^{(1)},\bar{\eta}^{(1)})$ as a convex combination of $(\bar{\yy}^{(2)},\bar{\eta}^{(2)})$ and $({\yy}^0,{\eta}^0)$ of the following form
$$(\bar{\yy}^{(1)},\bar{\eta}^{(1)})=(1-\frac{\alpha_1}{\alpha_2})({\yy}^0,{\eta}^0)+\frac{\alpha_1}{\alpha_2}(\bar{\yy}^{(2)},\bar{\eta}^{(2)}).$$
Convexity of $d^*$ implies that
\begin{align*}
    d^*(\alpha_1)=d^*(\bar{\yy}^{(1)},\bar{\eta}^{(1)})&\le (1-\frac{\alpha_1}{\alpha_2})d^*({\yy}^0,{\eta}^0)+\frac{\alpha_1}{\alpha_2}d^*(\bar{\yy}^{(2)},\bar{\eta}^{(2)})=\frac{\alpha_1}{\alpha_2}d^*(\bar{\yy}^{(2)},\bar{\eta}^{(2)})\\ &\le d^*(\bar{\yy}^{(2)},\bar{\eta}^{(2)})=d^*(\alpha_2),
\end{align*}
where the we have used $d^*({\yy}^0,{\eta}^0)=0$ because $({\yy}^0,{\eta}^0)\in\partial \epi$. Hence $d$ is monotonic.
}

Finally, we note that normalized distance functions introduced in Section \ref{sec:normalized-distance-functions} are well-defined Benders distance functions.
\begin{proposition}\label{prop:normalized-monotonic}
    Normalized distance function $d_h$ induces a monotonic epigraph distance function $d^*_h$ for any positive homogeneous function $h$.
\end{proposition}
\prf{Since $h$ is positive homogeneous, using Proposition \ref{prop:subproblem-transformation}, for any $(\bar{\yy},\bar{\eta})$ we may state $d^*(\bar{\yy},\bar{\eta})$ as 
$$d^*_h(\bar{\yy},\bar{\eta})=\max_{({\BFpi},\pi_0)\in\Pi_h} \;{\BFpi}^{\top} ({\bb}-B\bar{\yy})+\pi_0({\ff}^{\top}\bar{\yy}-\bar{\eta}),$$
where $\Pi_h=\{({\BFpi},\pi_0)\in\Pi:h({\BFpi},\pi_0)\le 1\}$. This implies that $d^*_h$ is convex, since it is the maximum of a number of linear functions. Therefore, by Theorem \ref{thrm:monotonic_convex}, $d_h$ is monotonic.}

\section{Separable Subproblems}\label{app:separable_subproblem}
A common application of BD arises when the constraint matrix is block diagonal, which means we can partition the continuous variables and the constraint matrices into $|K|$ independent and mutually exclusive subgroups and rewrite OP \eqref{generic-milp} as
\begin{equation}
    \begin{split}\label{generic-milp-separable}
        \min\quad & \sum_{k\in K}p_k{\cc}_{k}^{\top}{\xx}_{k}+{\ff}^{\top}{\yy}\\
    \text{s.t.} \quad & A_{k}{\xx}_{k}+B_{k}{\yy}\ge \bb_{k} \qquad \forall k\in K\\
    & {\xx}_{k} \ge {\vzero}, {\yy}\in Y \qquad \forall k\in K
    \end{split}
\end{equation}
Examples of this are two-stage stochastic optimization programs, where $\yy$ is the first stage decision variable and $\xx_k$ is the second stage decision variable under realization $k\in K$ which occurs with probability $p_k$ (i.e., $\sum_{k\in K}p_k=1$). Note that this formulation is sufficiently general to encompass general MILPs with a block diagonal structure, such as uncapacitated facility/network-type problems, where $k$ corresponds to commodity (or customer) $k$, $p_k=1/|K|$ for all $k$, and $p_k{\cc}_{k}$ is the cost vector of commodity $k$. 

Given that $\sum_{k\in K}p_k=1$, we may reformulate \eqref{generic-milp-separable} by introducing $|K|$ auxiliary variables $\eta_k$ as
\begin{align}
        \min\quad & \sum_{k\in K}p_k\eta_k\label{generic-milp-separable-epi-obj}\\
    \text{s.t.} \quad & \eta_k\ge {\cc}_{k}^{\top}{\xx}_{k}+{\ff}^{\top}{\yy} & \forall k\in K\,\label{generic-milp-separable-epi-ub}\\
    & A_{k}{\xx}_{k}+B_{k}{\yy}\ge \bb_{k} & \forall k\in K\,\label{generic-milp-separable-epi-c}\\
    & {\xx}_{k} \ge {\vzero}, {\yy}\in Y & \forall k\in K.\label{generic-milp-separable-epi-domain}
\end{align}
Hence, for given $\hat{\yy}$ and $\hat{\BFeta}=(\hat{\eta}_k)_{k\in K}$, CGSP \eqref{generic-milp-cgsp} becomes
\begin{align}
    \max_{\{(\BFpi_{k},\pi_{k0})\in \Pi_{k}\}_{k\in K}}\quad & \sum_{k\in K}{\BFpi}_{k}^{\top}({\bb_{k}}-B_{k}\hat {\yy})+\pi_{k0}({\ff}^{\top}\hat{\yy}-\hat{\eta}_k),\label{generic-milp-cgsp-separable}
\end{align}
where $\pi_{0k}$ and $\BFpi_{k}$ are the dual variables associated with \eqref{generic-milp-separable-epi-ub} and \eqref{generic-milp-separable-epi-c}, respectively, and
$$\Pi_{k}=\{({\BFpi}_{k},\pi_{k0}): {\BFpi}_{k}^{\top} A_{k}\le \pi_{0k} {\cc}_{k}^{\top},\; {\BFpi}_{k} \ge {\vzero}, \pi_{0k} \ge 0\}.$$ 
Note that, while \eqref{generic-milp-cgsp-separable} decomposes into $|K|$ subproblems, given that the cuts live in the space of $(\yy,\BFeta)$, the normalized separation problem \eqref{normalized-separation-problem} need not be separable for a general normalization function $h$.
Regardless, upon deriving the certificate $(\hat{\BFpi}_{k},\hat{\pi}_{k0})_{k\in K}$, one may produce a single cut
\begin{align}
    \sum_{k\in K}\hat{\BFpi}_{k}^{\top}({\bb_{k}}-B_{k} {\yy})+\hat{\pi}_{k0}({\ff}^{\top}{\yy}-{\eta}_k)\le 0,\label{separable-subproblem-single-cut}
\end{align}
or multiple cuts, one for each $k$, which imply \eqref{separable-subproblem-single-cut}:
\begin{align}
    \hat{\BFpi}_{k}^{\top}({\bb_{k}}-B_{k} {\yy})+\hat{\pi}_{k0}({\ff}^{\top}{\yy}-{\eta}_k)\le 0 \qquad & \forall k\in K.\label{separable-subproblem-multi-cut}
\end{align}

\subsection{Separable Normalization Functions}
To recover separability, it suffices to use a normalization function of the form
\begin{align}
    h((\BFpi_{k},\pi_{k0})_{k\in K})=\max_{k\in K}\left\{h_k(\BFpi_{k},\pi_{k0})\right\},\label{normalization-function-separable}
\end{align}
where $h_k(\BFpi_{k},\pi_{k0})$ is the block-level normalization function of block $k\in K$. Note that $h$ inherits positive homogeneity from the block-level normalization functions, hence is a valid normalization function. Moreover, $h((\BFpi_{k},\pi_{k0})_{k\in K})\le 1$ is equivalent to $h_k(\BFpi_{k},\pi_{k0})\le 1$ for each $k\in K$, thus \eqref{normalized-separation-problem}  decomposes into $|K|$ independent normalized separation problems.
As an example, setting $h_k(\BFpi_{k},\pi_{k0})=\|{\BFpi}_{k}^{\top}B_{k}-\pi_{0k}{\ff}^{\top}, \pi_{0k}\|_p$ implies that $h((\BFpi_{k},\pi_{k0})_{k\in K})$ is a composite of $\ell_{\infty}$ and $\ell_p$ norms.

\subsection{Stochastic Projections for Projective Normalization Functions}
Recall that a projective normalization function $g$ amounts to finding a point $(\tilde{\yy}, \tilde{\eta})\in \epi$ at which the cut(s) produced according to $d_g$ support $\epi$. Conversely, had we known $\tilde{\yy}$, we could produce the desired cuts by evaluating $|K|$ independent classical Benders subproblems at $\tilde{\yy}$. Given the stochastic nature of two-stage stochastic programs, it is reasonable to approximate $\tilde{\yy}$ using a sample $\tilde{K}\subset K$ instead of finding the true projection point $\tilde{\yy}$ using the full set $K$. Consequently, we can produce the cuts for the remaining blocks by simply solving classical Benders subproblems, for which we can exploit the efficiency of the oracle. Algorithm~\ref{alg:sampled-projections} provides an overview of this procedure. 

\begin{algorithm}[ht!]
	\caption{Stochastic Projections Algorithm}
	\label{alg:sampled-projections}
	\begin{algorithmic}[1]
        \State \textbf{STEP 0 (Sampling):} Take a sample $\tilde{K}$ of the blocks (scenarios). Let $\epi(\tilde{K})$ be the epigraph of OP when $K$ is restricted to $\tilde{K}$, and probabilities adjusted to $\tilde{p}_k = p_k/\sum_{l\in \tilde{K}}p_l$.
        \State  \textbf{STEP 1 (Projection):} Find projection of $(\hat{\yy},\hat{\BFeta})$ onto $\epi(\tilde{K})$. Let $(\bar{\yy},\bar{\BFeta})$ be this projection, and $(\bar{\BFpi}_k, \bar{\pi}_{0k})$ be the dual solution associated with block $k\in \tilde{K}$.
	    \For{$k \in K\setminus\tilde{K}$}
        	\State \textbf{STEP 2 (Cut Generation):} Attempt solving the following classical PSP:
	        \begin{align*}
	  [\text{PSP}]\quad\max\left\{ \cc^{\top}_k\xx_k: A_k\xx_k\ge  {\bb_k}-B_k\bar{\yy}, \xx_k \ge {\vzero}\right\}.
	        \end{align*}
  \State \textbf{if} PSP is feasible \textbf{then} set $(\bar{\BFpi}_{k},\bar{\pi}_{0k})=(\bar{\uu}_k,1)$ where $\bar{\uu}_k$ is an optimal dual solution. 
         \State \textbf{Else} set $(\bar{\BFpi}_{k},\bar{\pi}_{0k})=(\bar{\vv}_k,0)$, where $\bar{\vv}_k$ is a Farkas certificate.  
  \State \textbf{STEP 3 (Validation):} \textbf{if} $(\bar{\BFpi}_{k},\bar{\pi}_{0k})$ does not cut off $(\hat{\yy}, \hat{\eta}_k)$ \textbf{then} repeat Step 2 with $\hat{\yy}$ instead of $\bar{\yy}$.
	    \EndFor
	\end{algorithmic}
\end{algorithm}

A few remarks are in order. First, note that the projection point in Step 1 comes as a byproduct of solving the separation problem according to $g$. When using a blackbox solver, we can extract $\tilde{\yy}$ as the shadow price of $\BFtau = \sum_{k\in \tilde{K}}{\BFpi}_{k}^{\top}B_{k}-\pi_{0k}{\ff}^{\top}$. For the separable $d_g$, we can approximate $\tilde{\yy}=\sum_{k\in \tilde{K}}\tilde{p}_k\tilde{\yy}_k$, where $\tilde{p}_k$ is the adjusted probability of block $k\in \tilde{K}$, and $\tilde{\yy}_k$ is the projection point according to the $k^{\text{th}}$ block (i.e., shadow price of $\BFtau_k = {\BFpi}_{k}^{\top}B_{k}-\pi_{0k}{\ff}^{\top}$.) As detailed in Sections \ref{sec:directed-depth-maximizing} and \ref{sec:guided_projections_algorithm}, we can also extract the projection point in a combinatorial fashion using GPA and DDMA.

Second, given that the approximate projection point is not necessarily the true projection point for $k\notin \tilde{K}$, the cut produced in Step 2 need not cut off $(\hat{\yy}, \hat{\eta}_k)$. Therefore, to ensure convergence of the BD algorithm, we produce a classical Benders cut whenever necessary (Step 3).

Finally, while a larger sample size generally results in a better approximation of the projection point, a smaller sample size may result in an overall shorted computation time. In our experiments, we found that a sample size of as small as one randomly selected scenario may balance the overall computing time and quality of produced cuts.

\section{Description of Benchmark Instances}\label{app:instances}
\subsection{Capacitated Facility Location Problem}\label{app:cflp}
Facility location problems lie at the heart of network design and planning, and arise naturally in a wide range of applications such as supply chain management, telecommunications systems, urban transportation planning, health care systems and humanitarian logistics to count a few \citep[see e.g., ][]{drezner2001facility}.
\subsubsection{Formulation.}
Given a set of $m$ customers and $n$ potential locations for the facilities, CFLP in its simplest form as formulated below, consists of determining which facilities to open and how to assign customers to opened facilities to minimize cost, i.e.,
\begin{align}
    \text{[CFLP]}\quad \min\quad&\sum_{i=1}^{m} \sum_{j=1}^{n} x_{ij} d_{i}c_{ij}+\sum_{j=1}^n f_j y_j\label{cflp-obj}\\
    \text{s.t.}\quad & \sum_{j=1}^{n} x_{ij}\ge 1  &&\forall i\label{cflp-demand}\\
    & \sum_{i=1}^{m} x_{ij} d_{i}\le s_{j} y_{j} &&\forall j\label{cflp-cap}\\
    & x_{ij}\le y_{j}&&\forall i,j\label{cflp-valid}\\
    & \xx\ge \vzero, \yy\in Y,\label{cflp-domain}
\end{align}
where $f_j$ and $s_j$ are respectively the installation cost and capacity of facility $j=1,\dots,n$; $d_i$ is the demand of customer $i=1,\dots,m$; $c_{ij}$ is the cost of serving one unit of demand from customer $i$ using facility $j$; and $Y=\{\yy\in\{0,1\}^n:\sum_{j=1}^{n} s_j y_{j}\ge \sum_{i=1}^{m}d_{i}\}$ is the domain of the $\yy$ variables. 

The two-stage stochastic program with $|K|$ demand scenarios is formulated as follows
\begin{align}
    \text{[CFLP-S]}\quad \min\quad&\frac{1}{|K|}\sum_{k\in K}\sum_{i=1}^{m} \sum_{j=1}^{n} x_{ijk} d_{ik}c_{ij}+\sum_{j=1}^n f_j y_j\label{scflp-obj}\\
    \text{s.t.}\quad & \sum_{j=1}^{n} x_{ijk}\ge 1  &&\forall i, k\label{scflp-demand}\\
    & \sum_{i=1}^{m} x_{ijk} d_{ik}\le s_{j} y_{j} &&\forall j, k\label{scflp-cap}\\
    & x_{ijk}\le y_{j}&&\forall i,j,k\label{scflp-valid}\\
    & \xx\ge \vzero, \yy\in Y,\label{scflp-domain}
\end{align}
in which $Y$ is now defined as $Y=\{\yy\in\{0,1\}^n:\sum_{j=1}^{n} s_j y_{j}\ge \max_{k\in K}\sum_{i=1}^{m}d_{ik}\}$.

A core point $\bar{\yy}$ can be obtained by setting $\bar{y}_j=\frac{1}{r}+\epsilon$ for each $j$, where $\epsilon=10^{-3}$ and $r=\frac{\sum_{j=1}^{n}s_j}{\sum_{i=1}^{m}d_i}$ in the deterministic case and $r=\frac{\sum_{j=1}^{n}s_j}{\max_{k\in K}\sum_{i=1}^{m}d_{ik}}$ in the stochastic case.

\subsubsection{Instances.}
We used two sets of benchmark instances from the literature:

\noindent\texttt{\textbf{CAP}}: The famous \texttt{CAP} data set from the OR-Library \citep{orlib} consists of 24 small instances with $m=50$ customers and $n\in\{16,25,50\}$ facilities, and 12 large instances with $n=100$ facilities and $m=1000$ customers. The instances are denoted \texttt{CAPx1}--\texttt{CAPx4}, where $x\in\{6,7,9,10,12,13\}$ for the small instances and $x\in\{a,b,c\}$ for the large instances.

\noindent \texttt{\textbf{CST}}: This set, also known as the \texttt{GK} dataset \citep{gortz2012simple}, contains a set of randomly generated instances following the procedure proposed by \cite{cornuejols1991comparison}. We denote each instance by tuple $(n,m,r)$, where $(n,m)$ pairs were selected from $\{$(50,50), (50,100), (100,100), (100,200), (100,500), (500,500), (100,1000), (200,1000), (500,1000), (1000,1000)$\}$ and the scaling factor $r$ was selected from $\{5,10,15,20\}$. For each choice of $(n,m,r)$ we randomly generated 4 instances as follows. For each facility $j\in\{1,\dots,n\}$, we randomly drew $s_j$ and $f_j$ from $U[10,160]$ and $U[0,90]+U[100,110]\sqrt{s_j}$, respectively, where $U[a,b]$ represents the uniform distribution on $[a,b]$. For each customer $i\in\{1,\dots,m\}$, we randomly drew $d_{i}$ from $U[5,35]$. Finally, we scaled the facility capacities using parameter $r$ such that 
$r=\frac{\sum_{j=1}^{n}s_j}{\sum_{i=1}^{m}d_i}.$
To compute the allocation costs, we placed the customers and facilities in a unit square uniformly at random, and set $c_{ij}$ to 10 times the Euclidean distance of facility $j$ from customer $i$.

\subsubsection{Deriving a classical Benders cut.}\label{app:cflp-cut}
A classical Benders cut for separating a given (fractional) solution $\hat{\yy}$ can be derived efficiently as follows.
First, note that constraints \eqref{cflp-cap} can be treated as bounds on the primal variables $\xx$. Consequently, it suffices to update these bounds based on values of $\hat{\yy}$ and reduce the number of constraints in the primal subproblem from $n+m+nm$ to $n+m$ constraints. The resulting problem is a transportation problem, which can be solved efficiently for large instances using specialized algorithms. In our implementation, however, we have used \texttt{Cplex} for solving the transportation problems since it benefits from better warm starting.
Let $\hat{\uu}^{\text{D}}$ be the optimal dual solution associated with the demand constraints \eqref{cflp-demand} obtained by the solver. Given $\hat{\uu}^{\text{D}}$, as noted by several authors \citep[see e.g.,][]{cornuejols1991comparison, fischetti2016benders}, the Benders cut takes the form
\begin{align}
    \eta\ge \sum_{i=1}^{m}\hat{u}^{\text{D}}_i+\sum_{j=1}^{n}\left(f_j-\kappa_j(\hat{\uu}^{\text{D}})\right)y_j,\label{eq:cflp-cut}
\end{align}
where $\kappa_j(\hat{\uu}^{\text{D}})$ is the optimal value of the continuous knapsack problem
\begin{align*}
    \kappa_j(\hat{\uu}^{\text{D}})=\max\left\{\sum_{i=1}^{m}\alpha_i(\hat{u}^{\text{D}}_i-d_i c_{ij}): \sum_{i=1}^{m} \alpha_i d_{i}\le s_{j}, \BFalpha\in[0,1]^m\right\},
\end{align*}
which can be solved efficiently in $\calO(m)$ time by finding the weighted median of the ratios $\{\frac{\hat{u}^{\text{D}}_i}{d_i}-c_{ij}\}$ using the procedure described by  \cite{balas1980algorithm}.

\subsection{Uncapacitated Facility Location Problem}\label{app:uflp}
Uncapacitated facility location problem (UFLP), like its capacitated counterpart, is another important problem in operations research. UFLP is formulated as CFLP but without the capacity constraints \eqref{cflp-cap}.

\subsubsection{Instances.}
We used two sets of benchmark instances from the literature:

\noindent \texttt{\textbf{M*}}: This set from \cite{kratica2001solving} consists of 21 instances with $n=m\in \{100, 200, 300, 500, 1000\}$ facilities/customers.

\noindent \texttt{\textbf{KG}}: This set from \cite{ghosh2003neighborhood} consists of several instances of different parameter settings (symmetric/asymmetric and 3 subclasses representing different cost settings) with a total of 20 instances for each $n=m\in \{250, 500\}$. 

\subsubsection{Deriving a classical Benders cut.}\label{app:uflp-cut}
While UFLP shares a similar formulation with CFLP, deriving a classical Benders cut is fundamentally different in UFLP.  
A classical Benders cut for separating a given (fractional) solution $\hat{\yy}$ can be derived efficiently by solving the following continuous knapsack problem for each customer $i=1, \dots,m$:
\begin{align}
    \min\quad&\sum_{j=1}^{n} x_{ij} c_{ij}\label{uflp-sp-obj}\\
    \text{s.t.}\quad & \sum_{j=1}^{n} x_{ij}\ge 1  \label{uflp-sp-demand}\\
    & 0\le x_{ij}\le \hat{y}_{j}&& \forall j,\label{uflp-sp-valid}
\end{align}
which can be solved by simply sorting the facilities for each customer $i$ in the non-increasing order of $c_{ij}$. Let $k$ denote the index of the critical item (facility). Then, the shadow price of \eqref{uflp-sp-demand} is $c_{ik}$, while the shadow price of \eqref{uflp-sp-valid} is $\max\{0, c_{ik}-c_{ij}\}$ for each $j$.

\subsection{Multicommodity Capacitated Network Design Problem}\label{app:mcndp}
Network design problems concern selecting a set of arcs from a set of candidate arcs in a network so that the demand of a set of origin-destination pairs (commodities) is routed with minimum cost in the resulting network. Let $A$ and $V$ denote the set of arcs and vertices of the underlying network. Each commodity $i\in I$ has $d_i$ units of demand to be routed from node $O(i)\in V$ to $D(i)\in V$, and a unit cost $c_{ij}$ over arc $j\in A$. Each arc $j\in A$ has a fixed installation cost $f_j$ and capacity $s_j$. 

\subsubsection{Formulation.}
Let $y_j\in\{0,1\}$ be a binary variable determining if arc $j\in A$ is chosen, and let variable $x_{ij}$ measure the fraction of demand of commodity $i\in I$ routed on arc $j\in A$. We formulate MCNDP as
\begin{align}
    \text{[MCNDP]}\quad \min\quad&\sum_{i\in I} \sum_{j\in A} x_{ij} d_{i}c_{ij}+\sum_{j\in A} f_j y_j\label{mcndp-obj}\\
    \text{s.t.}\quad & \sum_{j\in A^+_v} x_{ij}-\sum_{j\in A^-_v} x_{ij}\ge b_{iv}  &&\forall i\in I, v\in V\label{mcndp-demand}\\
    & \sum_{i\in I} x_{ij} d_{i}\le s_{j} y_{j} &&\forall j\in A\label{mcndp-cap}\\
    & x_{ij}\le y_{j}&&\forall i\in I,j\in A\label{mcndp-valid}\\
    & \xx\ge \vzero, \yy\in Y,\label{mcndp-domain}
\end{align}
where $A^+_v$ and $A^-_v$ denote the set of outward and inward arcs to node $v\in V$, $b_{iv}=1$ for $v=O(i)$, $b_{iv}=-1$ for $v=D(i)$, $b_{iv}=0$ for $v\in V\setminus \{O(i), D(i)\}$, and $Y=\{0,1\}^{|A|}$ is the domain of $\yy$.

The stochastic counterpart can be formulated as follows:
\begin{align}
    \text{[MCNDP-S]}\quad \min\quad&\frac{1}{|K|}\sum_{k\in K}\sum_{i\in I} \sum_{j\in A} x_{ijk} d_{ik}c_{ij}+\sum_{j\in A} f_j y_j\label{smcndp-obj}\\
    \text{s.t.}\quad & \sum_{j\in A^+_v} x_{ijk}-\sum_{j\in A^-_v} x_{ijk}\ge b_{iv}  &&\forall i\in I, v\in V, k\in K\label{smcndp-demand}\\
    & \sum_{i\in I} x_{ijk} d_{ik}\le s_{j} y_{j} &&\forall j\in A, k\in K\label{smcndp-cap}\\
    & x_{ijk}\le y_{j}&&\forall i\in I,j\in A, k\in K\label{smcndp-valid}\\
    & \xx\ge \vzero, \yy\in Y.\label{smcndp-domain}
\end{align}

Producing a core point in MCNDP is not as straightforward as in CFLP. To produce a feasible core point, we first solve \eqref{mcndp-obj}--\eqref{mcndp-domain} by setting $y_j=1$ for all $j\in A$ to obtain the optimal solution $\xx^*$. We then set $\bar{y}_j=\max\left\{\epsilon,\frac{1}{s_j}\sum_{i\in I}x^*_{ij}d_i, \max_{i\in I}x^*_{ij}\right\}$ where $\epsilon=10^{-3}$ ensures that $\bar{y}_j>0$. The stochastic version requires repeating this procedure for each $k\in K$, and taking $\bar{y}_{j}=\max_{k\in K}\bar{y}_{jk}$, where $\bar{y}_{jk}=\max\left\{\epsilon,\frac{1}{s_j}\sum_{i\in I}x^*_{ijk}d_{ik}, \max_{i\in I}x^*_{ijk}\right\}$.

\subsubsection{Instances.}
We considered 10 classes of instances from the \texttt{R} set (\texttt{R01}--\texttt{R10}), each consisting of 9 instances with different capacities and fixed costs \citep{crainic2001bundle}. For the deterministic case, we used all 9 instances in each class. 
The number of arcs in these class range from 35 to 120, while the number of commodities range from 10 to 50. Nine out of 90 instances turned out to be infeasible when producing a core point as described above.
For the stochastic case, we selected instances $\{1,3,5,7,9\}$ from each class to capture different combinations of parameters. For each commodity $i$ with nominal demand $d_i$, we generated $|K|\in\{16, 32, 64\}$ demand scenarios from a triangular distribution with lower limit 0, upper limit $1.35d_i$ and mode $d_i$ as suggested in \cite{crainic2011progressive}.

\subsubsection{Deriving a classical Benders cut.}\label{app:mcndp-cut}
We may derive a classical Benders cut for separating a given (fractional) solution $\hat{\yy}$ in a two-step procedure similar to CFLP. As in CFLP, we note that constraints \eqref{mcndp-cap} can be treated as bounds on the primal variables $\xx$. Consequently, it suffices to update these bounds based on values of $\hat{\yy}$.
First, assume that the the resulting linear program is feasible, and $\hat{u}^{\text{D}}_{iv}$ be the optimal dual solution associated with the flow conservation constraints \eqref{mcndp-demand} obtained by the solver. It is not difficult to show that, given $\hat{\uu}^{\text{D}}$, the Benders optimality cut takes the form
\begin{align}
    \eta\ge \sum_{i\in I}\sum_{v\in V}\hat{u}^{\text{D}}_{iv}b_{iv}+\sum_{j\in A}\left(f_j-\kappa_j(\hat{\uu}^{\text{D}})\right)y_j,\label{eq:mcndp-cut-opt}
\end{align}
where $\kappa_j(\hat{\uu}^{\text{D}})$ is the optimal value of the continuous knapsack problem
\begin{align*}
    \kappa_j(\hat{\uu}^{\text{D}})=\min\left\{\sum_{i\in I}\alpha_i (\sum_{a\in A}(\hat{u}^{\text{D}}_{i, a^+}-\hat{u}^{\text{D}}_{i, a^-})-d_i c_{ij}): \sum_{i\in A} \alpha_i d_{i}\le s_{j}, \BFalpha\in[0,1]^{|I|}\right\},
\end{align*}
which can be solved efficiently in $\calO(|I|)$ time as before. Similarly, if the solver determines that subproblem is infeasible, with $\bar{u}^{\text{D}}_{iv}$ the components of the Farkas certificate associated with the flow conservation constraints \eqref{mcndp-demand}, we can express the Benders feasibility cut as 
\begin{align}
    0\ge \sum_{i\in I}\sum_{v\in V}\bar{u}^{\text{D}}_{iv}b_{iv}-\sum_{j\in A}\kappa_j(\bar{\uu}^{\text{D}})y_j,\label{eq:mcndp-cut-feas}
\end{align}
where $\kappa_j(\bar{\uu}^{\text{D}})$ is the optimal value of the continuous knapsack problem
\begin{align*}
    \kappa_j(\bar{\uu}^{\text{D}})=\min\left\{\sum_{i\in I}\alpha_i \sum_{a\in A}(\bar{u}^{\text{D}}_{i, a^+}-\bar{u}^{\text{D}}_{i, a^-}): \sum_{i\in A} \alpha_i d_{i}\le s_{j}, \BFalpha\in[0,1]^{|I|}\right\}.
\end{align*}

\subsection{Stochastic Network Interdiction Problem}\label{app:snip}
The stochastic network interdiction problem (SNIP), first introduced in \cite{pan2008minimizing}, concerns installing a predefined number of sensors on arcs of given network so as to minimize the expected probability that an intruder traverses the network undetected. As before, let $V$ and $A$ denote the set of vertices of the network. In the first stage, the interdictor is allowed to install at most $\sigma$ sensors on a subset $D\subseteq A$ of arcs. For each arc $a\in D$, the interdictor knows a priori the probabilities that the intruder avoids detection with or without a sensor installed on $a$, which are $r_a$ and $q_a$, respectively. What the interdictor does not know is the origin $s\in V$ and destination $t\in V$ of the intruder. In the second stage, under scenario $k$ which occurs with probability $p_k$, the intruder chooses a path from $s_k$ to $t_k$ that maximizes the probability of avoiding detection. 

\subsubsection{Formulation.}
Let $y_a\in \{0,1\}$ be a binary variable determining if a sensor is installed on arc $a\in D$, and let $x_{i,k}$ denote the probability that the intruder travels from $i\in V$ to $t_k$ undetected.
The formulation for SNIP is follows
\begin{align}
    \text{[SNIP]}\quad \min\quad&\sum_{k\in K}p_k x_{s_k, k}\label{snip-obj}\\
    \text{s.t.}\quad & x_{t_k,k} = 1 && \forall k\in K\\
    & x_{i,k} - q_a x_{j,k} \ge 0 && \forall a=(i,j)\in D, k\in K\\
    & x_{i,k} - r_a x_{j,k} \ge 0 && \forall a=(i,j)\in A\setminus D, k\in K\\
    & x_{i,k} - r_a x_{j,k} \ge -(r_a-q_q)\psi_{j,k}y_a && \forall a=(i,j)\in D, k\in K\\
    & \xx\ge \vzero, \yy\in Y,\label{snip-domain}
\end{align}
where $Y=\{\yy\in \{0,1\}^{|D|}: \sum_{a\in D} y_a\le \sigma\}$, $\sigma$ is the number of available sensors, and $\psi_{ik}$ is a parameter denoting the value of a maximum-reliability path from $i$ to $t_k$ when no sensors are placed, which can be obtained by solving a shortest path problem \citep{pan2008minimizing}.

A core point in SNIP can be produced by simply setting $\bar{y}_a=\frac{\sigma}{|D|}-\epsilon$ for some $\epsilon > 0$.

\subsubsection{Instances.}
The SNIP dataset, introduced by \cite{pan2008minimizing} and provided by \cite{bodur2016strengthened}, contains five network structures with 783 nodes and 2586 arcs. For each network, four parameter settings can be considered, out of which we considered the more challenging settings corresponding to \texttt{snipno3} and \texttt{snipno4}. We chose the budget parameter $\sigma\in\{30,40,50\}$, yielding a total of 30 instances, each with $|D|=320$ binary variables and $|K|=456$ scenarios.

\section{Implementation Details}\label{app:implementation_details}
We conducted our computational study on a Dell desktop equipped with Intel(R) Xeon(R) CPU E5-2680 v3 at 2.50GHz with 8 Cores and 32 GB of memory running a 64-bit Windows 10 operating system. 
We coded our algorithms in \texttt{C\#} and solved the linear/quadratic problems using the \texttt{ILOG Concert} library and \texttt{Cplex 12.10} solver. 
We implemented BD within B\&C in a modern fashion \citep{fortz2009improved,maher2021implementing} using the callback functionality of \texttt{Cplex}. Specifically, we treated Benders cuts for separating integer solutions as lazy constraints (invoked by \texttt{Cplex} using the \texttt{LazyConstraint} callback), and treated Benders cuts for separating fractional solutions as valid inequalities for the master problem (invoked by the \texttt{UserCut} callback).
In the following, we provide general implementation details for BD, which we believe are of technical value beyond the application of this paper.

\subsection{Coefficient Scaling}\label{app:coefficient-scaling}
A pitfall in implementing BD is that scales of master problem variables $\eta$ and $\yy$ are often unbalanced, meaning that the coefficient of $\eta$ in (optimality) cuts is often too small or too large compared to the coefficients of the $\yy$ variables. This imbalance poses numerical issues for the solver when handling the cuts. In addition, an imbalanced cut implies an imbalanced normalization function $h(\BFpi,\pi_0)=\|\BFpi^{\top}B-\pi_0\ff^{\top},\pi_0\|_p$, since the coefficients of $\eta$ and $\yy$ in the cut $\BFpi^{\top}\bb\le(\BFpi^{\top}B-\pi_0\ff^{\top})\yy+\pi_0\eta$ are $\pi_0$ and $\BFpi^{\top}B-\pi_0\ff^{\top}$, respectively.

Given that $\eta$ estimates $\cc^{\top}\xx + \ff^{\top}\yy$, scaling $\eta$ is equivalent to scaling the cost vectors $\cc$ and $\ff$. Therefore, to balance the scale of $\eta$ and other variables, we divide $(\cc,\ff)$ by $\beta>0$  in a pre-processing step. Note that this scaling does not affect the optimal solution $(\yy,\xx)$, but does affect $\eta$.
To choose a suitable value for $\beta$, we first solve DSP \eqref{sp-papadakos} using a core point $\bar{\yy}$ to obtain the dual solution $(\bar{\uu},1)$ and the optimality cut 
$\eta+(\bar{\uu}^{\top}B-\ff^{\top})\yy\ge \bar{\uu}^{\top}\bb$.
We then set $\beta$ as
$$\beta=\frac{1}{n}\|\bar{\uu}^{\top}B-\ff^{\top}\|_{1},$$
which is the average absolute coefficient value of the $\yy$ variables in the cut. 

\subsection{Reoptimizing the Separation Subproblems}\label{app:reoptimization}
Another important aspect in implementing the BD algorithm is being able to reoptimize the separation problems and retrieving the cuts quickly when a solver is used for solving the separation subproblems. Note that only the objective function in the separation problem \eqref{normalized-separation-problem} changes from one iteration of the BD Algorithm \ref{pseudo-code-bd-modified} to another. For linear separation subproblems, one can use the primal simplex algorithm by setting parameter \texttt{Cplex.Param.RootAlgorithm} to \texttt{Cplex.Algorithm.Primal} to leverage the reoptimization capabilities of this method.

Let us rearrange the objective function of the separation problem \eqref{normalized-separation-problem} at iteration $t$ of BD as
$${\BFpi}^{\top} {\bb}-\sum_{j=1}^{n}({\BFpi}^{\top}B_{.j}-\pi_0f_j){y}^{\itt}_j-\pi_0{\eta}^{\itt}.$$
Note that one needs to update the coefficient of ${\BFpi}^{\top}B_{.j}-\pi_0f_j$ (i.e., ${y}^{\itt}_j$) only 
when ${y}^{\itt}_j\ne{y}^{(t-1)}_j$. Therefore, we may additionally define $n$ auxiliary variables $\tau_{j}={\BFpi}^{\top}B_{.j}-\pi_0f_j$ to avoid changing the coefficients of all dual variables in the separation subproblems.
For instance, the separation problem \eqref{separation-problem-lp} for producing $\ell_p$-deepest cuts becomes
\begin{align*}
    \max \left\{{\BFpi}^{\top} {\bb}-\BFtau^{\top}\hat{\yy}-\pi_0\hat{\eta}: \|(\BFtau,\pi_0)\|_p\le 1,\; \BFtau=B^{\top}{\BFpi}-\pi_0{\ff},\; ({\BFpi},\pi_0)\in\Pi\right\}.
\end{align*}
Apart from updating the objective function coefficients in the separation problems, introducing the $\BFtau$ variables brings several advantages. First, it simplifies the expression for the projective normalization functions (e.g., in $\ell_p$-norm or in CW). Second, we can use these variables to retrieve the projection point for projective normalization function $g$. More precisely, adding the constraints $\BFtau=B^{\top}{\BFpi}-\pi_0{\ff}$, it is not difficult to verify that the $g$-projection of $\hat{\yy}$ is the vector of shadow prices of these constraints. Finally, whenever it is necessary to use the $\BFtau$ variables, after the subproblem is solved, one can save $\calO(mn)$ arithmetic operations in computing the cut coefficients by easily retrieving the value of the $\BFtau$ variables from the solver without having to recalculate the coefficients based on the $\BFpi$ variables.

\section{Supplementary Figures}\label{app:supplementary_figures}

Figure~\ref{fig:norms-truncate} illustrates the effect of $\ell_1$-, $\ell_2$- and $\ell_{\infty}$-norms on truncating the cone of dual solutions. For illustration, $\Pi$ is transformed from $\mathbb{R}^{m+1}$ to $\mathbb{R}^{n+1}$ as $\Gamma= \{(\BFtau,\pi_0)\in \mathbb{R}^{n+1}: \exists ({\BFpi},\pi_0)\in \Pi \text{ s.t. }\BFtau ={\BFpi}^{\top}B-\pi_0{\ff}^{\top}\}$. 

\begin{figure*}[h]
    \centering
    \includegraphics[clip,width=0.35\textwidth]{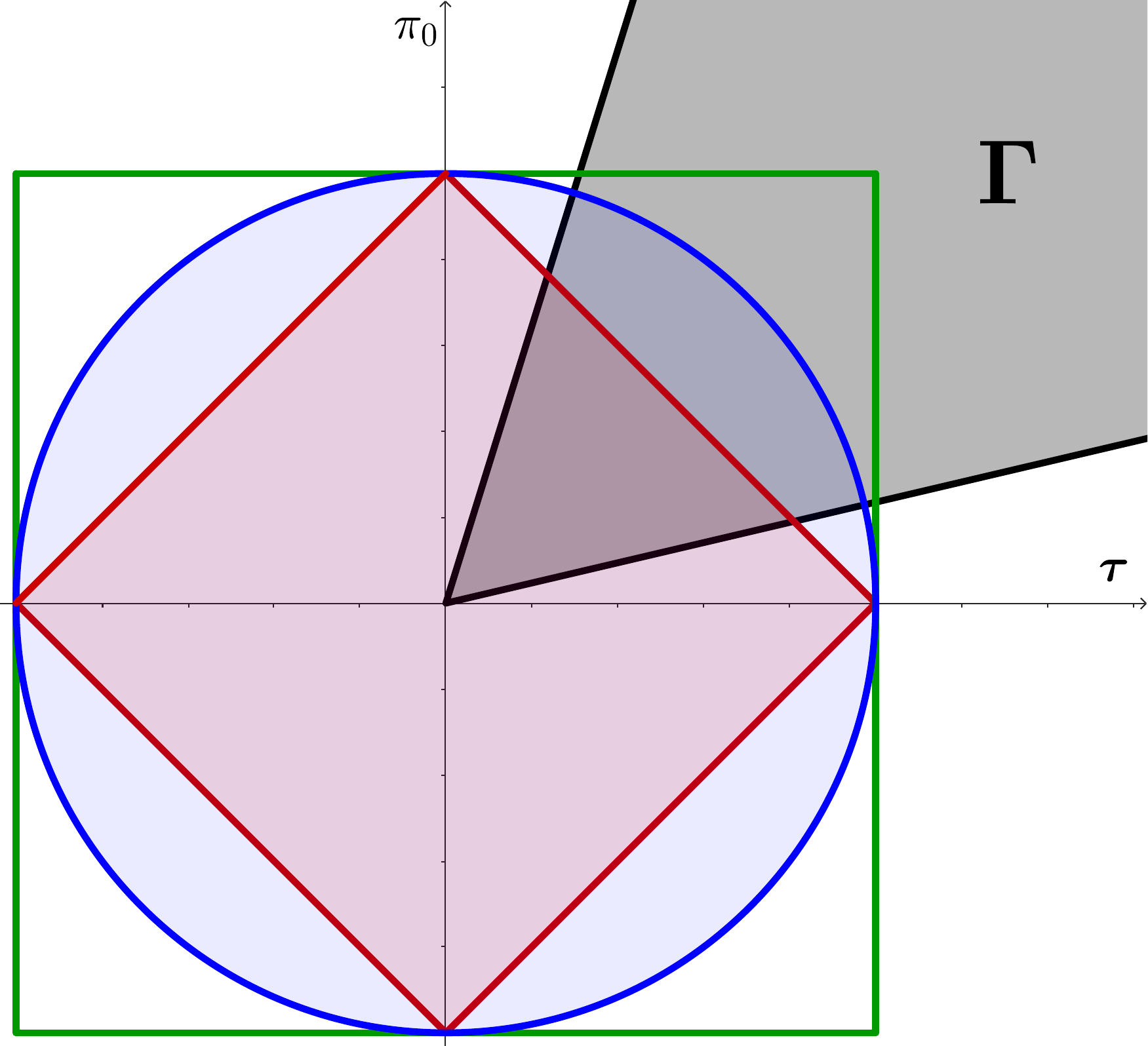}
    \caption{Effect of $\ell_{1}$-norm (red), $\ell_{2}$-norm (blue) and $\ell_{\infty}$-norm (green) on truncating the cone of dual solutions.}\label{fig:norms-truncate}
\end{figure*}

Figure~\ref{fig:guided-projections} provides a schematic view of steps of GPA. In this figure, $(\hat{\yy},\hat{\eta})$ is the point to be separated,
$(\tilde{\yy}^{(1)},\tilde{\eta}^{(1)})$ is the projection of $(\hat{\yy},\hat{\eta})$ onto $\cut({\BFpi}^{(0)},\pi_0^{(0)})$, and 
$(\tilde{\yy}^{(2)},\tilde{\eta}^{(2)})$ is the projection of $(\hat{\yy},\hat{\eta})$ onto $\cut({\BFpi}^{(0)},\pi_0^{(0)})\cap \cut({\BFpi}^{(1)},\pi_0^{(1)})$, which is in $\epi$, thus $\cut({\BFpi}^{(1)},\pi_0^{(1)})$ is the deepest cut in this example.

\begin{figure}[ht]
    \centering
    \includegraphics[width=0.4\textwidth]{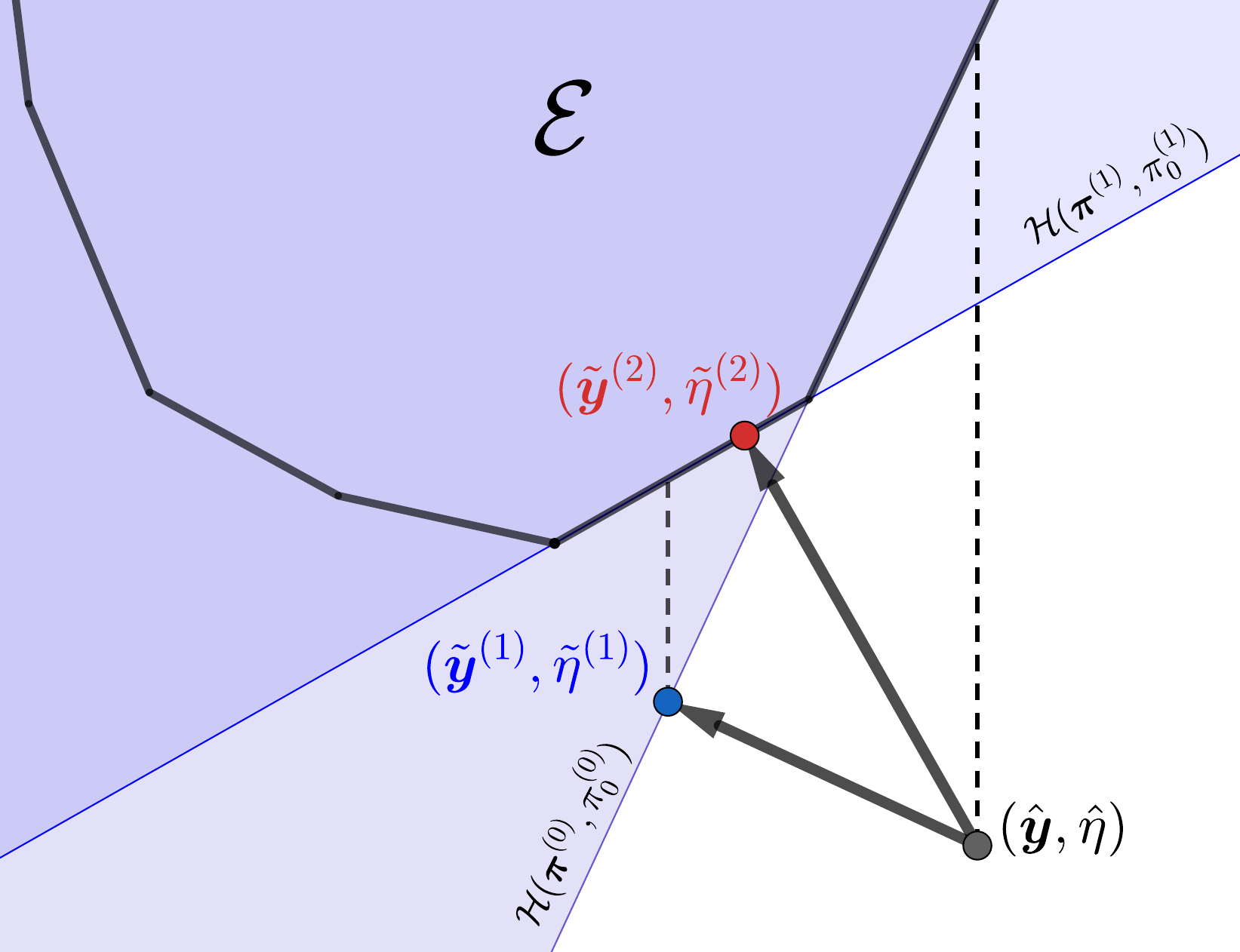}
    \caption{Guided projections algorithm. 
    }
    \label{fig:guided-projections}
\end{figure}

Figure~\ref{fig:gpa_bounds} illustrates iterations of GPA on instances from CFLP and MCNDP.  Lower bounds come from the depth of the cuts, while upper bounds come from distance of the incumbent point to the points identified on the boundary of $\epi$.
Note that for MCNDP, classical Benders cut (first cut) is a feasibility cut, whereas the deepest cut (last cut) is a significantly deeper optimality cut.

\begin{figure}[ht]
    \centering
    \subfigure[CFLP]{
 		\includegraphics[clip,width=0.45\textwidth]{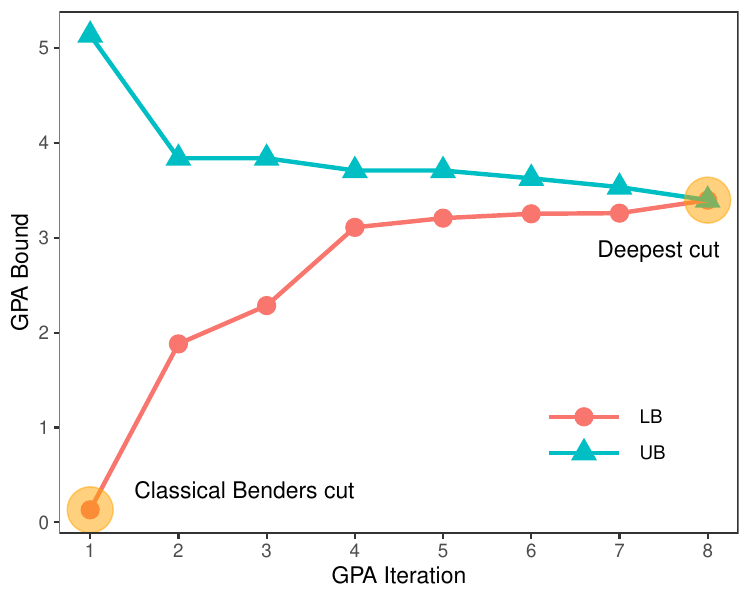}\label{fig:gpa_bounds-scflp}}
    \subfigure[MCNDP]{
 		\includegraphics[clip,width=0.45\textwidth]{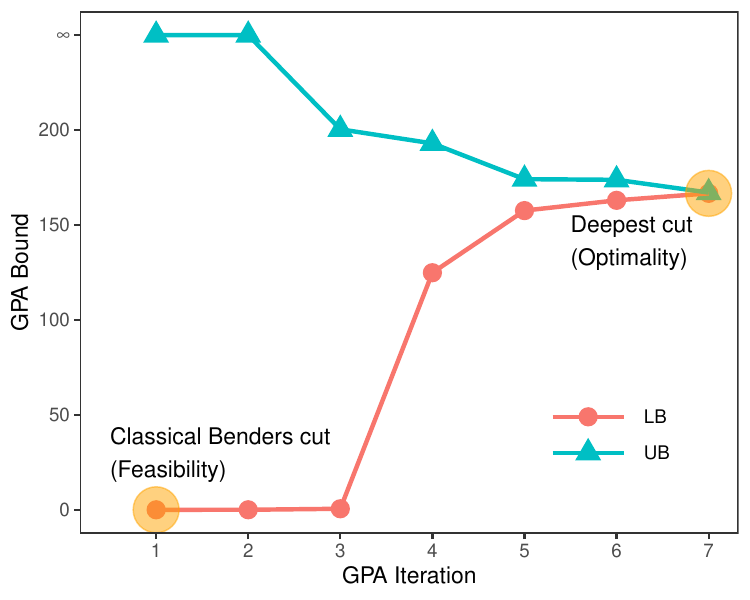}\label{fig:gpa_bounds-ndp}}\hspace{0.2cm}
    \caption{Illustrating GPA iterations on producing an $\ell_1$ deepest cut.}\label{fig:gpa_bounds}
\end{figure}

Figure~\ref{fig:lowerbounds_plot} illustrates the lower bounds produced by $\ell_1$, Conforti-Wolsey (\texttt{CW}) and Classical Benders (\texttt{CB}).

\begin{figure}[ht]
    \centering
    \includegraphics[clip,width=0.55\textwidth]{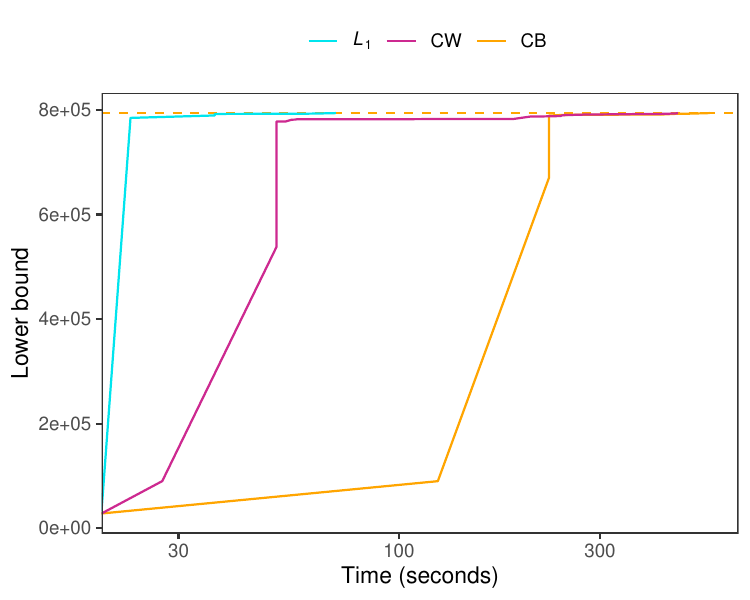}
    \caption{Comparison of lower bounds produced by $\ell_1$-deepest cuts, \texttt{CW} cuts and \texttt{CB} on the sample instance \texttt{cap121} with $|K|=1024$ stochastic demand scenarios. Time axis is logarithmically scaled.}\label{fig:lowerbounds_plot}
\end{figure}

\section{Supplementary Results}\label{app:results}
In this appendix we present supplementary computational results to support the experiments presented in the body of the paper. 

Tables \ref{tab:deterministic_solver_rootgap}, \ref{tab:deterministic_solver_time_solved} and \ref{tab:deterministic_solver_cuts} present detailed results for the first experiment where separation problems are treated as general LP/QPs.
Tables \ref{tab:deterministic_iterative_rootgap}, \ref{tab:deterministic_iterative_time_solved} and \ref{tab:deterministic_iterative_cuts} present detailed results for the second experiment where separation problems are solved using GPA and DDMA. Table \ref{tab:stochastic_cuts} presents number of cuts produced by each model when solving the stochastic instances. Figure~\ref{fig:pp-stochastic-supp} illustrates performance profiles of $\ell_1$, \texttt{CW} and \texttt{CB}.

\begin{table}[t]
\footnotesize
\caption{Deterministic instances: Comparing average percentage gap at the root node from different methods when separation problems are solved as general LP/QP.}
\label{tab:deterministic_solver_rootgap}
\noindent\begin{tabular*}{\columnwidth}{@{\extracolsep{\stretch{1}}}*{10}{@{}lrrrrrrrrr@{}}}
\toprule
Class       & \multicolumn{1}{c}{\texttt{MILP}} & \multicolumn{1}{c}{\texttt{CB}} & \multicolumn{1}{c}{\texttt{MISD}} & \multicolumn{1}{c}{\texttt{CW}} & \multicolumn{1}{c}{\texttt{MWP}} & \multicolumn{1}{c}{$R\ell_1$} & \multicolumn{1}{c}{$\ell_1$} & \multicolumn{1}{c}{$\ell_{\infty}$} & \multicolumn{1}{c}{$\ell_2$} \\ \midrule
CFLP        & & & & & & & & &\\
\;\;\; cap6        & 0.01  & 1.69  & 23.11 & 33.4  & 0.26  & 2.01  & 18.11 & 5.40   & 27.89 \\
\;\;\; cap7        & 0.00     & 29.07 & 1.67  & 26.72 & 26.72 & 4.12  & 23.59 & 37.56 & 1.43  \\
\;\;\; cap9        & 12.63 & 13.16 & 27.97 & 8.36  & 13.53 & 4.10   & 8.65  & 3.98  & 1.23  \\
\;\;\; cap10       & 0.00     & 10.51 & 3.03  & 0.85  & 0.85  & 22.11 & 8.92  & 0.51  & 4.00     \\
\;\;\; cap12       & 14.46 & 9.72  & 44.65 & 25.58 & 25.54 & 5.57  & 0.48  & 9.39  & 7.78  \\
\;\;\; cap13       & 0.00     & 45.99 & 7.60   & 2.31  & 2.31  & 20.81 & 3.94  & 12.42 & 9.59  \\\midrule
\;\;\; Average (\%)    & 4.33  & 17.48 & 17.02 & 15.49 & 10.98 & 9.48  & 10.34 & 10.93 & 8.3   \\\midrule
\;\;\; CST(50;50)   & 0.38  & 1.09  & 0.48  & 0.63  & 0.28  & 0.31  & 0.23  & 0.27  & 0.67  \\
\;\;\; CST(50;100)  & 0.36  & 0.59  & 0.34  & 0.32  & 0.37  & 0.23  & 0.15  & 0.25  & 0.21  \\
\;\;\; CST(100;100) & 0.24  & 0.4   & 0.16  & 0.21  & 0.37  & 0.18  & 0.12  & 0.24  & 0.30   \\
\;\;\; CST(100;200) & 0.54  & 0.59  & 0.54  & 0.29  & 0.66  & 0.43  & 0.29  & 0.38  & 1.52  \\\midrule
\;\;\; Average (\%)    & 0.38  & 0.67  & 0.38  & 0.36  & 0.42  & 0.29  & 0.2   & 0.28  & 0.67  \\\midrule
MCNDP       & & & & & & & & &\\
\;\;\; R1& 0.31  & 0.03  & 0.03  & 0.84  & 0.14  & 0.05  & 0.08  & 0.06  & 8.6   \\
\;\;\; R2& 0.34  & 2.95  & 0.29  & 0.66  & 0.27  & 0.69  & 0.21  & 0.37  & 0.82  \\
\;\;\; R3& 12.35 & 24.14 & 0.26  & 0.02  & 0.00     & 0.25  & 0.03  & 0.29  & 0.77  \\
\;\;\; R4& 1.09  & 6.03  & 1.78  & 9.91  & 2.31  & 2.09  & 1.40   & 2.30   & 8.41  \\
\;\;\; R5& 8.91  & 2.96  & 8.94  & 1.99  & 2.23  & 1.68  & 1.66  & 2.41  & 2.92  \\
\;\;\; R6& 1.82  & 5.92  & 7.82  & 3.09  & 3.81  & 2.87  & 2.44  & 3.06  & 8.99  \\
\;\;\; R7& 7.60   & 44.54 & 3.16  & 3.85  & 2.96  & 3.09  & 1.92  & 3.66  & 8.89  \\
\;\;\; R8& 1.91  & 56.43 & 3.80   & 2.93  & 3.83  & 3.01  & 2.50   & 3.84  & 4.85  \\
\;\;\; R9& 1.37  & 34.67 & 21.06 & 3.07  & 4.53  & 2.33  & 1.99  & 3.01  & 13.36 \\
\;\;\; R10         & 1.89  & 86.53 & 21.18 & 5.21  & 5.99  & 4.01  & 5.36  & 8.56  & 13.42 \\\midrule
\;\;\; Average (\%)     & 3.68  & 23.71 & 6.57  & 3.12  & 2.59  & 2.00     & 1.75  & 2.73  & 7.01  \\ \bottomrule
\end{tabular*}
\end{table}

\begin{table}[t]
\footnotesize
\caption{Deterministic instances: Comparing average computing time of different methods when separation problems are solved as general LP/QP.}
\label{tab:deterministic_solver_time_solved}
\noindent\begin{tabular*}{\columnwidth}{@{\extracolsep{\stretch{1}}}*{10}{@{}lrrrrrrrrr@{}}}
\toprule
Class       & \multicolumn{1}{c}{\texttt{MILP}} & \multicolumn{1}{c}{\texttt{CB}} & \multicolumn{1}{c}{\texttt{MISD}} & \multicolumn{1}{c}{\texttt{CW}} & \multicolumn{1}{c}{\texttt{MWP}} & \multicolumn{1}{c}{$R\ell_1$} & \multicolumn{1}{c}{$\ell_1$} & \multicolumn{1}{c}{$\ell_{\infty}$} & \multicolumn{1}{c}{$\ell_2$} \\ \midrule
CFLP        & & & & & & & & &\\
\;\;\; cap6        & 0.04  & 0.06   & 0.20   & 0.13   & 0.10   & 0.21   & 0.10   & 0.05   & 10.42  \\
\;\;\; cap7        & 0.01  & 0.17   & 0.14   & 0.10   & 0.11   & 0.14   & 0.03   & 0.05   & 18.86  \\
\;\;\; cap9        & 0.05  & 0.12   & 0.62   & 0.29   & 0.30   & 0.57   & 0.20   & 0.11   & 72.28  \\
\;\;\; cap10       & 0.00  & 0.11   & 0.59   & 0.28   & 0.28   & 0.43   & 0.13   & 0.10   & 68.15  \\
\;\;\; cap12       & 0.08  & 0.48   & 2.96   & 1.10   & 0.99   & 4.00   & 0.99   & 0.42   & 445.11 \\
\;\;\; cap13       & 0.01  & 0.32   & 1.73   & 1.15   & 1.10   & 1.76   & 0.58   & 0.39   & 386.23 \\\midrule
\;\;\; Average (s)     & 0.03  & 0.21   & 0.99   & 0.50   & 0.47   & 1.11   & 0.33   & 0.18   & 83.85  \\
\;\;\; Solved (24) & 24    & 24     & 24     & 24     & 24     & 24     & 24     & 24     & 20     \\\midrule
\;\;\; CST(50;50)   & 0.71  & 0.27   & 0.38   & 0.34   & 0.50   & 0.50   & 0.83   & 0.28   & 17.95  \\
\;\;\; CST(50;100)  & 2.35  & 1.73   & 4.30   & 2.87   & 4.57   & 3.70   & 10.60  & 4.29   & 6.63   \\
\;\;\; CST(100;100) & 5.02  & 5.68   & 12.23  & 8.60   & 15.27  & 12.65  & 47.96  & 14.78  & 36.69  \\
\;\;\; CST(100;200) & 25.21 & 48.49  & 279.28 & 181.46 & 280.73 & 240.60 & 44.12  & 37.26  & 159.04 \\\midrule
\;\;\;Average (s)     & 6.26  & 8.23   & 21.26  & 16.24  & 22.57  & 20.06  & 18.93  & 10.37  & 33.77  \\
\;\;\;Solved (64) & 64    & 64     & 56     & 62     & 59     & 58     & 64     & 64     & 61     \\\midrule
MCNDP      & & & & & & & & &\\
\;\;\; R1& 0.08  & 0.16   & 0.11   & 0.09   & 0.08   & 0.07   & 0.13   & 0.11   & 3.83   \\
\;\;\; R2& 0.14  & 0.20   & 0.28   & 0.20   & 0.20   & 0.22   & 0.31   & 0.20   & 10.98  \\
\;\;\; R3& 0.26  & 0.42   & 0.88   & 0.37   & 0.42   & 0.38   & 0.77   & 0.41   & 23.56  \\
\;\;\; R4& 0.35  & 0.69   & 0.58   & 0.55   & 0.72   & 0.44   & 1.05   & 0.91   & 12.85  \\
\;\;\; R5& 1.02  & 2.21   & 4.22   & 1.82   & 2.32   & 1.74   & 4.93   & 2.58   & 63.49  \\
\;\;\; R6& 9.93  & 43.14  & 69.78  & 48.29  & 49.02  & 51.58  & 70.55  & 56.05  & 161.59 \\
\;\;\; R7& 1.13  & 7.69   & 3.73   & 3.74   & 5.03   & 2.89   & 6.56   & 11.04  & 43.29  \\
\;\;\; R8& 5.79  & 33.86  & 28.27  & 20.27  & 25.91  & 16.09  & 34.15  & 31.68  & 172.20 \\
\;\;\; R9& 12.34 & 81.27  & 126.30 & 68.63  & 80.26  & 81.61  & 160.37 & 129.52 & 281.32 \\
\;\;\; R10         & 21.74 & 255.75 & 269.67 & 187.59 & 199.24 & 206.28 & 298.00 & 222.07 & 442.61 \\
\midrule
\;\;\;Average (s)     & 4.04  & 16.00  & 17.76  & 13.45  & 14.71  & 13.61  & 19.91  & 17.81  & 61.70  \\
\;\;\;Solved (81) & 81    & 72     & 71     & 74     & 73     & 74     & 70     & 69     & 59 \\
\bottomrule
\end{tabular*}
\end{table}

\begin{table}[t]
\footnotesize
\caption{Deterministic instances: Comparing number of cuts and fraction of feasibility cuts produced by different methods when separation problems are solved as general LP/QP.}
\label{tab:deterministic_solver_cuts}
\noindent\begin{tabular*}{\columnwidth}{@{\extracolsep{\stretch{1}}}*{9}{@{}lrrrrrrrr@{}}}
\toprule
Class       & \multicolumn{1}{c}{\texttt{CB}} & \multicolumn{1}{c}{\texttt{MISD}} & \multicolumn{1}{c}{\texttt{CW}} & \multicolumn{1}{c}{\texttt{MWP}} & \multicolumn{1}{c}{$R\ell_1$} & \multicolumn{1}{c}{$\ell_1$} & \multicolumn{1}{c}{$\ell_{\infty}$} & \multicolumn{1}{c}{$\ell_2$} \\ \midrule
CFLP        & & & & & & & & \\
\;\;\; ca6        & 32.5    & 32.7    & 16.9   & 18.7   & 13.6   & 6.1    & 7.5    & 16.0   \\
\;\;\; ca7        & 34.4    & 25.9    & 12.9   & 14.0   & 9.5    & 3.7    & 7.5    & 27.7   \\
\;\;\; ca9        & 52.5    & 120.7   & 35.9   & 44.2   & 48.2   & 14.3   & 18.3   & 44.9   \\
\;\;\; ca10       & 46.6    & 53.5    & 18.9   & 20.0   & 27.4   & 6.2    & 10.3   & 57.1   \\
\;\;\; ca12       & 77.2    & 240.5   & 39.1   & 41.4   & 92.5   & 13.2   & 26.6   & 57.6   \\
\;\;\; ca13       & 68.6    & 79.6    & 32.4   & 32.4   & 33.1   & 8.0    & 23.9   & 70.6   \\\midrule
\;\;\; Average (\#)    & 49.8    & 71.2    & 24.6   & 26.7   & 30.5   & 8.2    & 14.6   & 42.0   \\\midrule
\;\;\; CST(50;50)   & 16.6    & 8.3     & 6.3    & 14.1   & 6.3    & 5.4    & 5.3    & 6.2    \\
\;\;\; CST(50;100)  & 26.0    & 16.7    & 11.0   & 18.3   & 9.9    & 7.8    & 8.3    & 9.6    \\
\;\;\; CST(100;100) & 32.9    & 19.9    & 15.8   & 30.5   & 17.7   & 12.7   & 14.5   & 13.4   \\
\;\;\; CST(100;200) & 65.0    & 29.9    & 25.9   & 34.7   & 25.1   & 19.7   & 20.7   & 16.2   \\\midrule
\;\;\; Average (\#)    & 31.9    & 17.7    & 13.7   & 23.3   & 13.7   & 10.7   & 11.4   & 11.0   \\\midrule
MCNDP       & & & & & & & & \\
\;\;\; R1& 130.4   & 47.5    & 27.7   & 30.4   & 26.0   & 23.8   & 30.1   & 45.0   \\
\;\;\; R2& 106.2   & 60.4    & 35.8   & 34.4   & 42.1   & 28.8   & 39.7   & 54.1   \\
\;\;\; R3& 63.6    & 36.2    & 15.3   & 19.6   & 19.2   & 12.7   & 21.9   & 32.4   \\
\;\;\; R4& 172.6   & 84.2    & 41.3   & 48.6   & 40.8   & 40.5   & 48.4   & 65.4   \\
\;\;\; R5& 227.6   & 153.4   & 60.9   & 76.8   & 67.0   & 60.9   & 83.6   & 117.2  \\
\;\;\; R6& 808.4   & 417.9   & 244.0  & 254.4  & 243.1  & 174.9  & 210.0  & 72.1   \\
\;\;\; R7& 393.6   & 173.8   & 96.7   & 122.7  & 118.7  & 95.9   & 129.8  & 147.0  \\
\;\;\; R8& 727.2   & 419.4   & 200.8  & 260.1  & 198.0  & 183.8  & 268.8  & 143.8  \\
\;\;\; R9& 1162.9  & 597.4   & 293.7  & 338.9  & 357.0  & 229.2  & 364.1  & 40.8   \\
\;\;\; R10         & 2292.8  & 897.0   & 599.2  & 680.5  & 648.6  & 295.1  & 471.4  & 93.3   \\\midrule
\;\;\; Average (\#)    & 347.6   & 177.4   & 95.3   & 110.0  & 103.4  & 78.5   & 108.3  & 72.7   \\
\;\;\; Feasibility (\%)    & 28.57\% & 12.69\% & 6.02\% & 4.85\% & 5.41\% & 4.29\% & 3.39\% & 1.66\% \\
\bottomrule
\end{tabular*}
\end{table}

\begin{table}[t]
\footnotesize
\caption{Deterministic instances: Comparing average percentage gap at the root node from different methods when separation problems are solved using GPA and DDMA.}

\label{tab:deterministic_iterative_rootgap}
\noindent\begin{tabular*}{\columnwidth}{@{\extracolsep{\stretch{1}}}*{9}{@{}lrrrrrrrr@{}}}
\toprule
Class       & \multicolumn{1}{c}{\texttt{MILP}}& \multicolumn{1}{c}{\texttt{CB}} & \multicolumn{1}{c}{\texttt{CW}} & \multicolumn{1}{c}{\texttt{MWP}} & \multicolumn{1}{c}{$R\ell_1$} & \multicolumn{1}{c}{$\ell_1$} & \multicolumn{1}{c}{$\ell_{\infty}$} & \multicolumn{1}{c}{$\ell_2$} \\ \midrule
CFLP        & & & & & & & \\
\;\;\; cap6        & 0.01  & 22.68 & 14.64 & 15.37 & 1.17  & 0.05  & 0.16  & 10.18 \\
\;\;\; cap7        & 0.00  & 0.13  & 13.49 & 14.16 & 11.65 & 7.86  & 7.22  & 0.13  \\
\;\;\; cap9        & 12.63 & 30.34 & 33.43 & 35.10 & 4.03  & 2.23  & 4.59  & 39.00 \\
\;\;\; cap10       & 0.00  & 2.06  & 1.16  & 1.22  & 31.55 & 18.65 & 0.67  & 2.19  \\
\;\;\; cap12       & 14.46 & 47.81 & 25.52 & 26.80 & 5.61  & 1.98  & 14.22 & 3.64  \\
\;\;\; cap13       & 0.00  & 7.83  & 2.31  & 2.42  & 16.33 & 20.13 & 11.13 & 6.60  \\\midrule
\;\;\; Average (\%)    & 4.33  & 17.30 & 14.52 & 15.21 & 11.28 & 8.19  & 6.21  & 9.58  \\\midrule
\;\;\; CST(50;50)   & 0.38  & 0.68  & 0.34  & 0.37  & 0.37  & 0.30  & 0.27  & 0.26  \\
\;\;\; CST(50;100)  & 0.36  & 0.35  & 0.20  & 0.22  & 0.15  & 0.44  & 0.25  & 0.44  \\
\;\;\; CST(100;100) & 0.24  & 0.27  & 0.20  & 0.22  & 0.18  & 0.22  & 0.32  & 0.29  \\
\;\;\; CST(100;200) & 0.54  & 0.37  & 0.32  & 0.44  & 0.35  & 0.30  & 0.35  & 0.43  \\\midrule
\;\;\; Average (\%)    & 0.38  & 0.42  & 0.27  & 0.31  & 0.26  & 0.32  & 0.30  & 0.35  \\\midrule
MCNDP       & & & & & & & \\
\;\;\; R1& 0.31  & 0.31  & 0.32  & 0.05  & 0.58  & 0.59  & 0.35  & 7.53  \\
\;\;\; R2& 0.34  & 14.74 & 0.03  & 0.13  & 0.62  & 0.04  & 0.31  & 0.42  \\
\;\;\; R3& 12.35 & 12.33 & 1.71  & 1.71  & 0.01  & 0.31  & 0.22  & 11.94 \\
\;\;\; R4& 1.09  & 24.79 & 2.38  & 3.16  & 1.92  & 2.13  & 2.45  & 5.60  \\
\;\;\; R5& 8.91  & 18.63 & 2.73  & 2.76  & 1.47  & 1.27  & 1.06  & 1.98  \\
\;\;\; R6& 1.82  & 17.87 & 2.89  & 3.33  & 2.66  & 1.96  & 2.85  & 3.44  \\
\;\;\; R7& 7.60  & 18.33 & 3.79  & 3.67  & 3.15  & 9.81  & 2.82  & 3.12  \\
\;\;\; R8& 1.91  & 17.86 & 3.68  & 3.01  & 2.92  & 2.67  & 2.52  & 7.99  \\
\;\;\; R9& 1.37  & 57.55 & 3.83  & 3.76  & 3.07  & 1.95  & 2.22  & 9.83  \\
\;\;\; R10         & 1.89  & 45.48 & 4.91  & 5.45  & 4.06  & 2.75  & 4.40  & 16.32 \\\midrule
\;\;\; Average (\%)    & 3.68  & 21.84 & 2.62  & 2.69  & 2.04  & 2.31  & 1.91  & 6.72  \\
\bottomrule
\end{tabular*}
\end{table}

\begin{table}[t]
\footnotesize
\caption{Deterministic instances: Comparing average computing time and number of instances solved by different methods when separation problems are solved using GPA and DDMA.}
\label{tab:deterministic_iterative_time_solved}
\noindent\begin{tabular*}{\columnwidth}{@{\extracolsep{\stretch{1}}}*{9}{@{}lrrrrrrrr@{}}}
\toprule
Class       & \multicolumn{1}{c}{\texttt{MILP}}& \multicolumn{1}{c}{\texttt{CB}} & \multicolumn{1}{c}{\texttt{CW}} & \multicolumn{1}{c}{\texttt{MWP}} & \multicolumn{1}{c}{$R\ell_1$} & \multicolumn{1}{c}{$\ell_1$} & \multicolumn{1}{c}{$\ell_{\infty}$} & \multicolumn{1}{c}{$\ell_2$} \\ \midrule
CFLP        & & & & & & & \\
\;\;\;  cap6        & 0.04  & 0.08   & 0.03  & 0.04  & 0.08  & 0.04  & 0.04  & 0.27  \\
\;\;\;  cap7        & 0.01  & 0.02   & 0.04  & 0.04  & 0.05  & 0.02  & 0.02  & 0.13  \\
\;\;\;  cap9        & 0.05  & 0.16   & 0.09  & 0.09  & 0.36  & 0.12  & 0.12  & 0.79  \\
\;\;\;  cap10       & 0.00  & 0.09   & 0.07  & 0.07  & 0.21  & 0.08  & 0.06  & 0.49  \\
\;\;\;  cap12       & 0.08  & 0.57   & 0.19  & 0.21  & 1.13  & 0.35  & 0.39  & 1.58  \\
\;\;\;  cap13       & 0.01  & 0.31   & 0.21  & 0.22  & 0.63  & 0.21  & 0.18  & 1.00  \\\midrule
\;\;\; Average (s)     & 0.03  & 0.20   & 0.10  & 0.11  & 0.40  & 0.14  & 0.13  & 0.70  \\
\;\;\; Solved (24)     & 24    & 24     & 24    & 24    & 24    & 24    & 24    & 24    \\\midrule
\;\;\; CST(50;50)   & 0.71  & 0.10   & 0.06  & 0.06  & 0.04  & 0.05  & 0.06  & 0.25  \\
\;\;\; CST(50;100)  & 2.35  & 0.24   & 0.13  & 0.13  & 0.12  & 0.15  & 0.18  & 0.52  \\
\;\;\; CST(100;100) & 5.02  & 0.46   & 0.36  & 0.36  & 0.28  & 0.39  & 0.41  & 0.59  \\
\;\;\; CST(100;200) & 25.21 & 1.74   & 0.97  & 1.03  & 0.87  & 1.09  & 1.14  & 1.63  \\\midrule
\;\;\; Average (s)     & 6.26  & 0.62   & 0.37  & 0.39  & 0.32  & 0.41  & 0.44  & 0.74  \\
\;\;\; Solved (64)    & 64    & 64     & 64    & 64    & 64    & 64    & 64    & 64    \\\midrule
MCNDP       & & & & & & & \\
\;\;\; R1& 0.08  & 0.08   & 0.09  & 0.10  & 0.09  & 0.07  & 0.09  & 0.74  \\
\;\;\; R2& 0.14  & 0.09   & 0.10  & 0.14  & 0.12  & 0.08  & 0.11  & 0.85  \\
\;\;\; R3& 0.26  & 0.11   & 0.12  & 0.16  & 0.15  & 0.07  & 0.09  & 0.57  \\
\;\;\; R4& 0.35  & 0.41   & 0.49  & 0.42  & 0.42  & 0.41  & 0.40  & 2.20  \\
\;\;\; R5& 1.02  & 0.72   & 0.86  & 1.03  & 0.77  & 0.60  & 0.54  & 2.75  \\
\;\;\; R6& 9.93  & 49.23  & 18.24 & 18.94 & 20.03 & 18.77 & 19.58 & 20.97 \\
\;\;\; R7& 1.13  & 4.04   & 4.19  & 3.66  & 4.13  & 3.82  & 2.91  & 6.00  \\
\;\;\; R8& 5.79  & 20.10  & 15.12 & 14.39 & 12.32 & 11.54 & 11.80 & 15.60 \\
\;\;\; R9& 12.34 & 24.73  & 13.83 & 14.78 & 16.54 & 12.56 & 16.26 & 16.41 \\
\;\;\; R10         & 21.74 & 121.39 & 71.16 & 73.08 & 73.14 & 56.08 & 62.82 & 75.89 \\\midrule
\;\;\; Average (s)     & 4.04  & 10.50  & 7.22  & 7.28  & 7.32  & 6.44  & 6.81  & 8.80  \\
\;\;\; Solved (81)     & 81    & 75     & 78    & 77    & 77    & 81    & 79    & 79 \\
\bottomrule
\end{tabular*}
\end{table}

\begin{table}[t]
\footnotesize
\caption{Deterministic instances: Comparing number of cuts and fraction of feasibility cuts produced by different methods when separation problems are solved using GPA and DDMA.}
\label{tab:deterministic_iterative_cuts}
\noindent\begin{tabular*}{\columnwidth}{@{\extracolsep{\stretch{1}}}*{8}{@{}lrrrrrrrr@{}}}
\toprule
Class       & \multicolumn{1}{c}{\texttt{CB}} & \multicolumn{1}{c}{\texttt{CW}} & \multicolumn{1}{c}{\texttt{MWP}} & \multicolumn{1}{c}{$R\ell_1$} & \multicolumn{1}{c}{$\ell_1$} & \multicolumn{1}{c}{$\ell_{\infty}$} & \multicolumn{1}{c}{$\ell_2$} \\ \midrule
CFLP        & & & & & & & \\
\;\;\; cap6        & 40.2    & 17.8   & 18.7   & 15.8   & 8.2    & 11.5   & 13.8   \\
\;\;\; cap7        & 17.5    & 12.9   & 13.5   & 10.0   & 6.9    & 7.5    & 9.5    \\
\;\;\; cap9        & 84.6    & 30.6   & 32.1   & 57.6   & 22.7   & 32.9   & 51.1   \\
\;\;\; cap10       & 46.9    & 23.8   & 25.0   & 28.6   & 17.4   & 16.5   & 27.3   \\
\;\;\; cap12       & 154.7   & 38.6   & 40.5   & 84.5   & 29.1   & 58.7   & 81.7   \\
\;\;\; cap13       & 94.0    & 33.6   & 35.3   & 38.7   & 26.2   & 28.5   & 54.1   \\\midrule
\;\;\; Average (\#)     & 60.9    & 25.0   & 26.3   & 32.9   & 17.0   & 22.3   & 32.9 \\\midrule
\;\;\; CST(50;50)   & 17.3    & 8.9    & 8.6    & 5.8    & 6.6    & 8.4    & 8.9    \\
\;\;\; CST(50;100)  & 21.0    & 11.1   & 10.6   & 9.6    & 10.1   & 12.2   & 13.7   \\
\;\;\; CST(100;100) & 31.5    & 19.0   & 18.5   & 16.3   & 18.8   & 21.2   & 15.8   \\
\;\;\; CST(100;200) & 48.0    & 27.0   & 26.4   & 25.9   & 26.2   & 27.7   & 26.1   \\\midrule
\;\;\; Average (\#)     & 27.8    & 15.6   & 15.1   & 13.2   & 14.3   & 16.3   & 15.4   \\\midrule
MCNDP       & & & & & & & \\
\;\;\; R1& 88.5    & 30.6   & 33.3   & 28.3   & 21.4   & 33.5   & 23.5   \\
\;\;\; R2& 95.8    & 32.4   & 33.9   & 38.3   & 28.2   & 33.3   & 27.1   \\
\;\;\; R3& 50.4    & 15.1   & 17.7   & 21.2   & 11.4   & 15.7   & 15.7   \\
\;\;\; R4& 118.6   & 40.6   & 41.4   & 41.9   & 35.3   & 39.3   & 36.4   \\
\;\;\; R5& 177.2   & 54.5   & 64.0   & 61.6   & 47.7   & 50.3   & 51.1   \\
\;\;\; R6& 551.6   & 243.9  & 258.5  & 278.5  & 183.5  & 179.5  & 169.4  \\
\;\;\; R7& 289.6   & 107.1  & 104.0  & 107.7  & 93.1   & 97.3   & 81.9   \\
\;\;\; R8& 494.8   & 215.0  & 197.9  & 205.4  & 141.5  & 168.0  & 141.2  \\
\;\;\; R9& 850.4   & 329.1  & 296.5  & 384.0  & 239.9  & 281.4  & 225.6  \\
\;\;\; R10         & 1972.4  & 730.6  & 778.2  & 784.7  & 479.8  & 536.2  & 479.8  \\\midrule
\;\;\; Average (\#)    & 261.7   & 99.0   & 101.9  & 107.2  & 77.0   & 88.1   & 77.4   \\
\;\;\; Feasibility (\%)     & 22.78\% & 4.04\% & 5.06\% & 3.64\% & 3.11\% & 1.53\% & 2.03\% \\
\bottomrule
\end{tabular*}
\end{table}

\begin{table}[t]
\renewcommand{\arraystretch}{1.2}
\footnotesize
\caption{Stochastic instances: comparing average percentage gap and computing time at the root node.}
\label{tab:stochastic_rgaptime}
\noindent\begin{tabular*}{\columnwidth}{@{\extracolsep{\stretch{1}}}*{16}{@{}lrrrrrrrlrrrrrrr@{}}}
\toprule
\multicolumn{1}{c}{\multirow{2}{*}{Problem}} & \multicolumn{7}{c}{Root node gap (\%)}       &  & \multicolumn{7}{c}{Root node computing time (s)}       \\ \cline{2-8} \cline{10-16} 
    & \multicolumn{1}{c}{\texttt{MILP}} & \multicolumn{1}{c}{\texttt{CB}} & \multicolumn{1}{c}{\texttt{CW}} & \multicolumn{1}{c}{\texttt{MWP}} & \multicolumn{1}{c}{$R\ell_1$} & \multicolumn{1}{c}{$\ell_1$} & \multicolumn{1}{c}{$\ell_{\infty}$} &  & \multicolumn{1}{c}{\texttt{MILP}} & \multicolumn{1}{c}{\texttt{CB}} & \multicolumn{1}{c}{\texttt{CW}} & \multicolumn{1}{c}{\texttt{MWP}} & \multicolumn{1}{c}{$R\ell_1$} & \multicolumn{1}{c}{$\ell_1$} & \multicolumn{1}{c}{$\ell_{\infty}$}\\ \cline{1-8} \cline{10-16}
\multicolumn{1}{c}{CFLP (\texttt{CAP})} &      &    &    &    &     &    &      &  &      &    &    &    &     &    &      \\ \cline{1-1}
\;\;\; $|K|=256$    & 8.39  & 8.45& 11.74        & 12.17        & 5.49 & 4.77& 6.12  &  & 233.4  & 5.7 & 4.6 & 5.1 & 6.3 & 5.0   & 3.7   \\
\;\;\; $|K|=512$    & 4.47  & 16.69        & 12.00        & 12.45        & 4.07 & 3.03& 7.90  &  & 1046.6 & 5.8 & 4.2 & 3.9 & 5.5 & 6.4 & 6.4   \\
\;\;\; $|K|=1024$  & 54.24 & 19.25        & 7.05& 7.31& 3.26 & 1.46& 8.00  &  & 3080.8 & 4.5 & 3.6 & 4.0   & 5.7 & 4.7 & 5.0   \\
\;\;\; Average      & 20.43 & 14.70        & 10.24        & 10.62        & 4.27 & 3.08& 7.34  &  & 916.4  & 5.3 & 4.1 & 4.4 & 5.8 & 5.4 & 5.0   \\\midrule
\multicolumn{1}{c}{MCNDP (\texttt{R})} &      &    &    &    &     &    &      &  &      &    &    &    &     &    &      \\ \cline{1-1}
\;\;\; $|K|=16$    & 8.40  & 19.32        & 6.27& 7.66& 9.05 & 3.16& 2.95  &  & 25.6   & 1.6 & 2.1 & 1.7 & 1.9 & 1.8 & 1.8   \\
\;\;\; $|K|=32$    & 10.88 & 19.00        & 4.93& 4.95& 9.25 & 3.57& 3.15  &  & 62.5   & 2.1 & 2.4 & 3.0   & 2.4 & 2.1 & 2.0   \\
\;\;\; $|K|=64$    & 10.90 & 14.66        & 4.92& 5.02& 5.24 & 3.57& 3.64  &  & 164.1  & 3.3 & 3.5 & 3.4 & 3.7 & 3.2 & 2.9   \\
\;\;\; Average      & 10.06 & 17.64        & 5.37& 5.87& 7.83 & 3.43& 3.25  &  & 66.6   & 2.3 & 2.6 & 2.7 & 2.7 & 2.4 & 2.3   \\\midrule
\multicolumn{1}{c}{SNIP} &      &    &    &    &     &    &      &  &      &    &    &    &     &    &      \\ \cline{1-1}
\;\;\; \texttt{budget = 30} & 33.77 & 23.83        & 24.31        & 28.52        & 20.85 & 20.90        & 22.44 &  & 2367.6 & 4.0   & 4.4 & 3.3 & 5.4 & 6.1 & 4.5   \\
\;\;\; \texttt{budget = 40}& 36.43 & 32.35        & 27.85        & 27.44        & 27.57 & 23.29        & 29.15 &  & 2775.6 & 5.9 & 4.6 & 4.7 & 6.5 & 4.6 & 4.5   \\
\;\;\; \texttt{budget = 50} & 38.29 & 29.96        & 30.10        & 29.53        & 25.16 & 25.28        & 26.98 &  & 3009.8 & 4.1 & 3.5 & 3.8 & 5.3 & 4.0   & 4.8   \\
\;\;\; Average      & 36.15 & 28.66        & 27.40        & 28.50        & 24.50 & 23.14        & 26.16 &  & 2704.4 & 4.6 & 4.2 & 3.9 & 5.7 & 4.9 & 4.6  \\ \bottomrule
\end{tabular*}
\end{table}

\begin{table}[t]
\renewcommand{\arraystretch}{1.2}
\footnotesize
\caption{Stochastic instances: comparing average computing time (in seconds) and B\&B nodes.}
\label{tab:stochastic_timenodes}
\noindent\begin{tabular*}{\columnwidth}{@{\extracolsep{\stretch{1}}}*{16}{@{}lrrrrrrrlrrrrrrr@{}}}
\toprule
\multicolumn{1}{c}{\multirow{2}{*}{Problem}} & \multicolumn{7}{c}{Computing time (s)}       &  & \multicolumn{7}{c}{B\&B nodes}       \\ \cline{2-8} \cline{10-16} 
    & \multicolumn{1}{c}{\texttt{MILP}} & \multicolumn{1}{c}{\texttt{CB}} & \multicolumn{1}{c}{\texttt{CW}} & \multicolumn{1}{c}{\texttt{MWP}} & \multicolumn{1}{c}{$R\ell_1$} & \multicolumn{1}{c}{$\ell_1$} & \multicolumn{1}{c}{$\ell_{\infty}$} &  & \multicolumn{1}{c}{\texttt{MILP}} & \multicolumn{1}{c}{\texttt{CB}} & \multicolumn{1}{c}{\texttt{CW}} & \multicolumn{1}{c}{\texttt{MWP}} & \multicolumn{1}{c}{$R\ell_1$} & \multicolumn{1}{c}{$\ell_1$} & \multicolumn{1}{c}{$\ell_{\infty}$}\\ \cline{1-8} \cline{10-16}
\multicolumn{1}{c}{CFLP (\texttt{CAP})} &      &    &    &    &     &    &      &  &      &    &    &    &     &    &      \\ \cline{1-1}

\;\;\; $|K|=256$    & 235  & 57   & 28   & 28   & 48   & 10   & 23   &  & ---  & 488    & 257    & 243     & 662    & 38     & 231    \\
\;\;\; $|K|=512$    & 1048 & 135  & 56   & 58   & 111  & 19   & 52   &  & ---  & 463    & 192    & 192     & 511    & 36     & 200    \\
\;\;\; $|K|=1024$  & 3081 & 366  & 160  & 168  & 293  & 42   & 126  &  & ---  & 456    & 222    & 222     & 425    & 42     & 170    \\
\;\;\; Average      & 919  & 144  & 65   & 67   & 118  & 21   & 55   &  & ---  & 469    & 222    & 218     & 524    & 39     & 199    \\
\;\;\; Solved (48)     & 37   & 48   & 48   & 48   & 46   & 48   & 48   &  &      &        &        &         &        &        &        \\\midrule
\multicolumn{1}{c}{MCNDP (\texttt{R})} &      &    &    &    &     &    &      &  &      &    &    &    &     &    &      \\ \cline{1-1}
\;\;\; $|K|=16$    & 56   & 19   & 15   & 16   & 16   & 12   & 14   &  & 31   & 248    & 101    & 111     & 103    & 80     & 109    \\
\;\;\; $|K|=32$    & 111  & 26   & 23   & 25   & 21   & 19   & 19   &  & 19   & 290    & 123    & 118     & 114    & 85     & 83     \\
\;\;\; $|K|=64$    & 237  & 38   & 35   & 38   & 34   & 31   & 33   &  & 10   & 274    & 148    & 128     & 143    & 122    & 121    \\
\;\;\; Average      & 115  & 27   & 23   & 25   & 23   & 20   & 21   &  & 19   & 270    & 123    & 119     & 119    & 94     & 103    \\
\;\;\; Solved (133)      & 111  & 126  & 129  & 128  & 131  & 131  & 130  &  &      &        &        &         &        &        &        \\\midrule
\multicolumn{1}{c}{SNIP$^{\dagger}$} &      &    &    &    &     &    &      &  &      &    &    &    &     &    &      \\ \cline{1-1}
\;\;\; \texttt{budget = 30} & 3601 & 434  & 412  & 365  & 366  & 381  & 397  &  & ---  & 188 & 220 & 193  & 169 & 175 & 181 \\
\;\;\; \texttt{budget = 40}& 3603 & 2632 & 2397 & 2702 & 2607 & 2248 & 2474 &  & ---  & 872 & 909 & 1045 & 728 & 731 & 820 \\
\;\;\; \texttt{budget = 50} & 3606 & 3127 & 3216 & 2973 & 2773 & 2309 & 2892 &  & ---  & 644 & 527 & 534  & 529 & 522 & 535 \\
\;\;\; Average & 3603 & 1534 & 1475 & 1438 & 1389 & 1260 & 1421 &  & ---  & 473 & 472 & 476  & 405 & 406 & 430 \\
\;\;\; Solved  (30)   & 0    & 15   & 15   & 16   & 18   & 23   & 17   &  &      &        &        &         &        &        & \\\bottomrule
\multicolumn{10}{l}{\small $^{\dagger}$ B\&B nodes are in thousands for SNIP.} \\
\end{tabular*}
\end{table}

\begin{table}[t]
\renewcommand{\arraystretch}{1.2}
\footnotesize
\caption{Stochastic instances: comparing average number of cuts.}
\label{tab:stochastic_cuts}
\noindent\begin{tabular*}{\columnwidth}{@{\extracolsep{\stretch{1}}}*{16}{@{}lrrrrrrrlrrrrrrr@{}}}
\toprule
\multicolumn{1}{c}{\multirow{2}{*}{Problem}} & \multicolumn{6}{c}{Root node}       &  & \multicolumn{6}{c}{Total}       \\ \cline{2-7} \cline{9-14} 
     & \multicolumn{1}{c}{\texttt{CB}} & \multicolumn{1}{c}{\texttt{CW}} & \multicolumn{1}{c}{\texttt{MWP}} & \multicolumn{1}{c}{$R\ell_1$} & \multicolumn{1}{c}{$\ell_1$} & \multicolumn{1}{c}{$\ell_{\infty}$} & & \multicolumn{1}{c}{\texttt{CB}} & \multicolumn{1}{c}{\texttt{CW}} & \multicolumn{1}{c}{\texttt{MWP}} & \multicolumn{1}{c}{$R\ell_1$} & \multicolumn{1}{c}{$\ell_1$} & \multicolumn{1}{c}{$\ell_{\infty}$}\\ \cline{1-7} \cline{9-14}
\multicolumn{1}{c}{CFLP (\texttt{CAP})} &      &    &    &    &     &    &      &  &     &    &     &    &      \\ \cline{1-1}
\;\;\; $|K|=256$   & 5314  & 3581  & 3905  & 2039 & 2385 & 2579 &  & 22562 & 8246  & 9015  & 13039 & 2965  & 6081 \\
\;\;\; $|K|=512$   & 9770  & 7971  & 8622  & 4194 & 4248 & 4871 &  & 45234 & 16175 & 17370 & 22062 & 5400  & 12461 \\
\;\;\; $|K|=1024$  & 19695 & 14765 & 16061 & 7953 & 9187 & 9295 &  & 88101 & 30246 & 33258 & 40823 & 11468 & 21911 \\
\;\;\; Average     & 10076 & 7499  & 8149  & 4083 & 4534 & 4889 &  & 44801 & 15920 & 17335 & 22731 & 5685  & 11843 \\\midrule
\multicolumn{1}{c}{MCNDP (\texttt{R})} &      &    &    &    &     &    &    &    &    &    &     &    &      \\ \cline{1-1}
\;\;\; $|K|=16$    & 341   & 294   & 290   & 323  & 249  & 274  &  & 733   & 520   & 523   & 498   & 353   & 414   \\
\;\;\; $|K|=32$    & 783   & 608   & 617   & 610  & 482  & 494  &  & 1500  & 1149  & 1197  & 1008  & 716   & 713   \\
\;\;\; $|K|=64$    & 1877  & 1234  & 1192  & 1150 & 878  & 993  &  & 3046  & 2144  & 2150  & 1960  & 1430  & 1519   \\
\;\;\; Average     & 797   & 606   & 599   & 611  & 474  & 514  &  & 1497  & 1088  & 1106  & 996   & 714   & 767   \\\midrule
\multicolumn{1}{c}{SNIP} &      &    &    &    &     &  &      &    &    &    &     &    &      \\ \cline{1-1}
\;\;\; \texttt{budget = 30} & 11816 & 4624  & 6908  & 3255 & 3134 & 5202 &  & 27413 & 5642  & 10968 & 5501  & 5621  & 10559 \\
\;\;\; \texttt{budget = 40} & 10948 & 5617  & 9028  & 3727 & 3638 & 5784 &  & 29556 & 7133  & 13407 & 7632 & 7218  & 12778 \\
\;\;\; \texttt{budget = 50} & 13386 & 6406  & 9835  & 3631 & 3592 & 6472 &  & 29216 & 8445  & 16243 & 7883 & 7524  & 13346 \\
\;\;\; Average     & 12009 & 5500  & 8497  & 3531 & 3447 & 5796 &  & 28712 & 6979  & 13367 & 6917 & 6733  & 12166 \\ \bottomrule
\end{tabular*}
\end{table}

\begin{figure*}[ht]
    \centering
\begin{adjustbox}{minipage=\linewidth,scale=0.93}
    \subfigure[CFLP: B\&B nodes.]{
 		\includegraphics[clip,width=0.29\textwidth]{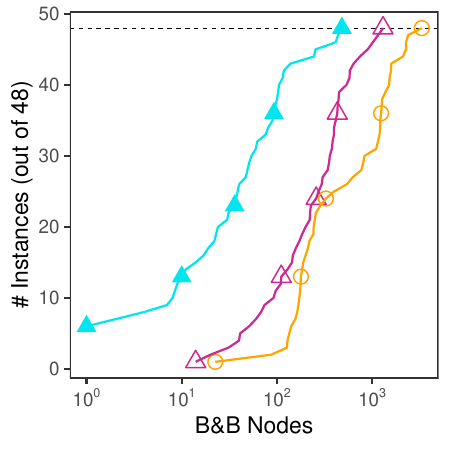}\label{fig:pp-scflp-nodes}}
    \subfigure[CFLP: Root node cuts.]{
 		\includegraphics[clip,width=0.29\textwidth]{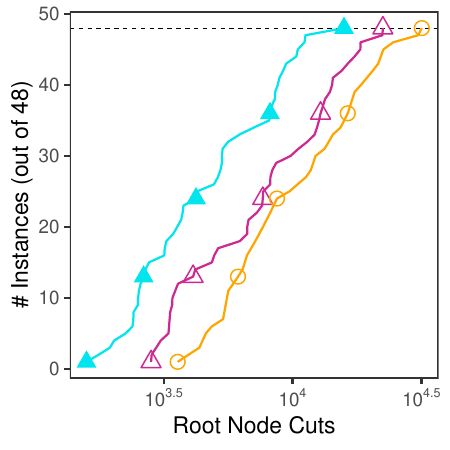}\label{fig:pp-scflp-rncut}}
    \subfigure[CFLP: Root node time.]{
 	    \includegraphics[clip,width=0.365\textwidth]{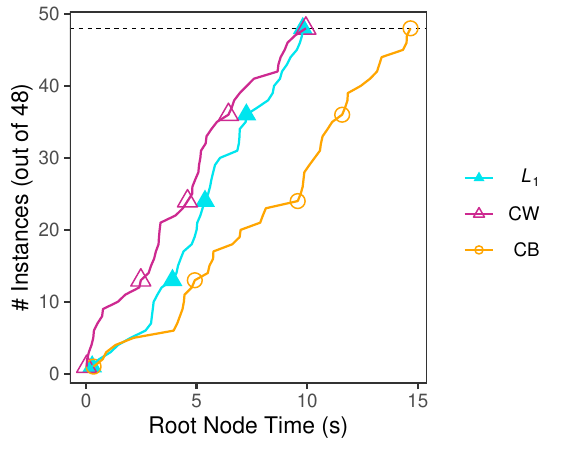}\label{fig:pp-scflp-rntime}}
    \subfigure[MCNDP: B\&B nodes.]{
 		\includegraphics[clip,width=0.29\textwidth]{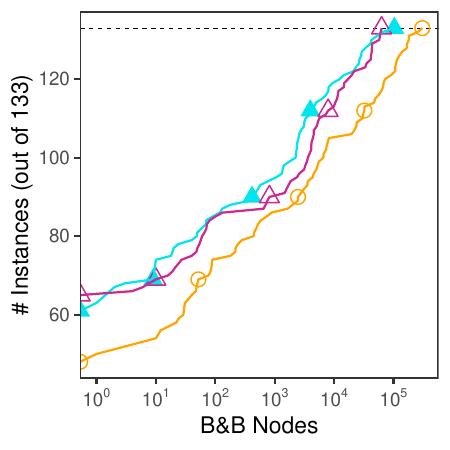}\label{fig:pp-sndp-nodes}}
    \subfigure[MCNDP: Root node cuts.]{
 		\includegraphics[clip,width=0.29\textwidth]{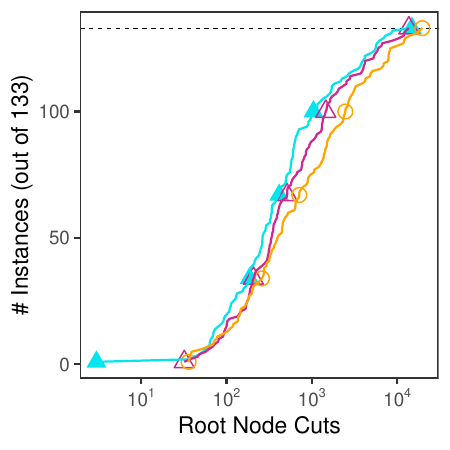}\label{fig:pp-sndp-rncut}}
    \subfigure[MCNDP: Root node time.]{
 	    \includegraphics[clip,width=0.365\textwidth]{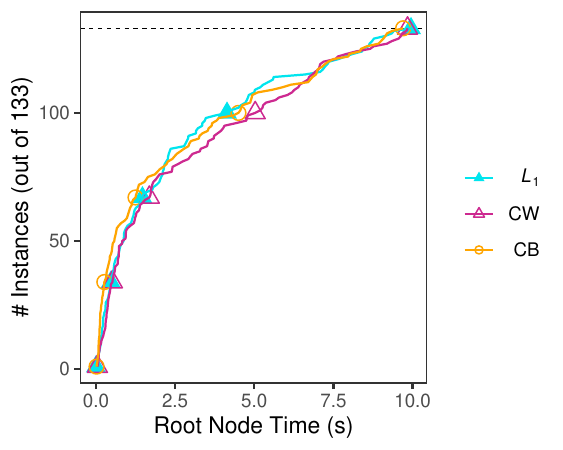}\label{fig:pp-sndp-rntime}}
    \subfigure[SNIP: B\&B nodes.]{
 		\includegraphics[clip,width=0.29\textwidth]{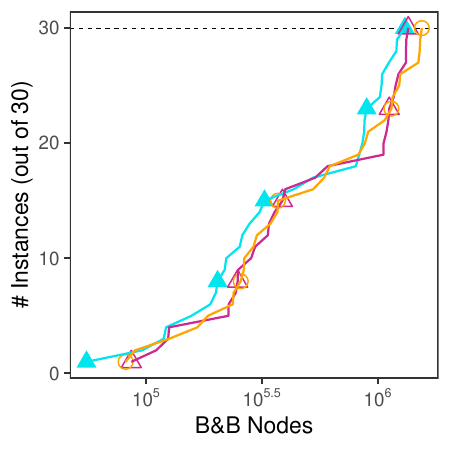}\label{fig:pp-snip-nodes}}
    \subfigure[SNIP: Root node cuts.]{
 		\includegraphics[clip,width=0.29\textwidth]{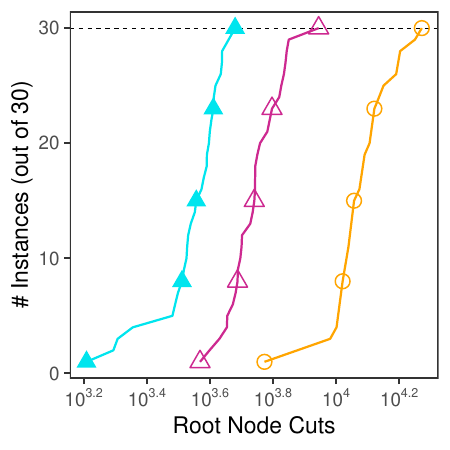}\label{fig:pp-snip-rncut}}
    \subfigure[SNIP: Root node time.]{
 	    \includegraphics[clip,width=0.365\textwidth]{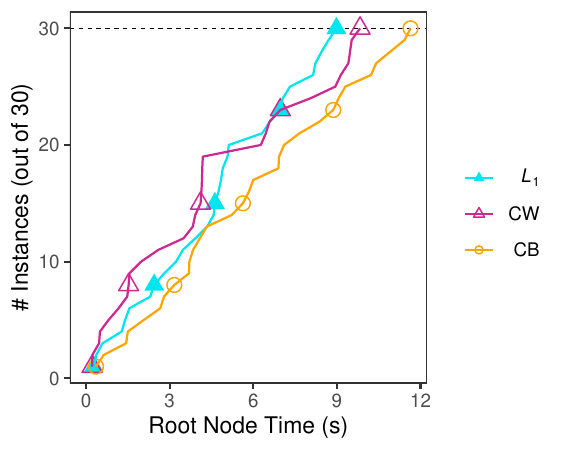}\label{fig:pp-snip-rntime}}
\end{adjustbox}
    \caption{Supplementary results for $\ell_1$-deepest cuts, CW cuts and classical Benders cuts on instances of stochastic CFLP (first row) and stochastic MCNDP (second row) and SNIP (third row).}\label{fig:pp-stochastic-supp}
\end{figure*}

\end{APPENDICES}

\end{document}